\documentclass[aop,MSNbibl,seceqn,nameyear,dvips]{arximspdf}

\doi{10.1214/13-AOP870} %kopijuoti is PTS
\volume{42}
\issue{5}
\pubyear{2014}
\firstpage{2032}
\lastpage{2112}

\makeatletter

\newcommand{\rrvert}{\vert}
\newcommand{\llvert}{\vert}
\def\cal{\mathcal}
\newtheorem{theorem}{Theorem}[section]
\newtheorem{lemma}[theorem]{Lemma}
\newtheorem{cor}[theorem]{Corollary}
\newproclaim{definition}[theorem]{Definition}
\newtheorem{prop}[theorem]{Proposition}
\newtheorem{claim}[theorem]{Claim}
\newproclaim{remark}[theorem]{Remark}

\newproclaim{example}{Example}

\newcommand{\R}{{\mathbb{R}}}
\newcommand{\Q}{{\mathbb{Q}}}
\newcommand{\N}{{\mathbb{N}}}
\newcommand{\Z}{{\mathbb{Z}}}
\newcommand{\CG}{\mathcal{G}}
\newcommand{\CO}{\mathcal{O}}
\newcommand{\CF}{\mathcal{F}}
\newcommand{\CH}{\mathcal{H}}
\def\NN{\mathbb{N}}
\def\1{\mathbf{1}}
\newcommand{\tU}{\widetilde{U}}
\newcommand{\tV}{\widetilde{V}}
\newcommand{\tW}{\widetilde{W}}
\newcommand{\tM}{\widetilde{M}}
\newcommand{\bU}{\bar{U}}
\newcommand{\bV}{\bar{V}}
\newcommand{\bW}{\bar{W}}
\newcommand{\bM}{\bar{M}}
\newcommand{\btau}{\bar{\tau}}
\newcommand{\tgamma}{\tilde{\gamma}}
\newcommand{\vep}{\varepsilon}
\newcommand{\ep}{\varepsilon}
\newcommand{\eqref}[1]{(\ref{#1})}
\makeatother

\begin{document}
\begin{frontmatter}

\title{Nonuniqueness for a parabolic SPDE with $\frac{3}{4}-\vep$-H\"older diffusion coefficients}
\runtitle{Nonuniqueness for a parabolic SPDE}

\begin{aug}
\author[a]{\fnms{Carl} \snm{Mueller}\corref{}\thanksref{t1}\ead[label=e1]{carl.2013@outlook.com}},
\author[b]{\fnms{Leonid} \snm{Mytnik}\thanksref{t2}\ead[label=e2]{leonid@ie.technion.ac.il}}
\and
\author[c]{\fnms{Edwin} \snm{Perkins}\thanksref{t3}\ead[label=e3]{perkins@math.ubc.ca}}

\thankstext{t1}{Supported by an NSF grant.}
\thankstext{t2}{Supported by an ISF grant.}
\thankstext{t3}{Supported by an NSERC Discovery grant.}
\runauthor{C. Mueller, L. Mytnik and E. Perkins}

\affiliation{University of Rochester,
Technion and University of British Columbia}

\address[a]{C. Mueller\\
Department of Mathematics\\
University of Rochester\\
Rochester, New York 14627\\
USA\\
\printead{e1}}

\address[b]{L. Mytnik\\
Faculty of Industrial Engineering\\
\quad and Management\\
Technion---Israel Institute of Technology\\
Haifa 32000\\
Israel\\
\printead{e2}}

\address[c]{E. Perkins\\
Department of Mathematics\\
University of British Columbia\\
Vancouver, British Columbia V6T 1Z2\\
Canada\\
\printead{e3}}
\end{aug}

% HISTORY:
\received{\smonth{8} \syear{2012}}
\revised{\smonth{4} \syear{2013}}

% ABSTRACT
\medskip
\begin{abstract}
Motivated by Girsanov's nonuniqueness examples for SDEs, we
prove nonuniqueness for the parabolic stochastic partial differential equation
(SPDE)
\[
\frac{\partial u}{\partial t}=\frac{\Delta}{2} u(t,x) + \bigl|u(t,x)\bigr|^\gamma
\dot{W}(t,x),\qquad  u(0,x) = 0.
\]
Here $\dot W$ is a space--time white noise on $\R_+\times\R$. More precisely,
we show the above stochastic PDE has a nonzero solution for $0<\gamma<3/4$.
Since \mbox{$u(t,x)=0$} solves the equation, it follows that
solutions are neither unique in law nor pathwise unique. An analogue of
Yamada--Watanabe's famous theorem for SDEs was recently shown in
Mytnik and Perkins [\textit{Probab. Theory Related Fields} \textbf{149}
(2011) 1--96]
for SPDE's by establishing pathwise uniqueness of solutions to
\[
\frac{\partial u}{\partial t}=\frac{\Delta}{2} u(t,x) + \sigma \bigl(u(t,x)\bigr)
\dot{W}(t,x)
\]
if $\sigma$ is H\"older continuous of index $\gamma>3/4$. Hence our examples
show this result is essentially sharp. The situation for the above
class of
parabolic SPDE's is therefore similar to their finite dimensional
counterparts, but with the index $3/4$ in place of $1/2$. The case
$\gamma=1/2$ of the first equation above is particularly interesting
as it
arises as the scaling limit of the signed mass for a system of annihilating
critical branching random walks.
\end{abstract}

% KEYWORDS

\end{frontmatter}
\newpage
%s1 #&#
\section{Introduction}
\label{intro}
\setcounter{equation}{0}

This work concerns uniqueness theory for parabolic semilinear
stochastic partial differential equations (SPDE)
of the form
%
%e1.1 #&#
\begin{eqnarray}
\label{generalized} \frac{\partial u}{\partial t} (t,x)&=& \frac{\Delta}{2} u(t,x) + \sigma
\bigl(x,u(t,x)\bigr) \dot{W}(t,x),
\nonumber
\\[-8pt]
\\[-8pt]
\nonumber
u(0,x) &=& u_0(x),
\nonumber
\end{eqnarray}
where $\dot{W}(t,x)$ is two-parameter white noise on $\R_+\times\R$,
and $\sigma\dvtx \R^2\to\R$ is $\gamma$-H\"older continuous in $u$ and also
has at most linear growth at $\infty$ in $u$. See (2.1)$'$ in \citet
{shi94} or \eqref{SPDEMP} below
for a precise definition of a solution.
Weak existence of solutions in the appropriate function space is then
standard; see, for example, Theorems 1.1 and 2.6 of \citet{shi94} or
Theorem~1.1 of
\citet{mp11}. If $\gamma=1$, then $\sigma$ is Lipschitz in $u$,
and pathwise uniqueness of solutions follows from standard fixed-point
arguments; see Chapter~3 in \citet{wal86}.
A~natural question is then:
\[
\mbox{If $\gamma<1$, are solutions pathwise unique?}
\]
The motivation for this problem comes from a number of models arising from
branching processes and population genetics for which $\gamma=1/2$.

Next we give some examples. In the first three, we only consider nonnegative
solutions, while in the fourth example we allow solutions to take negative
values. If $E\subset\R$, we write $C(E)$ for the space of continuous
functions on $E$ with the topology of uniform convergence on compact sets.

\begin{example}\label{ex1}
If $\sigma(u)=\sqrt u$ and we assume $u\ge
0$, then
a solution to \eqref{generalized} corresponds to the density
$u(t,x)\,dx=X_t(dx)$, where
$X_t$ is the one-dimensional super-Brownian motion. The super-Brownian motion
is a meas\-ure-valued process which
arises as the rescaled limit of branching random walks; see \citet
{rei89} and
\citet{ks88}. More precisely, assume that particles occupy sites
in $\Z/\sqrt N$. With Poisson rate $N/2$, each particle produces
offspring at
a randomly chosen nearest neighbor site. Finally, particles die at rate $N/2$.
For $x\in\Z/\sqrt N$ and $t\geq0$, set
\[
U^N(t,x)=N^{-1/2}\times(\mbox{number of particles at $x$ at
time $t$}).
\]
If the initial ``densities'' converge in the appropriate state space,
then $U^N$
will converge weakly on the appropriate function space to the
solution of \eqref{generalized}, with $\sigma$ as above; see
\citet{rei89} for a proof of this result using nonstandard analysis.
Furthermore, this solution is unique in law.
Uniqueness in law is established by the well-known exponential duality between
$u(t,x)$ and solutions $v(t,x)$ of the semilinear PDE
\[
\frac{\partial v}{\partial t}=\frac{\Delta v}{2}-\frac{1}{2}v^2.
\]
One of us [\citet{myt98w}] extended this exponential duality and hence
proved uniqueness in law for $\sigma(u)=u^p,  u\ge0$ where $1/2<p<1$.
The dual process is then a solution to an SPDE driven by a one-sided
stable process.
Pathwise uniqueness among nonnegative solutions remains unsolved for
$0< p\le3/4$; see below for $p>3/4$.
\end{example}

\begin{example}\label{ex2} If $\sigma(x,u)=\sqrt{g(x,u)u}, u\ge0$,
where $g$
is smooth, bounded, and bounded away from $0$, then any kind of uniqueness
for solutions to~\eqref{generalized} is unresolved except when $g$ is
constant. Such
equations arise as weak limit points of the branching particle systems
as in
Example~\ref{ex1}, but where the branching and death rates of a particle at $x$ in
population $u^N$ is $Ng(x,u^N)/2$.
\end{example}

\begin{example}\label{ex3}
If $\sigma(x,u)=\sqrt{u(1-u)}$, $u\in
[0,1]$, then
solutions to \eqref{generalized} are population densities for the stepping
stone model on the line. That is, $u(t,x)$ is the proportion of a particular
allele type at location $x$ in a population undergoing Brownian
migration and
resampling between generations. For this model, uniqueness in law holds
by a moment
duality argument [see \citet{shi88}], and pathwise uniqueness remains
unresolved.
\end{example}

\begin{example}\label{ex4}
In this example, we no longer require $u$ to be
nonnegative. Consider $\sigma(u)=\sqrt{|u|}$ for $u\in\R$;
that is, consider the SPDE
%
%e1.2 #&#
\begin{equation}
\label{annbranch} \frac{\partial u}{\partial t}(t,x) =
\frac{\Delta}{2} u(t,x) + \sqrt {\bigl|u(t,x)\bigr|}
\dot{W}(t,x).
\end{equation}
This equation arises as a weak limit of the signed particle density of
two branching random walks, one with positive mass and one with
negative mass,
which annihilate each other upon collision. More precisely, consider two
particle systems on $\Z/\sqrt N$, one with positive mass and the other with
negative mass. Each particle independently produces offspring of the
same sign at a
randomly chosen nearest neighbor at rate $N/2$ and dies at rate $N/2$. The
systems interact when particles collide, and then there is pairwise
annihilation. Define $U^{N,\pm}(t,x)$ as in Example~\ref{ex1} where one considers
separately the positive and negative masses. Extend these functions by linear
interpolation to $x\in\R$. If $U^{N,\pm}(0,\cdot)\to u^{\pm
}(0,\cdot)$
uniformly for
some limiting cadlag (right-continuous with left limits) functions with
compact support satisfying
$u^+(0,\cdot)u^-(0,\cdot)\equiv0$, then $\{(U^{N,+},U^{N,-})\dvtx N\in\N
\}$ is
tight in the Skorokhod space of cadlag $C(\R)$-valued paths, where the
latter space of continuous functions has the topology of uniform convergence
on compact sets.
Any weak limit point $(u^+,u^-)$ will satisfy
%
%e1.3 #&#
\begin{eqnarray}
\label{annbranch2} \frac{\partial u^{\pm}}{\partial t}(t,x) &=&\frac{\Delta}{2} u^{\pm
}(t,x) +
\sqrt{u^\pm(t,x)} \dot{W_\pm}(t,x) -\dot K_t,
\nonumber
\\[-8pt]
\\[-8pt]
\nonumber
u^+(t,x)u^-(t,x)&\equiv& 0,
\nonumber
\end{eqnarray}
where $\dot W_+$ and $\dot W_{-}$ are independent space--time white
noises and
$K_t$ is a continuous nondecreasing process taking values in the space of
finite measures on the line with the topology of weak convergence. The
space--time
measure $K(dt,dx)$ records the time and location of the killing
resulting from the particle collisions. It is then easy to check that
$u=u^+-u^-$ satisfies \eqref{annbranch}. No results about uniqueness were
known for this process. The above convergence was proved in an earlier
draft of this
article but we have not included it as the details are a bit lengthy, if
routine. The convergence will only be used to help our intuition in what
follows.
\end{example}

In general, pathwise uniqueness of solutions, that is, the fact that
two solutions
with the same white noise and initial condition must coincide a.s., implies
the uniqueness of their laws; see, for example, \citet{kur07}. Quite
different duality arguments give uniqueness in law in Examples \ref{ex1} and \ref{ex3}, at
least among nonnegative solutions. But this kind of duality argument is
notoriously
nonrobust, and the interest in pathwise uniqueness stems in part from the
hope that such an approach would apply to a broader class of examples,
including perhaps Examples \ref{ex2} and \ref{ex4}.

It has long been hoped that pathwise uniqueness holds in \eqref{generalized}
if $\sigma$ is $\gamma$-H\"older continuous in the solution $u$ for
$\gamma\ge1/2$, since \citet{yw71} showed the corresponding
result holds for finite-dimensional stochastic differential equations (SDEs).
They proved that if $\sigma_i\dvtx \R\to\R$ is H\"older continuous of index
$1/2$ and
$b_i\dvtx \R^d\to\R$ is Lipschitz continuous, then solutions to
\[
dX^i_t=\sigma_i\bigl(X^i_t
\bigr)\,dB^i_t+b_i(X_t)\,dt, \qquad i=1,
\ldots,d
\]
are pathwise unique. Note that \eqref{generalized} has the same ``diagonal
form'' as the above SDE albeit in infinitely many dimensions. It was
\citet{vio75b} who first noted Yamada and Watanabe's proof extends to
infinite dimensional equations such as \eqref{generalized} if the
noise is
white in time but has a bounded covariance kernel in the spatial variable.
This proof breaks down for noise that is white in both time and space,
since in the $t$ variable, solutions are
H\"older continuous of index $(1/4)-\ep$ for all $\ep>0$, but not
H\"older continuous of index $1/4$. Hence, solutions are too rough in the
time variable to be semimartingales. Nonetheless in
\citet{mp11} a more involved extension of the Yamada--Watanabe
argument was established which proved pathwise uniqueness in
\eqref{generalized} if $\sigma(x,\cdot)$ is H\"older continuous of index
$\gamma>3/4$, uniformly in $x$.

This leads to the natural question of sharpness in this last result,
that is:
%
%e1.4 #&#
\begin{eqnarray}\label{mainqu}
\begin{tabular}{p{300pt}@{}}
Does pathwise uniqueness fail in general for
\eqref{generalized} if $\sigma(x,\cdot)=\sigma(\cdot )$
is $\gamma$-H\"older continuous for $\gamma\le3/4$, and in
particular for $\gamma=1/2$?
\end{tabular}
\end{eqnarray}
For the corresponding SDE, the Yamada--Watanabe result is shown to be
essentially sharp by Girsanov's equation
%
%e1.5 #&#
\begin{equation}
\label{girseq} X_t=\int_0^t\bigl|X_s\bigr|^\gamma\,
dB_s
\end{equation}
for which one solution is $X_t=0$. If $\gamma<1/2$, there are nonzero
solutions to
(\ref{girseq}), and so solutions are neither pathwise unique nor
unique in
law; see Section V.26 in \citet{rw87}. This suggests we
consider the SPDE
%
%e1.6 #&#
\begin{eqnarray}
\label{spde} \frac{\partial u}{\partial t}(t,x) &=& \frac{\Delta}{2} u(t,x) +
\bigl|u(t,x)\bigr|^\gamma\dot{W}(t,x),
\nonumber
\\[-8pt]
\\[-8pt]
\nonumber
u(0,x) &=& 0.
\nonumber
\end{eqnarray}

To state our main result we need some notation.
A superscript $k$, respectively~$\infty$, indicates that
functions are in addition $k$ times, respectively infinitely often, continuously
differentiable. A subscript $b$, respectively $c$, indicates that they
are also
bounded (together with corresponding derivatives), respectively have
compact support.
Let
$\langle f,g\rangle=\int_{\R}f(x)g(x) \,dx$ denote the $L^2$ inner product.
Set
\[
\|f\|_{\lambda}:=\sup_{x \in\R}\bigl |f(x)\bigr| e^{\lambda|x|},
\]
and define $C_{\mathrm{rap}}:=\{f \in C(\R)\dvtx \|f\|_{\lambda}
< \infty\mbox{ for any } \lambda>0\}$, endowed with the topology induced
by the norms $\|\cdot\|_{\lambda}$ for $\lambda>0$. That is,
$f_n\to f$ in $C_{\mathrm{rap}}$ if and only if $d(f,f_n)=\sum_{k=1}^\infty
2^{-k}(\Vert f-f_n\Vert_{k}\wedge1)\to0$ as $n\to\infty$.
Then $(C_{\mathrm{rap}},d)$ is a Polish space. The space $C_{\mathrm{rap}}$ is a
commonly used state space
for solutions to \eqref{generalized}; see \citet{shi94}.

We assume in \eqref{generalized} that $\dot W$ is a white noise on the
filtered probability space
$(\Omega,\CF,\CF_t,P)$, where $\CF_t$ satisfies the usual
hypotheses. This
means $W_t(\phi)$ is an $\CF_t$-Brownian motion with variance
$\Vert\phi\Vert_2^2t$ for each $\phi\in L^2(\R,dx)$ and $W_t(\phi
_1)$ and
$W_t(\phi_2)$ are independent if
$\langle\phi_1,\phi_2\rangle=0$. A stochastic
process $u\dvtx \Omega\times\R_+\times\R\to\R$ which is
$\CF_t\mbox{-previsible}\times\mbox{Borel}$ measurable will be
called a
solution to the
SPDE~\eqref{generalized} with initial condition $u_0\dvtx \R\to\R$ if
for each $\phi\in C_c^\infty(\R)$,
%
%e1.7 #&#
\begin{eqnarray}\label{SPDEMP}
\langle u_t,\phi\rangle&=&\langle u_0,\phi\rangle+\int
_0^t \biggl\langle u_s,
\frac{\Delta}{2}\phi \biggr\rangle \,ds
\nonumber
\\[-8pt]
\\[-8pt]
\nonumber
 &&{}+\int_0^t\int\sigma
\bigl(x,u(s,x)\bigr)\phi(x)W(ds,dx) \qquad\mbox{for all $t\ge0$ a.s.}
\end{eqnarray}
(The existence of all the integrals is of course part of the
definition.) We
often write $u_t$ for $u(t,\cdot)$. We use the framework of
\citet{wal86} to define stochastic integrals with respect to
$W(ds,dx)$. For
$u_0\in C_{\mathrm{rap}}$, we say $u$ is a $C_{\mathrm{rap}}$-valued solution if, in
addition, $t\to u(t,\cdot)$ has continuous $C_{\mathrm{rap}}$-valued paths
for all
$\omega$.

Here is our main result which answers question \eqref{mainqu},
at least for $\gamma<3/4$.

%th1.1 #&#
\begin{theorem} \label{thmmainresult} If $0<\gamma<3/4$, there is a
$C_{\mathrm{rap}}$-valued solution $u(t,x)$ to \eqref{spde} such that with positive
probability, $u(t,x)$ is not identically zero. In particular,
uniqueness in
law and pathwise uniqueness fail for \eqref{spde}.
\end{theorem}

This leaves open the state of affairs for $\gamma=3/4$ where, based on analogy
with the SDE, one would guess that uniqueness holds. Our theorem does,
however, dampen the hope of handling many of the SPDE's in the above examples
through a Yamada--Watanabe type theorem. It also shows that the SPDE in
Example~\ref{ex4} does not specify a unique law.

%Following on work of Mytnik, Perkins, and Sturm \citet{mps06}, Mytnik
%and
%Perkins \citet{mp11} showed uniqueness for (\ref{generalized}) if $f$
%is H\"older continuous of index
%$\alpha>\frac{3}{4}$. The proof drew on the Yamada-Watanabe
%method, and used the following intuition, which we illustrate for the
%less
%general case (\ref{spde}). The bad part of $|u|^\alpha$ is near $u=0$,
%since
%that is where the function fails to be Lipschitz. On the other hand,
%for
%values of $(t,x)$ near points where $u=0$, the noise term
%$|u(t,x)|^\alpha\dot{W}(t,x)$ is small, which makes the solution
%smoother than
%usual there. But if $u$ is smooth then it cannot escape quickly from
%0, so
%the region where $u$ is close to 0 must be larger than we expect.
%Repeating
%the argument, we see that there is an even larger region where $u$ is
%close to
%0 and hence smoother than expected. This leads to an even slower rate
%of
%growth for $u$ near the points where $u=0$, and so on. Iterating this
%reasoning a countable number of times yields the critical value
%$\alpha=\frac{3}{4}$.

A standard construction of a nonzero solution to Girsanov's SDE
proceeds as
follows. Start an ``excursion'' from $\pm\varepsilon$, run it until
it hits $0$, and then proceed to the next excursion, starting with the
opposite sign. The process consisting of $\pm\varepsilon$ jumps will disappear
as $\varepsilon\to0$ due to the alternating signs. For $\gamma<1/2$, a~diffusion calculation shows that the rescaled return time of the
diffusion is
in the domain of attraction of a stable subordinator of index
$(2(1-\gamma))^{-1}<1$, and the limiting jumps will lead to nontrivial
excursions in the scaling limit. With a bit of work one can
do the same in \eqref{spde} for $\gamma<1/2$. That is, one can seed randomly
chosen bits of mass of size $\pm\varepsilon$ and run the SPDE until it
hits $0$
and try again. Theorem~4 of \citet{bmp10} carries out
this argument and gives Theorem~\ref{thmmainresult} for $\gamma<1/2$.
Therefore, in the rest of this work we will assume
%
%e1.8 #&#
\begin{equation}
\label{gammch}1/2\le\gamma<3/4.
\end{equation}

When $\gamma\ge1/2$ the above excursion argument breaks down as the
time to
construct a nontrivial excursion will explode. Instead we start excursions
which overlap in time and deal with the potential spatial overlap of
positive and negative excursions. As Example~\ref{ex4} suggests we will
annihilate mass when the overlap occurs. Much of the challenge will be
to show that this overlap can be quite small if $\gamma<3/4$.

We now outline our strategy for constructing a nontrivial solution to
\eqref{spde}. Let $M_F(E)$ denote the space of finite measures on the metric
space $E$ with the weak topology. We will also use $\mu(\phi)$ and
$\langle\mu,\phi\rangle$ to denote integral of a function $\phi$
against a
measure $\mu$. Below we will construct
$\eta^+_{\ep}, \eta^-_{\ep}\in M_F([0,1]^2)$, both of which
converge to
Lebesgue measure $dt\,dx$ on the unit square as $\ep\downarrow0$, and
we will
also construct nonnegative solutions $U^\ep(t,x)$ and $V^\ep(t,x)$
with $0$
initial conditions to the equations
%
%e1.9 #&#
%e1.10 #&#
\begin{eqnarray}
\label{pmn-spde} \frac{\partial U^{\ep}}{\partial t}(t,x) &=& \dot\eta^+_{\ep}(t,x) +
\frac{\Delta}{2} U^{\ep}(t,x)+ U^\ep(t,x)^\gamma
\dot W^+(t,x)-\dot K^\ep_t,
\\
\label{pmn-spdeV} \frac{\partial V^{\ep}}{\partial t}(t,x) &=& \dot \eta ^-_{\ep}(t,x) +
\frac{\Delta}{2} V^{\ep}(t,x)+ V^\ep(t,x)^\gamma
\dot W^-(t,x)-\dot K_t^\ep.
\end{eqnarray}
Here $\dot W^+$ and $\dot W^-$ are independent white noises, and $t\to
K_t^\ep$
is a nondecreasing $M_F(\R)$-valued process. As suggested by
\eqref{annbranch2}, $K^\ep(dt,dx)$ will record the locations of the
pairwise annihilations resulting from the collisions between our two
annihilating populations. This construction will lead to the condition
\[
U^\ep(t,\cdot)V^\ep(t,\cdot)\equiv0.
\]
Note that $\eta^\pm_{\ep}$ are immigration terms. We will always
assume that
$\ep\in(0,1]$. If $\eta_\ep=\eta_\ep^+-\eta_\ep^-$, it is easy
to check that
$u_\ep=U^\ep-V^\ep$ satisfies
%
%e1.11 #&#
\begin{equation}
\label{n-spde} \frac{\partial u_{\ep}}{\partial t}(t,x) = \dot\eta_{\ep}(t,x) +
\frac
{\Delta}{2} u_{\ep}(t,x)+ \bigl|u_{\ep}(t,x)\bigr|^\gamma
\dot W(t,x)
\end{equation}
for an appropriately defined white noise $\dot W$. We will show that there
exists a subsequence $\ep_k$ such that as $k\to\infty$, $u_{\ep_k}(t,x)$
converges weakly in the Skorokhod space of $C_{\mathrm{rap}}$-valued paths
to a solution $u(t,x)$ of (\ref{spde}); see
Proposition~\ref{prop21}. $U^\ep$ is the positive part of $u_\ep$,
and so
Theorem~\ref{thmmainresult} will then follow easily from the following
assertion:
%
%cl1.2 #&#
\begin{claim}
\label{claim1}
There exists $\delta>0$ such that for all $\ep\in(0,1]$,
\[
P \biggl(\sup_{t\in[0,1]}\int U^{\ep}(t,x) \,dx>\delta
\biggr) > \delta.
\]
\end{claim}

If $N_\ep=\lfloor\ep^{-1}\rfloor$ (the greatest integer less than
$\ep
^{-1}$), the
measure $\eta_{\ep}$ will be obtained by smearing out spatial mass
using the time grid
%rid of jumps later?}}
%
%e1.12 #&#
\begin{equation}
\label{main-prob-est} \CG_{\ep} = \{k\ep/2 \dvtx 1\leq k\leq2N_\ep
\}.
\end{equation}
We further denote by $\CG_{\ep}^{\mathrm{odd}}$ the points of $\CG_{\ep}$
for which $k$ is odd, where $k$ is
%are the numerators
in the definition of $\CG_{\ep}$ above. We also define $\CG_{\ep
}^{\mathrm{even}}$
to be those grid points for which $k$ is even and let
%
%e1.13 #&#
\begin{equation}
\label{defJ} J_{\ep}^{x}(z)= \ep^{1/2} J
\bigl((x-z)\ep^{-1/2}\bigr),\qquad x,z\in\R,
\end{equation}
where $J$ is a nonnegative even continuous function bounded by $1$
with support in $[-1,1]$, and such that \mbox{$\int_{\R} J(z) \,dz=1$}.
Now let us enumerate points in $\CG_{\ep}^{\mathrm{odd}}$ and $\CG_{\ep
}^{\mathrm{even}}$ as follows:
\[
\{s_i, i\in\N_{\ep}\} = \CG_{\ep}^{\mathrm{odd}},\qquad
\{t_i, i\in \N_{\ep}\} = \CG_{\ep}^{\mathrm{even}}
,
\]
where $s_i=(2i-1)\frac{\ep}{2}$ and $t_i=2i\frac{\ep}{2}$ for
$i\in\N_\ep=\{1,\ldots, N_\ep\}$. Let $x_i, y_i, i=1,2,\ldots,$
be a sequence
of independent random variables distributed uniformly on $[0,1]$.

We define $\eta_{\ep}$
to be the signed measure
\begin{eqnarray*}
\eta_{\ep}(A) &=& \biggl[ \sum_{s_i\in\CG_{\ep}^{\mathrm{odd}}}\int
J_{\ep}^{x_i}(y)1_A(s_i,y) \,dy - \sum
_{t_i\in\CG_{\ep}^{\mathrm{even}}}\int J_{\ep
}^{y_i}(y)1_A(t_i,y)
\,dy \biggr]
\\
&\equiv&\eta_{\ep}^+(A)-\eta_{\ep}^-(A).
\end{eqnarray*}
It is easy to check that $\eta^\pm_\ep$ are as claimed above.

To simplify the outline of our proof, we will take $\gamma=1/2$ so that
we can
appeal to Example~\ref{ex4} for intuition. In later sections we do not make this
restriction on $\gamma$. We can then decompose $U^\ep=\sum_{i=1}^{N_\ep}U^i$
into descendants of the $i$th immigrant at $(s_i,x_i)$ (type $i$ particles)
and similarly write $V^\ep=\sum_{j=1}^{N_\ep}V^j$. We will suppress
$\varepsilon$
in the notation for clusters $U^i$ and $V^j$. We can also keep track of the
killed mass and, by adding these ghost particles back in, dominate
$U^\ep$ by
a super-Brownian motion $\bar U$ with immigration $\eta^+_\ep$, and dominate
the $\{U^i\}$ by independent super-Brownian motions $\{\bar U^i\}$
which sum
to $\bar U$. Similar processes $\bar V$ and $\{\bar V^j\}$ may be built to
bound the $V^\ep$ and $\{V^j\}$, respectively. We also can decompose
$K=\sum_iK^{i,U}=\sum_j K^{j,V}$ according to the type of individual being
killed. From hitting probabilities of Feller's branching diffusion
$\bar U^i(1)=\langle\bar U^i,1\rangle$, we know that with reasonably large
probability one of the $\bar U^i$ clusters does hit $1$, and we
condition on
such an event for a fixed choice of $i$, denoting the conditional law
by $Q_i$.
We now proceed in three steps:

\textit{Step} 1. $K^{i,U}_{s_i+t}(1)\le t^{3/2-\ep}$ for small $t$ with
reasonably large probability (see Lemma~\ref{thetabnd} below),
uniformly in $\ep$.

This step uses a modulus of continuity for the support of the
dominating
super-Brownian motions which states that they can spread locally no faster
than $t^{1/2}$ with some logarithmic corrections which we omit for the
purposes of this outline; see
Theorem~3.5 in \citet{mp92} for a more general version which
we will need for the general $\gamma$ case. This means both $\bar U^i$ and
$\bar V^j$ are constrained to lie inside a growing space--time parabola rooted
at their space--time birth points and hence the same is true for the
dominated processes
$U^i$ and $V^j$. If $\tau_j$ is the lifetime of $\bar V^j$ then, using the
known law of $\tau_j$ (it is the hitting time of zero by Feller's branching
diffusion starting from $\ep$) and a bit of geometry to see how large
$\tau_j$
has to be for the parabola of $\bar V^j$ to intersect with that of
$\bar U^i$
from $s_i$ to $s_i+t$, one can easily deduce that with reasonably large
probability
the only $\bar V^j$ clusters which can intersect with the $\bar U^i$ cluster
we have singled out are those born in the space--time rectangle
$[s_i,s_i+t]\times[x_i-2t^{1/2},x_i+2t^{1/2}]$. This means these are the
only $K^{j,V}$'s [killing by descendants of $(t_j,y_j)$] that can
contribute to
$K^{i,U}$ on $[s_i,s_i+t]$ since other $V$ particles will not collide
with the
$U^i$ mass. In particular, with reasonably large probability none of
the $V^j$
clusters born before $s_i$ can affect the mass of $U^i$ on
$[s_i,s_i+t]$; see
Lemma~\ref{presimass} for the proof of this last assertion for general
$\gamma$.
The mean amount of killing by these $V^j$'s can be no more than the mean
amount of immigration which fuels these populations. More precisely if one
integrates out the version of \eqref{pmn-spdeV} for $V^j$ over space, sums
over the above indices $j$ and bring the sum of the resulting $K^j$ to the
left-hand side, then one finds that if
\[
R_i=[s_i,s_i+t]\times\bigl[x_i-2t^{1/2},x_i+2t^{1/2}
\bigr],
\]
then
\[
E \biggl[\sum_{(t_j,y_j)\in R_i}K^j_{s_i+t}(1)
\biggr] \le E\bigl(\eta^-_\ep\bigl( [s_i,s_i+t]
\times\bigl[x_i-2t^{1/2},x_i+2t^{1/2}
\bigr]\bigr)\bigr)\le ct^{3/2}.
\]
A standard interpolation argument now shows the integrand on the
left-hand side is bounded by $ct^{3/2-\ep}$ for small enough $t$ a.s.,
and the claimed result follows from the above and the fact that any
killing by $K^{i,U}$ is matched by a killing on $V$ by one of the
$K^{j,V}$'s. It will turn out that for $\gamma<3/4$ one can get the
same bound on $K^i_t(1)$.

\textit{Step} 2. Under $Q_i$, which was the conditional law defined before step
1, $4\bar U^i_{s_i+\cdot}(1)$ is a $4$-dimensional $\mbox
{Bess}^2$-process and so
$\bar U^i_{s_i+t}(1)\ge t^{1+\ep}$ for small $t$ a.s.

This follows from a standard change of measure argument; see
Lemma~\ref{lembaresc} and its proof below. For general $\gamma<3/4$,
the mass\vadjust{\goodbreak}
$4\bar U^i_{s_i+\cdot}(1)$ will be a time change of a $4$-dimensional
$\mbox{Bess}^2$-process,
and one will be able to show that $\bar U^i_{s_i+t}(1)\ge t^\beta$ for
small $t$
a.s. for some $\beta<3/2$.

\textit{Step} 3. There is a reasonably large $Q_i$-probability (uniform in
$\ep$) that $U^i_{s_i+t}(1)\ge t^{1+\ep}$ for small $t$.

To see this, note that the above steps set up a competition
between the conditioning which
gives $\bar U^i(1)$ a positive linear drift and the killing which is limited
by step 1. To decide which effect wins when considering $U^i(1)$, we will
consider the ratio
\[
R_t=\frac{\bar U^i_{s_i+t}(1)-U^i_{s_i+t}(1)}{\bar U^i_{s_i+t}(1)}\in[0,1]
\]
of ghost particles to total population (alive and dead). An application of
It\^o's lemma will show that $R$ is a submartingale satisfying
\[
R_t=N_t+\frac{K_{s_i+t}(1)}{\bar U_{s_i+t}(1)},
\]
where $N_t$ is a continuous martingale. The last term is at most
$t^{1/2-2\ep}$
for small $t$ with reasonably large $Q_i$ probability by steps 1 and 2. We
localize to get the above behavior almost surely up to a stopping time,
take means
and use Kolmogorov's inequality for martingales to see that $R_t$ is
less than $1/2$ with
reasonably large probability, uniformly in $\ep$. By step 2 we can
conclude that on
this set $U^i_{s_i+t}(1)\ge(1/2) t^{1+\ep}$ for small $t$, and so
$U^i_{s_i+t}(1)$ is bounded away
from $0$ for small $t$ with reasonably large $Q_i$-probability
uniformly in $t$, as
required. This step is carried out in the proof of Proposition~\ref
{prop2} in
Section~\ref{secstochanal} below.

There are a number of problems when carrying out the above argument. In
step~1 we should
pay attention to the fact that the underlying probability is $Q_i$. In
addition, the argument for general $\gamma$ is more involved. For
example, the
clusters of the dominating processes $\bar V^j$ will no longer be
independent as they are when $\gamma=1/2$ due to the branching property
of solutions. Also, the rate of propagation results in
\citet{mp92} only apply for solutions where there is an
underlying historical process which records the ancestral histories of the
surviving population members. We could extend the construction of our
solutions to \eqref{pmn-spde} and \eqref{pmn-spdeV} to include such processes,
but this gets a bit unwieldy. Instead we prove a comparison theorem for
supports of solutions of parabolic SPDE's (Proposition~\ref{prop1}) which
allows us to derive these results from the corresponding property of solutions
of $\eqref{generalized}$ with $\sigma(u)=u^\gamma$. The latter
property holds
for any solution since these solutions are known to be unique in law by
\mbox{\citet{myt98w}}.

%re1.3 #&#
\begin{remark} The condition that $\gamma<3/4$ is required in step 1
to ensure that with
reasonably large probability, the $V$ particles born before time $s_i$
do not contribute to
the killing. Such killing, if it occurred, could lead to the
immediate annihilation of the $i$th seed with high probability. The
bound on
$\gamma$ is also used in steps~2 and 3 since otherwise the lower bound on
$\bar U^i_{s_i+t}(1)$ near $0$ will be $t^\beta$ for some $\beta>3/2$
which will be
of no use in keeping $R_t$ small for $t$ small.
\end{remark}

Here is an outline of the paper. Section~\ref{secsetup}
gives a careful description of the approximating solutions arising in
\eqref{pmn-spde}, \eqref{pmn-spdeV} and the various decompositions of these
processes. The actual construction of these approximate solutions is
carried out in
Appendix~\ref{secconstr}, while the fact that limit points of these
approximating solutions provide
actual solutions to \eqref{spde} is given in Appendix~\ref{secspdegr},
along with some standard moment bounds. In Section~\ref{secinclexcl} an
inclusion--exclusion argument reduces the
nonuniqueness result to a pair of Propositions (\ref{prop2} and
\ref{prop3}) which correspond to step 3 and an amalgamation of steps
1 and
2, respectively. In Section~\ref{secprop3} Proposition~\ref{prop3}
is then
reduced to a sequence of 5 lemmas, the main ones being Lemma~\ref{lembaresc}
and Lemma~\ref{thetabnd}, corresponding to steps 2 and 1, respectively.
Section~\ref{secstochanal} deals with the main parts of
the proof rooted in stochastic analysis including the proofs of
Lemma~\ref{lembaresc} and Proposition~\ref{prop2}. Sections~\ref{seclem44} and \ref{secKgrowth}
deal with the main parts of the proof involving qualitative properties
of the
clusters including the proof of Lemma~\ref{thetabnd} (the growth rate
of the killing measure) in
Section~\ref{secKgrowth}. Section~\ref{secKgrowth} also gives the proof
of the comparison theorem for supports of solutions of certain SPDE's.

%s2 #&#
\section{Set-up of equations}\label{secsetup}
\setcounter{equation}{0}

In what follows we assume that $\gamma\in[1/2,3/4)$. We will carry out
the method outlined in the \hyperref[intro]{Introduction}.

Recall that $\N_\ep=\{1,\ldots,N_\ep\}$ where $N_\ep=\lfloor\ep
^{-1}\rfloor$.
For any Polish space $\mathbf{E}$, let $D(\R_+,\mathbf{E})$ be the
Skorokhod space of cadlag $\mathbf{E}$-valued paths with left limits in
$\mathbf{E}$, and define
\begin{eqnarray*}
D^{\ep}(\R_+,\mathbf{E})&=& D(\R_+,\mathbf{E})\cap C(\R _+\setminus\CG
_{\ep},\mathbf{E})
\\
&=& \mbox{the space of cadlag $\mathbf{E}$-valued functions on $\R_+$, whose
paths}
\\
&& \mbox{are continuous on any time interval $\biggl[\displaystyle\frac{(i-1)\ep}{2},
\frac{i\ep}{2}\biggr),1\le i\le 2N_\ep$, % $[s_i, t_i), [t_i. s_{i+1}), i\in\NN_{\ep}$
}
\\
&&\mbox{and on $[N_{\ep}\ep,\infty)$.}
\end{eqnarray*}
We will construct a sequence of processes
% $(x_i, s_i)\in\CG_{\ep}^{\mathrm{odd}}$ (resp. $(y_i, t_i)\in\CG_{
$\{(U^{i,\ep},V^{i,\ep}), i\in\N_{\ep}\}$ with sample paths in
$(C(\R_+\setminus\CG_{\ep},C^+_{\mathrm{rap}})\cap D^\ep(\R_+,L^1(\R))^2$.
For each $\phi\in C_b^2(\R)$, w.p.1, $U^i,V^j$ (we will suppress
$\ep$ in our notation) will satisfy the
following equations for all $t\ge0$ and all $i,j\in\NN_\ep$.
Recall that $J^{x_i}$ was defined in (\ref{defJ}):
%
%e2.1 #&#
\begin{eqnarray}
\label{UVdefn} %\nonumber
\qquad\cases{ %
 U^i_t(\phi) = \bigl\langle
J^{x_i},\phi\bigr\rangle\1(t\geq s_i)
\vspace*{2pt}\cr
\hspace*{40pt}\displaystyle{} + \int_0^t \int_{\R}U(s,x)^{\gamma-1/2}
U^i(s,x)^{1/2} \phi(x) W^{i,U}(ds,dx)
\vspace*{2pt}\cr
 \hspace*{40pt}\displaystyle{}+ \int_0^t
U^{i}_s\biggl(\frac{1}{2}\Delta\phi\biggr) \,ds -
K^{i,U}_t(\phi),
\vspace*{6pt}\cr
V^j_t(\phi) = \bigl\langle J^{y_j},\phi\bigr
\rangle\1(t\geq t_j)
\vspace*{2pt}\cr
\hspace*{40pt}\displaystyle{}+ \int_0^t\int_{\R}
V(s,x)^{\gamma-1/2} V^j(s,x)^{1/2} \phi(x)
W^{j,V}(ds,dx)
\vspace*{2pt}\cr
 \hspace*{40pt}\displaystyle{}+ \int_0^t
V^{j}_s\biggl(\frac{1}{2}\Delta\phi\biggr) \,ds -
K^{j,V}_t(\phi), %\label{UVdefn}
\vspace*{6pt}\cr
 %\end{eqnarray}
\displaystyle\hspace*{180pt}\mbox{with } U_t = \sum
_i U^i_t, %\\
V_t =\sum_i V^i_t,}\hspace*{-10pt}
\end{eqnarray}
where, as will be shown in Proposition~\ref{thm11}, $U$ and $V$ have
paths in $D^\ep(\R_+,\break C^+_{\mathrm{rap}})$.
Here $W^{i, U}, W^{j,V},\/ i,j\in\NN_{\ep}$ are independent space time
white noises. $K^{i,U}, K^{j,V}$ and hence $K_t$ below, are all
right-continuous nondecreasing\break  $M_F(\R)$-valued processes representing
the mutual killing of the two kinds of particles, such that
%
%e2.2 #&#
\begin{equation}
\label{eq22} \sum_i K^{i,U}_t
= \sum_j K^{j,V}_t =:
K_t
\end{equation}
and
%
%e2.3 #&#
\begin{equation}
\label{eq23} U_t(x)V_t(x)= 0 \qquad\forall t\geq0, x\in\R.
\end{equation}
That is, $U$ and $V$ have disjoint supports and hence the same is true
of $U^i$ and $V^j$ for all $i,j\in\NN_\ep$.
It follows from \eqref{UVdefn} with $\phi\equiv1$ that for $t<s_i$,
$K^{i,U}_t(1)+U^i_t(1)$ is a continuous nonnegative local martingale,
hence supermartingale, starting at~$0$. Therefore
$K^{i,U}_t=U^i_t = 0,   t< s_i$ and similarly $K^{j,V}_t=V^j_t=0,
t<t_j$ for all $i,j \in\NN_{\ep}$.
One can think of $U$ and $V$ as two populations with initial masses
immigrating at times $s_i, i\in\N_\ep$ and $t_j, j\in\N_\ep $,
respectively.
Condition (\ref{eq23}) implies the presence of a ``hard killing'' mechanism
in which representatives of both populations annihilate each
other whenever they meet. The meaning of the ``hard killing'' notion will
become clearer when we will explain the construction of the equations as
limits of so-called soft-killing models.

We can regard $K^{i,U}$ and $K^{j,V}$ as the ``frozen'' mass that was killed
in corresponding populations due to the hard killing. If we reintroduce this
mass back we should get the model without killing. To this end let us
introduce the equations for ``killed'' populations which we denote by
$\tU^i, \tV^j$. These will take values in the same path space as
$U^i$, $V^j$.
For each $\phi\in C_b^2(\R)$, we require the following equations hold almost
surely for all $t\ge0$ and $i,j\in\NN_\ep$:
%
%e2.4 #&#
\begin{eqnarray}
\label{tUVdefn}\qquad \cases{ %
\displaystyle\tU^i_t(
\phi) = \int_0^t \int_{\R}
\bigl[ \bigl(\tU(s,x)+U(s,x) \bigr)^{2\gamma} - U(s,x)^{2\gamma}
\bigr]^{1/2}
\vspace*{2pt}\cr
\hspace*{65pt}\displaystyle{} \times\sqrt{\frac{\tU^i(s,x)}{\tU(s,x)}} \phi(x) \tW ^{i,U}(ds,dx)
\vspace*{2pt}\cr
\hspace*{40pt}\displaystyle{} + \int_0^t \tU^{i}_s
\biggl(\frac{1}{2}\Delta\phi\biggr) \,ds + K^{i,U}_t(\phi
),
\vspace*{2pt}\cr
\displaystyle\tV^j_t(\phi) = \int_0^t
\int_{\R} \bigl[ \bigl(\tV(s,x)+V(s,x) \bigr)^{2\gamma}
- V(s,x)^{2\gamma} \bigr]^{1/2}
\vspace*{2pt}\cr
\displaystyle\hspace*{65pt}{} \times\sqrt{\frac{\tV^j(s,x)}{\tV(s,x)}} \phi(x) \tW ^{j,V}(ds,dx)
\vspace*{2pt}\cr
\displaystyle\hspace*{40pt}{} + \int_0^t \tV^{j}_s
\biggl(\frac{1}{2}\Delta\phi\biggr) \,ds + K^{j,V}_t(
\phi), %\end{array}
\vspace*{2pt}\cr
 \hspace*{180pt}\displaystyle\mbox{with } \tU_t = \sum_i
\tU^i_t, \tV_t = \sum
_j \tV^j_t, }\hspace*{-14pt}
\end{eqnarray}
where, as will be shown in Proposition~\ref{thm11}, $\tU$ and $\tV$
have paths in $D^\ep(\R_+,\break C^+_{\mathrm{rap}})$ and we define $\sqrt{0/0}=0$
in the stochastic integral. The white noises $\tW^{i, U}$, $\tW^{j,V}$,
$i,j\in\NN_{\ep} $, are independent and also independent of
$\{ W^{i, U}, W^{j,V},\break  i,j\in\NN_{\ep}\}$. Again it is easy to see that
%
%e2.5 #&#
\begin{equation}
\label{0early} \tU^i_t=0\qquad\mbox{for }t<s_i\quad
\mbox{and}\quad\tV_t^j=0\qquad\mbox{for }t<t_j, i,j
\in\NN_\ep.
\end{equation}
Then using stochastic calculus, we deduce that the
processes defined by $\bU^i_t\equiv U^i_t+\tU^i_t, \bV^i\equiv
V^i_t+\tV^i_t$
satisfy the following equations for each $\phi$ as above, w.p.1 for all
$t\ge0$, $i,j\in\NN_\ep$:
%
%e2.6 #&#
{\fontsize{10.2}{12.2}\selectfont
\begin{eqnarray}
\label{eq25}\qquad\hspace*{2pt} %\nonumber
\cases{ %
\displaystyle\bU^i_t(\phi) =\bigl\langle J^{x_i},\phi\bigr
\rangle\1(t\geq s_i)+ \int_0^t \bU
^{i}_s\biggl(\frac{1}{2}\Delta\phi\biggr) \,ds
\vspace*{2pt}\cr
\displaystyle\hspace*{37pt}{} +\int_0^t\int_{\R}
\sqrt{U(s,x)^{2\gamma-1}U^i(s,x)+ \bigl(\bU(s,x)^{2\gamma} -
U(s,x)^{2\gamma} \bigr)\frac{\tU
^i(s,x)}{\tU
(s,x)}}
\vspace*{2pt}\cr
\displaystyle\hspace*{75pt}{}\times \phi(x) \bW^{i,U}(ds,dx),
\vspace*{2pt}\cr
\displaystyle\bV^j_t(\phi) = \bigl\langle
J^{y_j},\phi\bigr\rangle\1(t\geq t_j)+ \int
_0^t \bV ^{j}_s\biggl(
\frac{1}{2}\Delta\phi\biggr) \,ds
\vspace*{2pt}\cr
\displaystyle\hspace*{37pt}{} +\int_0^t\int
_{\R} \sqrt{V(s,x)^{2\gamma-1}V^j(s,x)+ \bigl(
\bV(s,x)^{2\gamma} - V(s,x)^{2\gamma} \bigr)\frac{\tV
^j(s,x)}{\tV
(s,x)}}
\vspace*{2pt}\cr
\displaystyle\hspace*{75pt}{}\times \phi(x) \bW^{j,V}(ds,dx),
\vspace*{2pt}\cr
\hspace*{180pt}\displaystyle\mbox{with } \bU_t = \sum_i
\bU^i_t, %\\
\bV_t = \sum
_j \bV^j_t,}
\end{eqnarray}}
\hspace*{-3pt}where, $\{\bW^{i,U},\bW^{j,V}, i,j\in\NN_{\ep}\}$ is again a
collection of independent white noises.
In spite of the complicated appearance of (\ref{eq25}), for $\bU,
\bV
$ we easily~get
%
%e2.7 #&#
\begin{eqnarray}
\label{eq28} \cases{ %
\displaystyle\bU_t(\phi)
= \int_0^t \int\phi(x)\eta^+_{\ep}(ds,dx)
+ \int_0^t \bU _s\biggl(
\frac{1}{2}\Delta\phi\biggr) \,ds
\vspace*{2pt}\cr
\displaystyle\hspace*{37pt}{} +\int_0^t\int_{\R}
\bU(s,x)^{\gamma} \phi(x) \bW^{U}(ds,dx),\qquad t\geq0,
\vspace*{6pt}\cr
\displaystyle\bV_t(\phi) = \int_0^t\int
\phi(x)\eta^-_{\ep}(ds,dx)+ \int_0^t
\bV _s\biggl(\frac{1}{2}\Delta\phi\biggr) \,ds
\vspace*{2pt}\cr
\displaystyle\hspace*{37pt}{} +\int_0^t\int_{\R}
\bV(s,x)^{\gamma} \phi(x) \bW^{V}(ds,dx),\qquad t\geq0, }
\end{eqnarray}
for independent white noises $\bW^{U}$ and $\bW^{V}$.
%It follows from
One can easily derive from the proof of Theorem~1 of \citet{myt98w} that
$(\bar U,\bar V)$ is unique in law (see Remark~\ref{rem09081} below).

Our next proposition establishes existence of solutions to the above systems
of equations. The filtration $(\mathcal{F}_t)$ will always be
right-continuous and such that $\mathcal{F}_0$ contains the $P$-null
sets in
$\mathcal{F}$. For any $T\geq1$, the space $D^{\ep}([0,T],\mathbf
{E})$ is
defined in the same way as $D^{\ep}(\R_+,\mathbf{E})$, but for
$\mathbf{E}$-valued functions on $[0,T]$.

For any function $f\in D(\R_+, \R)$, we set $\bolds{\Delta}
f(t)\equiv f(t)-f(t-)$, for any $t\geq0$.

%pr2.1 #&#
\begin{prop}
\label{thm11}
There exists a sequence $(U^i,V^i,\tU^i, \tV^i,\bU^i,\bV^i,
K^{i,U},\break  K^{i,V})_{i\in\NN_{\ep}}$ of processes
in
\begin{eqnarray*}
&&\bigl(\bigl(C\bigl([0,T]\setminus\CG_{\ep}, C^+_{\mathrm{rap}}\bigr)
\\
&&{}\quad\cap D^{\ep}\bigl([0,T], L^1(\R)\bigr)\bigr)^4
\times D^{\ep}\bigl(\R_+,C^+_{\mathrm{rap}}\bigr)^2\times
D^{\ep}\bigl(\R_+, M_F(\R)\bigr)^2
\bigr)^{N_{\ep}},
\end{eqnarray*}
which satisfies
(\ref{UVdefn})--(\ref{eq28}). Moreover,
$(U,V,\tU,\tV)\in D^\ep(\R_+,C^+_{\mathrm{rap}})^4$, and the following
conditions hold:
\begin{longlist}[(a)]
\item[(a)] For any $i\in\NN_{\ep}$,
$\bU^i_{s_i+\cdot}\in C(\R_+,C^+_{\mathrm{rap}})$, $\bV^i_{t_i+\cdot
}\in
C(\R_+,C^+_{\mathrm{rap}})$ and
\[
\bU^i(s,\cdot)=0,\qquad s<s_i,\qquad \bV^i(s,\cdot)=0,\qquad
s<t_i.
\]
\item[(b)]
$K_\cdot$ only has jumps at times in $\CG_{\ep}$, and
\end{longlist}
%
%e2.8 #&#
\begin{equation}
\label{Kjumps} \sup_t\bolds{\Delta} K_t(1)\le
\ep.
\end{equation}
\end{prop}
In what follows we will call $\bU^i$, $\bV^i$ (resp., $U^i$, $V^i$) the
clusters of the processes~$\bU$, $\bV$
(resp., $U$, $V$).

Now with all the processes in hand, let us state the results which will
imply the nonuniqueness
in~(\ref{spde}) with zero initial conditions.
First define
%
%e2.9 #&#
\begin{equation}
u_{\ep}(t):= U_t-V_t\in C_{\mathrm{rap}}
\end{equation}
and recall that $U_t,V_t$ implicitly depend on $\ep$.
Then it is easy to see from the above construction that $u_{\ep}$
satisfies the following SPDE:
%
%e2.10 #&#
\begin{eqnarray}
\label{eqspde-u-ep} \bigl\langle u_{\ep}(t),\phi\bigr\rangle&=&\sum
_{i}\bigl\langle J^{x_i}\1(t\geq s_i),
\phi\bigr\rangle- \sum_{j}\bigl\langle
J^{y_j}\1(t\geq t_j),\phi\bigr\rangle
\nonumber
\\[-8pt]
\\[-8pt]
\nonumber
&&{} +\int_0^t \frac{1}{2}\bigl\langle
u_{\ep}(s),\Delta\phi\bigr\rangle \,ds + \int_0^t
\int\bigl|u_{\ep}(s,x)\bigr|^{\gamma} \phi(x) W(ds,dx)
\end{eqnarray}
for $\phi\in C^2_b(\R)$.

The following two propositions will imply Theorem~\ref{thmmainresult}.
%
%pr2.2 #&#
\begin{prop}
\label{prop21} Let $\ep_n=\frac{1}{n}$. Then
$\{ u_{\ep_n}\}_{n}$ is tight in $D(\R_+,C_{\mathrm{rap}})$. If $u$ is any
limit point
as $\ep_{n_k}\downarrow0$, then $u$ is a $C_{\mathrm{rap}}$-valued
solution of the SPDE (\ref{spde}).
\end{prop}
The next proposition is just a restatement of Claim~\ref{claim1}.
%
%pr2.3 #&#
\begin{prop}
\label{prop22}
There exists $\delta_{\scriptsize{\ref{prop22}}},\ep_{\scriptsize{\ref{prop22}}}>0$ such that
for all $\ep\in(0,\ep_{\scriptsize{\ref{prop22}}}]$,
\[
P \biggl(\sup_{t\in[0,1]}\int U^\ep_t(x)
\,dx>\delta_{\scriptsize{\ref
{prop22}}} \biggr) > \delta_{\scriptsize{\ref{prop22}}}.
\]
\end{prop}
The proof of Proposition~\ref{prop21} will be standard and may be
found in Appendix~\ref{secspdegr}. Most of the paper is devoted to the
proof of
Proposition~\ref{prop22}.

%s3 #&#
\section{Outline of the proof of Proposition~\texorpdfstring{\protect\ref
{prop22}}{2.3}}\label{secinclexcl}
\setcounter{equation}{0}

We analyze the behavior of the clusters $U^i, V^i$ and show that with
positive probability at least one of them survives. As in the previous
section, we suppress dependence on the parameter $\ep\in(0,1]$.\vadjust{\goodbreak}

To make our analysis precise we need to introduce the event $A_i$ that
the mass of the cluster
$\bar U^i$ reaches $1$ before the cluster dies. Define
\begin{eqnarray*}
\btau_i &=&\inf\bigl\{t\dvtx \bU^i_{s_i+t}(1)=1
\bigr\},
\\
A_i&\equiv& \{\btau_i<\infty\},
\end{eqnarray*}
so that $\btau_i$ is an $(\CF_{s_i+t})$-stopping time.
Since we will often assume that one of $A_i$ occurs with positive
probability, we define the conditional probability measure $Q_i$,
%{\it Carl, Ed: We assume that there is a fixed grid in time with mesh $
%of every type and the central point of immigrated mass is distributed
%uniformly in %space interval
% $[0,1]$}
%
%e3.1 #&#
\begin{equation}
\label{Qdef} Q_i(A) = P(A | A_i) \qquad\forall A\in\CF.
\end{equation}

We need the following elementary lemma whose proof is given in
Section~\ref{secstochanal}.

%le3.1 #&#
\begin{lemma}
\label{lem32}
For all $1\le i,j\le N_\ep$, the events $A_i=A_i(\ep)$ satisfy:
\begin{longlist}[(a)]
\item[(a)]
$P(A_i) = \ep$;
\item[(b)]
$P(A_i\cap A_j) = \ep^2,   i\neq j$.
\end{longlist}
\end{lemma}

A simple inclusion--exclusion lower bound on $P(\bigcup_{i=1}^{\lfloor
2^{-1}\ep^{-1}\rfloor}A_i)$ shows that
for $\ep\le1/4$, with probability at least $3/16$,
at least one cluster of $\bar U^i$ survives until it attains mass $1$.
We will focus on the corresponding $U^i$ and to show it is nonzero with
positive probability (all uniformly in $\ep$), and we will
establish a uniform (in~$\ep$) escape rate. Set
%
%e3.2 #&#
\begin{equation}
\label{betadef}\beta=\frac{3/2-\gamma}{2(1-\gamma)},
\end{equation}
and
note that $\beta<3/2$ for $\gamma<3/4$. Our escape rate depends on a
parameter $\delta_1\in(0,1)$ (which will eventually be taken small
enough depending on $\gamma$) and is given in the event
\[
B_i(t) = \bigl\{ U^i_{s_i+s}(1)\geq
\tfrac{1}{2} s^{\beta+\delta
_1}, \forall s\in\bigl[\ep^{2/3}, t\bigr]
\bigr\}.\]

Denote the closed support of a measure $\mu$ on $\R$ by $S(\mu)$. Let
\[
T_R =\inf\bigl\{t\dvtx \bigl\| \bU_t(\cdot)
\bigr\|_{\infty}\vee\bigl\|\bV_t(\cdot)\bigr\| _{\infty
}> R\bigr\},
\]
so that $(T_R-s_i)^+$ is an $\CF_{s_i+t}$-stopping time.
To localize the above escape rate we let $\delta_0\in(0,1/4]$ and
define additional $(\CF_{s_i+t})$-stopping times ($\inf\varnothing
=\infty
$) by
\begin{eqnarray*}
\rho_i^{\delta_0,\ep}&=& \rho_i=\inf\bigl\{ t\dvtx S
\bigl(\bar U^i_{s_i+t}\bigr)\not \subset \bigl[x_i-
\ep^{1/2}- t^{1/2-\delta_0}, x_i +\ep^{1/2}+t^{1/2-\delta_0}
\bigr]\bigr\},
\\
H_i^{\delta_1,\ep}&=&H_i=\inf\bigl\{t\ge0\dvtx
\bU^i_{t+s_i}(1)<(t+\ep )^{\beta
+\delta_1}\bigr\},
\\
\theta^{\delta_0,\ep}_i&=&\theta_i=\inf\bigl\{t\dvtx
K^{i,U}_{t+s_i}(1)> (t+\ep )^{3/2-2\delta_0}\bigr\},
\\
v_i^{\delta_0,\delta_1,\ep}&=&v_i=\btau_i\wedge
H_{i}\wedge\theta _i\wedge\rho_i
\wedge(T_R-s_i)^+.
\end{eqnarray*}

We now state the two key results and show how they lead to
Proposition~\ref{prop22}. The first result
is proved in Section~\ref{secstochanal} below using some stochastic
analysis and change of measure arguments. The second is reduced to a
sequence of lemmas in Section~\ref{secprop3}.
%
%pr3.2 #&#
\begin{prop}
\label{prop2} There are $\delta_{\scriptsize{\ref{prop2}}}(\gamma)>0$ and
$p=p_{\scriptsize{\ref{prop2}}}(\gamma)\in(0,1/2]$ such that if $0<2\delta_0\le
\delta_1\le\delta_{\scriptsize{\ref{prop2}}}$, then
\[
Q_i\bigl(B_i(t\wedge v_i)\bigr)
\geq1-5t^{p}\qquad \mbox{for all }t>0 \mbox{ and }\ep\in(0,1].
\]
\end{prop}

%pr3.3 #&#
\begin{prop}
\label{prop3}
For each $\delta_1\in(0,1)$ and small enough $\delta_0>0$, depending on
$\delta_1$ and $\gamma$,
there exists a nondecreasing function $\delta_{\scriptsize{\ref{prop3}}}(t)$, not
depending on $\ep$, such that
\[
\lim_{t\downarrow0} \delta_{\scriptsize{\ref{prop3}}}(t)=0,
\]
and for all $\ep,t\in(0,1]$,
\[
P \Biggl(\bigcup_{i\geq1}^{tN_\ep} \bigl(
\{v_i <t\}\cap A_i \bigr) \Biggr) \leq t
\delta_{\scriptsize{\ref{prop3}}}(t).
\]
\end{prop}

With these two propositions we can give the following:
%pa3.subsection.subsubsection.1 #&#
\begin{pf*}{Proof of Proposition~\protect\ref{prop22}}
Let $p=p_{\scriptsize{\ref{prop2}}}$ and $\delta(t)=\delta_{\scriptsize{\ref{prop3}}}(t)$.
Assume $t=t_{\scriptsize{\ref{prop22}}}\in(0,1]$ is chosen so that
$5t^{p}+t+\delta
(t)\leq1/2$. We claim that
%
%e3.3 #&#
\begin{eqnarray}
\label{Bilbnd} P \Biggl(\bigcup_{i=1}^{tN_\ep}
B_i(t) \Biggr)\geq\frac{t}{4}\qquad \forall\ep\in(0,t/8].
\end{eqnarray}
Choose $\delta_1>0$ as in Proposition~\ref{prop2}, then $\delta
_0\in
(0,\delta_1/2]$ as in Proposition~\ref{prop3} and finally $t=t_{\scriptsize{\ref
{prop22}}}$ as above. Then we have
\begin{eqnarray*}
P \Biggl(\bigcup_{i=1}^{tN_\ep}
B_i(t) \Biggr) &\geq& P \Biggl(\bigcup_{i=1}^{tN_\ep}
B_i(t\wedge v_i)\cap A_i\cap\{v_i
\geq t\} \Biggr)
\\
&\geq& P \Biggl(\bigcup_{i=1}^{tN_\ep}
B_i(t\wedge v_i)\cap A_i \Biggr)- P \Biggl(
\bigcup_{i=1}^{tN_\ep} A_i\cap
\{v_i< t\} \Biggr)
\\
&\geq& \sum_{i=1}^{tN_\ep} P \bigl(
B_i(t\wedge v_i)\cap A_i \bigr) - \sum
_{i=1}^{tN_\ep} \sum
_{j=1, j\neq i}^{tN_\ep}P ( A_i\cap A_j )
\\
&&{} - P \Biggl(\bigcup_{i=1}^{tN_\ep}
A_i\cap\{v_i< t\} \Biggr).
\end{eqnarray*}
Recall the definition of the conditional law $Q_i$, and use Lemma~\ref
{lem32}(b) to see that the above is at least
\begin{eqnarray*}
&&\sum_{i=1}^{tN_\ep} Q_i
\bigl( B_i(t\wedge v_i) \bigr)P(A_i) -
t^2N_\ep^2\ep^2 \mbox{} - P
\Biggl(\bigcup_{i=1}^{tN_\ep} A_i
\cap\{v_i< t\} \Biggr)
\\
&&\qquad\geq\ep(tN_\ep-1)- 5N_\ep\ep t^{1+p_{\scriptsize{\ref{prop2}}}} -
t^2 - t\delta _{\scriptsize{\ref{prop3}}}(t)
\\
&&\qquad\geq t\bigl[1-5t^{p_{\scriptsize{\ref{prop2}}}}-t-
\delta_{\scriptsize{\ref{prop3}}}(t)\bigr]-2\ep,
\end{eqnarray*}
where the next to last inequality follows by Lemma~\ref{lem32}(a) and
Propositions~\ref{prop2}, \ref{prop3}.
Our choice of $t=t_{\scriptsize{\ref{prop22}}}$ shows that for
$\ep\le t_{\scriptsize{\ref{prop22}}}/8$. The above is at least
$\frac{t}{2}-\frac{t}{4}=\frac{t}{4}$. It follows from the final
part of
\eqref{UVdefn} that for all $t\ge0$,
$\int U^\ep_t(x) \,dx\ge\max_i\int U^{i,\ep}_t(x) \,dx$. The proposition
follows
immediately from \eqref{Bilbnd}.
\end{pf*}

%s4 #&#
\section{Lower bounds on the stopping times: Proof of
Proposition~\texorpdfstring{\protect\ref{prop3}}{3.3}}\label{secprop3}
\setcounter{equation}{0}

In this section we reduce the proof of Proposition~\ref{prop3} to five
lemmas which will be proved in Sections~\ref{secstochanal}--\ref
{secKgrowth} below.
The bounds in this section may depend on the parameters $\delta_0$ and
$\delta_1$, but not $\ep$.
We introduce
%
%e4.1 #&#
\begin{equation}
\label{bardef} \bar\delta=\bar\delta(\gamma)=\tfrac{1}{3} \bigl(
\tfrac
{3}{2}-2\gamma \bigr)\in(0,1/6].
\end{equation}

%le4.1 #&#
\begin{lemma}\label{lembaresc} For $\delta_0>0$ sufficiently small, depending
on $\delta_1,\gamma$, there is a function
$\eta_{\scriptsize{\ref{lembaresc}}}\dvtx \R_+\to[0,1]$ so
that $\eta_{\scriptsize{\ref{lembaresc}}}(t)\rightarrow0$ as $t\downarrow0$, and
for all $t>0$ and $\ep\in(0,1]$,
\[
Q_i(H_i\leq\btau_i\wedge\rho_i
\wedge t)\leq\eta_{\scriptsize{\ref
{lembaresc}}}(t)+8\ep^{\delta_1}.
\]
\end{lemma}

%le4.2 #&#
\begin{lemma}\label{bartaubnd} For all $t>0$ and $\ep\in(0,1]$,
\[
Q_i\bigl(\btau_i\leq t\wedge(T_R-s_i)^+
\bigr)\leq2\gamma R^{2\gamma-1}t+\ep.
\]
\end{lemma}

%le4.3 #&#
\begin{lemma}\label{thetabnd} If $0<\delta_0\le\bar\delta$, there
is a
constant $c_{\scriptsize{\ref{thetabnd}}}$, depending on $\gamma$ and $\delta_0$,
so that
\[
Q_i(\theta_i< \rho_i\wedge t)\leq
c_{\scriptsize{\ref{thetabnd}}}(t\vee\ep )^{\delta
_0}\qquad\mbox{for all }\ep,t\in(0,1] \mbox{ and }s_i\le t.
\]
\end{lemma}
%
%analysis
%up to $\rho^i$ so we can focus on $K^i$ measure of parabola. It uses
%the modulus
%of continuity of $\bV^j$'s under $Q^i$ via comparison lemma and
%Mueller-Perkins modulus
%bounds. More specifically it will use the following result with $Q_i$
%in place of $P$ and $\rho_j^V$
%in place of $\rho_i^U=\rho_i$ which holds since the $\{\bar V_j\}$
%satisfy the same martingale problem under $Q_i$.

It remains to handle the $\rho_i$ and $T_R$. This we do under the
probability $P$.

%le4.4 #&#
\begin{lemma} \label{rhobnd} There is a constant $c_{\scriptsize{\ref{rhobnd}}}\ge
1$, depending on $\gamma$ and $\delta_0$, so that
\[
P \Biggl(\bigcup_{i=1}^{pN_\ep}\{
\rho_i\le t\} \Biggr)\le c_{\scriptsize{\ref
{rhobnd}}}(t\vee\ep)p\1(p\ge\ep)
\quad\mbox{for all $\ep,p,t\in(0,1]$}.
\]
%
%and $\eta_{\ref{rhobnd}}(t)\downarrow0$ as $t \downarrow0$.
\end{lemma}

%Finally we use the following which is immediate from the convergence of
%$\overline U^\vep_t(\cdot)$ in $C(\R_+,C^+_{\mathrm{rap}})$ as $\vep
%requires extending the convergence in Proposition~\ref{prop21} to the
%dominating processes which I suspect we will need in its proof. ????

%le4.5 #&#
\begin{lemma}\label{TRbnd} For any $\ep_0>0$ there is a function
$\delta_{\scriptsize{\ref{TRbnd}}}\dvtx (0,2]\to\R_+$ so that
$\lim_{t\to0}\delta_{\scriptsize{\ref{TRbnd}}}(t)=0$ and
\[
P\Bigl(\sup_{s< t,x\in\R}\bar U(s,x)\vee\bar V(s,x)>t^{-2-\ep_0}
\Bigr)\le t\delta _{\scriptsize{\ref{TRbnd}}}(t)\qquad\mbox{for all }\ep\in(0,1], t\in(0,2].
\]
%
%$\sup_{0<\vep\le1}P(T_R<2)=\eta_R\downarrow0$ as $R\to\infty$.
\end{lemma}
%
%If we could show
%then by Markov's inequality, the above probability is at most $2t^{1+(

%One suspects bounds like the above (or stronger first moment bounds)
%should be a byproduct of the proof of Proposition~\ref{prop21}.
%Shiga (CMS 94) should help--eg. Thm. 2.5 should help as maximal moment
%results do come out of Kolmogorov type estimates.

%Another approach would be to argue that for $it\le1/2$,
%(which may be true). Then an easy argument reduces the result to the
%old nonquantitative version:
%%\qed

Assuming the above five results it is now very easy to give the
following:

\begin{pf*}{Proof of Proposition~\ref{prop3}} For $\delta_1\in
(0,1)$ choose $\delta_0>0$ small enough so that the conclusion of
Lemmas~\ref{lembaresc} and \ref{thetabnd} hold. Then for $0<t\le
1\le
R$ and $0<\ep\le1$, using Lemma~\ref{rhobnd} with $p=t$, we have
\begin{eqnarray*}
&& P\Biggl(\bigcup_{i=1}^{tN_\ep}
\{v_i< t\}\cap A_i\Biggr)
\\
&&\qquad\le P\Biggl(\bigcup_{i=1}^{tN_\ep}
\{T_R< t+s_i\}\Biggr)+P\Biggl(\bigcup
_{i=1}^{tN_\ep}\bigl\{ \btau _i<t
\wedge(T_R-s_i)^+\bigr\}\cap A_i\Biggr)
\\
&&\qquad\quad{}+P\Biggl(\bigcup_{i=1}^{tN_\ep}\{
\rho_i< t\}\Biggr) +P\Biggl(\bigcup_{i=1}^{tN_\ep}
\{H_i< \btau_i\wedge\rho_i\wedge t\}\cap
A_i\Biggr)
\\
&&\qquad\quad{}+P\Biggl(\bigcup_{i=1}^{tN_\ep}\{
\theta_i< \rho_i\wedge t\}\cap A_i\Biggr)
\\
&&\qquad\le P(T_R< 2t)+\sum_{i=1}^{tN_\ep}Q_i
\bigl(\btau_i\le t\wedge (T_R-s_i)^+
\bigr)P(A_i)+c_{\scriptsize{\ref{rhobnd}}}(t\vee\ep)t\1(t\ge\ep)
\\
&&\qquad\quad{}  +\sum_{i=1}^{tN_\ep}Q_i(H_i<
\btau_i\wedge\rho_i\wedge t)P(A_i)+\sum
_{i=1}^{tN_\ep} Q_i(
\theta_i< \rho_i\wedge t)P(A_i).
\end{eqnarray*}
Now apply Lemma~\ref{lem32} and Lemmas~\ref{lembaresc}--\ref
{thetabnd} to bound the above by
\begin{eqnarray*}
&&P \Bigl(\sup_{s<2t,x\in\R}\bar U(s,x)\vee\bar V(s,x)> R \Bigr)+2
\gamma R^{2\gamma-1}t^2 +\ep t+c_{\scriptsize{\ref{rhobnd}}}t^2
\\
&&\qquad{}+t\eta_{\scriptsize{\ref{lembaresc}}}(t)+t8\ep^{\delta_1}+tc_{\scriptsize{\ref
{thetabnd}}}(t\vee
\ep)^{\delta_0}.
\end{eqnarray*}
We may assume without loss of generality that $\eta_{\scriptsize{\ref
{lembaresc}}}$ is
nondecreasing and $t\ge\ep$ (or else the left-hand side is $0$). Set
$R=t^{-2-\ep_0}$, where $\ep_0>0$ is
chosen so that $3-4\gamma-\ep_0(2\gamma-1)>0$ and use Lemma~\ref
{TRbnd} to
obtain the required bound with
\[
\delta_{\scriptsize{\ref{prop3}}}(t)=2\delta_{\scriptsize{\ref{TRbnd}}}(2t)+2\gamma
(2t)^{3-4\gamma-\ep_0(2\gamma-1)}+2c_{\scriptsize{\ref{rhobnd}}}t+\eta_{\scriptsize{\ref
{lembaresc}}}(t)+8t^{\delta_1}+c_{\scriptsize{\ref{thetabnd}}}t^{\delta_0}.
\]
This finishes the proof of Proposition~\ref{prop3}.
\end{pf*}

%s5 #&#
\section{Change of measure and stochastic analysis: Proofs of
Proposition~\texorpdfstring{\protect\ref{prop2}}{3.2} and Lemmas~\texorpdfstring{\protect\ref
{lembaresc}}{4.1} and \texorpdfstring{\protect\ref{bartaubnd}}{4.2}}\label{secstochanal}
\setcounter{equation}{0}

Define
\[
\bar\tau_i(0)=\inf\bigl\{t\ge0\dvtx \bU^i_{s_i+t}(1)=0
\bigr\}
\]
and
\[
\bar\tau_i(0,1)=\bar\tau_i(0)\wedge\bar
\tau_i,
\]
where $\bar\tau_i$ was defined at the beginning of Section~\ref{secinclexcl}.

It follows from \eqref{eq25} that
%
%e5.1 #&#
\begin{equation}
\label{barmassmart} \bU^i_{t+s_i}(1)=\vep+\bM^i_t,
\end{equation}
where $\bM^i$ is a continuous local $(\CF_{s_i+t})$-martingale starting
at $0$ at $t=0$ and satisfying
%
%e5.2 #&#
\begin{eqnarray}
\label{barMsqfn} \bigl\langle\bM^i\bigr\rangle_t&=&\int
_{s_i}^{s_i+t}\int U(s,x)^{2\gamma
-1}U^i(s,x)
\nonumber
\\[-8pt]
\\[-8pt]
\nonumber
&&\hspace*{37pt}{}+\bigl(
\bU(s,x)^{2\gamma
}-U(s,x)^{2\gamma}\bigr)\frac
{\widetilde U^i(s,x)}{\widetilde U(s,x)} \,dx \,ds.
\end{eqnarray}

%le5.1 #&#
\begin{lemma}\label{hittime}
There is a $c_{\scriptsize{\ref{hittime}}}=c_{\scriptsize{\ref{hittime}}}(\gamma)>0$ so that
\[
P\bigl(\bar\tau_i(0)>t\bigr)\le c_{\scriptsize{\ref{hittime}}}
\ep^{2-2\gamma}t^{-1}\qquad\mbox {for all }t>0.
\]
\end{lemma}
%
%pa5.subsection.subsubsection.1 #&#
\begin{pf}
It follows from \eqref{barMsqfn} that
%
%e5.3 #&#
\begin{eqnarray}
\label{barsqfnbnd} && \frac{d\langle\bM^i\rangle(t)}{dt}
\nonumber\\
&&\qquad=\int U(s_i+t,x)^{2\gamma-1}U^i(s_i+t,x)
\nonumber
\\
&&\hspace*{10pt}\qquad\quad{} +\bigl(\bar U(s_i+t,x)^{2\gamma}-U(s_i+t,x)^{2\gamma}
\bigr) \frac{\widetilde U^i(s_i+t,x)}{\widetilde U(s_i+t,x)} \,dx
\nonumber
\\[-8pt]
\\[-8pt]
\nonumber
&&\qquad\ge \int U(s_i+t,x)^{2\gamma-1}U^i(s_i+t,x)+
\widetilde U(s_i+t,x)^{2\gamma-1}\widetilde U^i(s_i+t,x)
\,dx
\\
&&\qquad\ge \int U^i(s_i+t,x)^{2\gamma}+\widetilde
U^i(s_i+t,x)^{2\gamma} \,dx
\nonumber
\\
&&\qquad\ge 2^{1-2\gamma}\int\bU^i(s_i+t,x)^{2\gamma}
\,dx.
\nonumber
\end{eqnarray}
If $\gamma>1/2$, the result now follows from Lemma~3.4 of \citet{mp92}.

If $\gamma=1/2$, then one can construct a time scale $\tau_t$ satisfying
$\tau_t\le t$ for $\tau_t\le\bar\tau_i(0)$, under which
$t\to U^i_{s_i+\tau_t}(1)$ becomes Feller's continuous state branching
diffusion. The required result then follows from well-known bounds on the
extinction time for the continuous state branching process; for
example, see equation~(II.5.12)
in \mbox{\citet{per02}}.
\end{pf}
%
%pr5.2 #&#
\begin{prop}\label{girs}
\[
Q_i(A)=\int_A\frac{\bU^i_{s_i+(\bar\tau_i\wedge t)}(1)}{\ep} \,dP\qquad \mbox{for all }A\in\CF_{s_i+t}, t\ge0.
\]
\end{prop}
%
%pa5.subsection.subsubsection.2 #&#
\begin{pf}
Since $\bar\tau_i(0,1)<\infty$ a.s. (by the previous
lemma) and $\bar U^i(1)$ remains at $0$ when it hits $0$, we have
%
%e5.4 #&#
\begin{equation}
\label{tau1}\1(\bar\tau_i<\infty)=\bU^i_{s_i+\bar
\tau
_i(0,1)}(1)\qquad
\mbox{a.s.}
\end{equation}
By considering $\bar\tau_i(0,1)\le t$ and $\bar\tau_i(0,1)>t$
separately we see that
%
%e5.5 #&#
\begin{equation}
\label{tau2} \bU^i_{s_i+(\bar\tau_i(0,1)\wedge t)}(1)=\bU^i_{s_i+(\bar\tau
_i\wedge
t)}(1)\qquad
\mbox{a.s. on }\{\bar\tau_i>t\}.
\end{equation}
If $A\in\CF_{s_i+t}$, then
%
%e5.6 #&#
\begin{eqnarray}
\label{pdec} P(A,\bar\tau_i<\infty) &=&P(A,\bar\tau_i
\le t)+P(A,t<\bar\tau_i<\infty)
\nonumber\\
&=&\int\1(A,\bar\tau_i\le t)\bU^i_{s_i+(\bar\tau_i\wedge t)}(1) \,dP
\\
&&{}+E\bigl(\1(A,\bar\tau_i>t)P(\bar\tau_i<\infty|
\CF_{s_i+t})\bigr).
\nonumber
\end{eqnarray}
By \eqref{tau1} and \eqref{tau2} on $\{\bar\tau_i>t\}$,
\begin{eqnarray*}
P(\bar\tau_i<\infty|\CF_{s_i+t})&=&E\bigl(
\bU^i_{s_i+\bar\tau
_i(0,1)}(1)|\CF _{s_i+t}\bigr)
\\
&=&\bU^i_{s_i+(\bar\tau_i(0,1)\wedge t)}(1)
\\
&=&\bU^i_{s_i+(\bar\tau_i\wedge t)}(1).
\end{eqnarray*}
Then from \eqref{pdec} we conclude that
%
%e5.7 #&#
\begin{equation}
\label{girs1} P(A,\bar\tau_i<\infty)=\int_A
\bU^i_{s_i+(\bar\tau_i\wedge
t)}(1) \,dP.
\end{equation}
If $A=\Omega$, we get
%
%e5.8 #&#
\begin{equation}
\label{girs2} P(\bar\tau_i<\infty)=E\bigl(\bU^i_{s_i+(\bar\tau_i\wedge t)}(1)
\bigr)=\bU ^i_{s_i}(1)=\ep.
\end{equation}
Taking ratios in the last two equalities, we see that
\[
Q_i(A)=\int_A\bU^i_{s_i+(\bar\tau_i\wedge t)}(1)/
\ep \,dP
\]
as required.
\end{pf}

\begin{pf*}{Proof of Lemma~\ref{lem32}} (a) Immediate from
\eqref{girs2}.
%$+s_i$'s.}

(b) Assume $i<j$. The orthogonality of the bounded continuous
$(\CF_t)$-martin\-gales
$\bar U^i_{t\wedge(s_i+\bar\tau^i(0,1))}(1)$ and $\bar U^j_{t\wedge
(t_j+\bar\tau^j(0,1))}(1)$ [see \eqref{eq25}] shows that
%
%e5.9 #&#
\begin{eqnarray}
\label{orthoU} && E \bigl[\bar U^i_{s_i+\bar\tau^i(0,1)}(1)\bar
U^j_{s_j+\bar\tau
^j(0,1)}(1) |\CF_{s_j} \bigr]\mathbf{1}
\bigl(s_i+\bar\tau_i(0,1)>s_j\bigr)
\nonumber
\\[-8pt]
\\[-8pt]
\nonumber
&&\qquad= \bar U^i_{s_j}(1)\ep\mathbf{1}\bigl(s_i+
\bar\tau_i(0,1)>s_j\bigr).
\end{eqnarray}
By first using \eqref{tau1} and then \eqref{orthoU}, we have
\begin{eqnarray*}
&&P(A_i\cap A_j)
\\
&&\qquad=E \bigl[\bar U^i_{s_i+\bar\tau^i(0,1)}(1)\bar U^j_{s_j+\bar\tau
^j(0,1)}(1)
\bigr]
\\
&&\qquad=E \bigl[\bar U^i_{s_i+\bar\tau^i(0,1)}(1)\mathbf{1}
\bigl(s_i+\bar\tau _i(0,1)\le s_j\bigr)E
\bigl[\bar U^j_{s_j+\bar\tau^j(0,1)}(1) |\CF _{s_j} \bigr] \bigr]
\\
&&\qquad\quad{}+E \bigl[E \bigl[\bar U^i_{s_i+\bar\tau^i(0,1)}(1)\bar
U^j_{s_j+
\bar\tau^j(0,1)}(1) |\CF_{s_j} \bigr]\mathbf{1}
\bigl(s_i+\bar\tau _i(0,1)>s_j\bigr) \bigr]
\\
&&\qquad=E \bigl[\bar U^i_{(s_i+\bar\tau^i(0,1))\wedge s_j}(1)\mathbf {1}
\bigl(s_i+\bar\tau_i(0,1)\le s_j\bigr)\ep
\bigr]
\\
&&\qquad\quad{}+E \bigl[\bar U^i_{s_j}(1)\ep\mathbf{1}
\bigl(s_i+\bar\tau _i(0,1)>s_j\bigr) \bigr]
\\
&&\qquad=E \bigl[\bar U^i_{(s_i+\bar\tau^i(0,1))\wedge s_j}(1) \bigr]\ep =\ep^2.
\end{eqnarray*}
\upqed\end{pf*}

\begin{pf*}{Proof of Lemma~\ref{lembaresc}} Clearly $\bar
M^i_{t\wedge\bar\tau_i}$ is a bounded $(\CF_{s_i+t})$-martingale under~$P$. Girsanov's theorem [see Theorem VIII.1.4 of \citet{ry99}] shows that
%
%e5.10 #&#
\begin{equation}
\label{girsdec}\bar M^i_{t\wedge\bar\tau_i}=\bar M^{i,Q}_t+
\int_0^{t\wedge\btau_i}\bar U^i_{s_i+s}(1)^{-1}
\,d\bigl\langle \bar M^i \bigr\rangle_s,
\end{equation}
where $\bar M^{i,Q}$ is an $(\CF_{s_i+t})$-local martingale under $Q_i$
such that $\langle\bar M^{i,Q}\rangle_t=\langle\bar M^i\rangle
_{t\wedge\bar\tau_i}$.

If $\bar X_t=\bar U^i_{s_i+(t\wedge\btau_i)}(1)$, for
\[
t\le\int_0^{\btau_i}\bar X_s^{-1}\,d
\bigl\langle\bar M^i\bigr\rangle_s\equiv
R_i,
\]
define $\tau_t$ by
%
%e5.11 #&#
\begin{equation}
\label{taudef}\int_0^{\tau_t}\bar
X_s^{-1}\,d\bigl\langle \bar M^i\bigr
\rangle_s=t.
\end{equation}
Since $\bar X_s>0$ and $\frac{d\langle\bar M^i\rangle_s}{ds}>0$ for
all $0\le s\le\btau_i$ $Q_i$-a.s. [see \eqref{barMsqfn}] this uniquely
defines $\tau$ under $Q_i$ as a strictly increasing continuous function
on $[0,R_i]=[0,\tau^{-1}(\btau_i)]$. By differentiating \eqref{taudef}
we see that
%
%e5.12 #&#
\begin{equation}
\label{diffid1} \frac{d}{dt}\bigl(\bigl\langle\bar M^i \bigr
\rangle\circ\tau\bigr) (t)=\bar X(\tau _t),\qquad t\le\tau^{-1}(
\btau_i).
\end{equation}
Let $N_t=\bar M^{i,Q}(\tau_t)$, so that
\[
Z_t\equiv\bar X (\tau_t)=\ep+N_t+t\qquad \mbox{for }t\le\tau ^{-1}(\btau_i),
\]
and by \eqref{diffid1} for $t$ as above,
\[
\langle N\rangle_t=\bigl\langle\bar M^i\bigr\rangle(
\tau_t)=\int_0^t Z_s
\,ds.
\]
Therefore we can extend the continuous local martingale $N(t\wedge\tau
^{-1}(\btau_i))$ for $t>\tau^{-1}(\btau_i)$ so that $4Z_t$ is the
square of a $4$-dimensional Bessel process; see Section~XI.1 of \citet
{ry99}. By the escape rate for $4Z$ [see Theorem~5.4.6 of \citet{kni81}]
and a comparison theorem for SDE [Theorem V.43.1 of \citet{rw87}] there
is a nondecreasing $\eta_{\delta_0}\dvtx \R_+\to[0,1]$ so that $\eta
_{\delta
_0}(0+)=0$ and if $T_Z=\inf\{t\dvtx Z_t=1\}$, and
\[
\Gamma(\vep,\delta_0)=\inf_{0<t\le T_Z}
\frac{Z(t)}{t^{1+\delta_0}},
\]
then
%
%e5.13 #&#
\begin{equation}
\label{gambnd} \sup_{0<\ep\le1}Q_i\bigl(\Gamma(\ep,
\delta_0)\le r\bigr)\le\eta_{\delta_0}(r).
\end{equation}
Clearly $T_Z=\tau^{-1}(\bar\tau_i)$ and so
\[
\inf_{0<u\le\btau_i}\frac{\bar X(u)}{\tau^{-1}(u)^{1+\delta
_0}}=\inf_{0<t\le T_Z}
\frac{\bar X(\tau_t)}{t^{1+\delta_0}}=\Gamma(\ep,\delta_0).
\]
That is,
%
%e5.14 #&#
\begin{equation}
\label{barXbnd1} \bar X(u)\ge\Gamma(\vep,\delta_0)
\tau^{-1}(u)^{1+\delta_0} \qquad\mbox {for all }0<u\le\btau_i.
\end{equation}

To get a lower bound on $\tau^{-1}(u)$, use \eqref{barsqfnbnd} to see
that for $s<\rho_i\wedge\btau_i$,
\begin{eqnarray*}
\frac{d\langle\bar M^i\rangle_s}{ds} %&\ge& \int U(s_i+s,x)^{2\gamma-1}U^i(s_i+s,x)+\widetilde U(s_i+s,x)^{2
%&\ge& \int U^i(s_i+s,x)^{2\gamma}+\widetilde U^i(s_i+s,x)^{2\gamma} dx\\
&\ge& 2^{1-2\gamma}
\int\1\bigl(x_i-\ep^{1/2}-s^{(1/2)-\delta_0}\le x\le
x_i+\ep^{1/2}+s^{(1/2)-\delta_0}\bigr)\\
&&\hspace*{38pt}{}\times
 \bar U^i(s_i+s,x)^{2\gamma} \,dx
\\
&\ge& 2^{1-2\gamma} \bigl[2\bigl(\ep^{1/2}+s^{(1/2)-\delta_0}\bigr)
\bigr]^{1-2\gamma}\bar X(s)^{2\gamma},
\end{eqnarray*}
where the bound on $s$ is used in the last line. Therefore for $\ep
/2\le s<\rho_i\wedge\btau_i$ there is a $c_1(\gamma)>0$ so that
\begin{eqnarray*}
\frac{d\tau^{-1}(s)}{ds}&=&\bar X_s^{-1}\frac{d\langle\bar
M^i\rangle
_s}{ds}
\\
&\ge&c_1(\gamma)s^{((1/2)-\delta_0)(1-2\gamma)}\bar X_s^{2\gamma
-1}
\\
&\ge& c_1(\gamma)\Gamma(\ep,\delta_0)^{2\gamma-1}s^{((1/2)-\delta
_0)(1-2\gamma)}
\tau^{-1}(s)^{(2\gamma-1)(1+\delta_0)},
\end{eqnarray*}
where \eqref{barXbnd1} is used in the last line. Therefore if $\ep\le
t\le\rho_i\wedge\btau_i$, then
\[
\int_{\ep/2}^t\frac{d\tau^{-1}(s)}{\tau^{-1}(s)^{(2\gamma
-1)(1+\delta
_0)}}\ge
c_1(\gamma)\Gamma(\ep,\delta_0)^{2\gamma-1}\int
_{\ep/2}^t s^{((1/2)-\delta_0)(1-2\gamma)} \,ds.
\]
If $\delta'_0=\delta_0(2\gamma-1)$, this in turn gives
\begin{eqnarray*}
\tau^{-1}(t)^{2-2\gamma-\delta_0'}&\ge& c_1(\gamma)\Gamma(\ep,
\delta _0)^{2\gamma-1} \biggl[t^{1+({1}/{2}-\delta_0)(1-2\gamma
)}- \biggl(
\frac{\ep}{2} \biggr)^{1+({1}/{2}-\delta_0)(1-2\gamma)} \biggr]
\\
&\ge&c_2(\gamma)\Gamma(\ep,\delta_0)^{2\gamma-1}t^{(3/2)-\gamma
+\delta'_0}.
\end{eqnarray*}
We have shown that if $\beta(\delta_0)=\frac{(3/2)-\gamma+\delta
'_0}{2-2\gamma-\delta'_0}$, then for $\ep\le t\le\rho_i\wedge
\btau_i$,
\begin{eqnarray*}
\tau^{-1}(t)&\ge&c_2(\gamma)^{1/(2-2\gamma-\delta'_0)}\Gamma(\ep,
\delta _0)^{{(2\gamma-1)}/{(2-2\gamma-\delta'_0)}}t^{\beta(\delta_0)}
\\
&\ge&c_2(\gamma)^{1/(2-2\gamma-\delta'_0)} \bigl(\Gamma(\ep,\delta
_0)\wedge 1\bigr)^2 t^{\beta(\delta_0)},
\end{eqnarray*}
where $\delta'_0<1/4$ is used in the last line.

Recall the definition of the constant $\beta\in[1,\frac{3}{2})$ from
\eqref{betadef}. Use the above in~\eqref{barXbnd1} to see that there is
a $c_3(\gamma)\in(0,1)$ so that for $\ep\le t\le\rho_i\wedge\btau
_i\wedge1$,
\begin{eqnarray*}
\bar X(t)&\ge&\bigl[c_3(\gamma) \bigl(\Gamma(\ep,
\delta_0)\wedge 1\bigr)\bigr]^4t^{\beta
(\delta_0)(1+\delta_0)}
\\
&>& (2t)^{\beta+\delta_1},
\end{eqnarray*}
provided that $c_3(\gamma)(\Gamma(\ep,\delta_0)\wedge1)>
2t^{\delta
_0}$, and $\delta_0$ is chosen small enough depending on $\delta_1$ and
$\gamma$. By \eqref{gambnd} we conclude that for $t\le1$, and $\ep
\in(0,1]$,
%
%e5.15 #&#
\begin{eqnarray}
\label{escaway0} && Q_i\bigl(\bar X_s
\le(2s)^{\beta+\delta_1}\mbox{ for some }\ep\le s\le \rho _i\wedge\bar
\tau_i\wedge t\bigr)
\nonumber\\
&&\qquad\le Q_i\bigl(\Gamma(\ep,\delta_0)\wedge1\le2
t^{\delta
_0}/c_3(\gamma)\bigr)
\\
\nonumber
&&\qquad\le \eta_{\delta_0}\bigl( 2t^{\delta_0}/c_3(
\gamma)\bigr)+\1 \bigl(2t^{\delta
_0}\ge c_3(\gamma)\bigr)\equiv
\eta_{\scriptsize{\ref{lembaresc}}}(t).
\end{eqnarray}
The above inequality is trivial for $t>1$ as then the right-hand side
is at least $1$.

Next note that since $Z_t=\bar X(\tau_t)$ for $t\le T_Z$, $\bar
X_u\equiv1$ for $u\ge\bar\tau_i$, and $4Z$ has scale function
$s(x)=-x^{-1}$ [see (V.48.5) in \citet{rw87}], we see that for $\ep
^{\delta_1}\le2^{-\beta-\delta_1}$,
%
%e5.16 #&#
\begin{eqnarray}
\label{escnr0}
Q_i\bigl(\bar X_t
\le(2\ep)^{\beta+\delta_1}\mbox{ for some }t\ge0\bigr)&\le& Q_i\bigl(4Z
\mbox{ hits }4(2\ep)^{\beta+\delta_1}\mbox{ before }4\bigr)
\nonumber\\
&=&\frac{s(4)-s(4\ep)}{s(4)-s(4\cdot2^{\beta+\delta_1}\ep^{\beta
+\delta
_1})}
\nonumber\\
&=&\frac{1-\ep}{2^{-\beta-\delta_1}\ep^{1-\beta-\delta
_1}-\ep
}
\nonumber
\\[-8pt]
\\[-8pt]
\nonumber
&=&\frac{1-\ep}{2^{-\beta-\delta_1}\ep^{-\delta_1}(\ep
^{1-\beta
}-2^{\beta+\delta_1}\ep^{\delta_1+1})}
\\
\nonumber
&\le& \frac{1-\ep}{2^{-\beta-\delta_1}\ep^{-\delta
_1}(\ep
^{1-\beta}-\ep)}
\\
\nonumber
&\le& 2^{\beta+\delta_1}\ep^{\delta_1}\le8\ep^{\delta_1}.
\end{eqnarray}
The above bound is trivial if $\ep^{\delta_1}>2^{-\beta-\delta_1}$.

We combine \eqref{escaway0} and \eqref{escnr0} to conclude that
\begin{eqnarray*}
&&Q_i\bigl(\bar X_s\le(s+\ep)^{\beta+\delta_1}
\mbox{ for some }0\le s\le \rho _i\wedge\bar\tau_i\wedge t
\bigr)
\\
&&\qquad\le Q_i\bigl(\bar X_s\le(2s)^{\beta+\delta_1}\mbox{
for some }\ep\le s\le \rho_i\wedge\bar\tau_i\wedge t
\bigr)
\\
&&\qquad\quad{} +Q_i\bigl(\bar X_s\le(2\ep)^{\beta+\delta_1}\mbox{
for some }0\le s\le\ep \bigr)
\\
&&\qquad\le \eta_{\scriptsize{\ref{lembaresc}}}(t)+8\ep^{\delta_1}.
\end{eqnarray*}
The result follows.
\end{pf*}

\begin{pf*}{Proof of Lemma~\ref{bartaubnd}} As in the previous proof
we set
\[
\bar X_t=\bar U_{s_i+(t\wedge\btau_i)}(1)=\ep+\bar M^i_{t\wedge\bar
{\tau}_i}.
\]
From \eqref{girsdec} we have under $Q_i$
%
%e5.17 #&#
\begin{equation}
\label{Xdec1} \bar X_t=\ep+\bar M_t^{i,Q}+\int
_0^{t\wedge\btau_i}\bar X_s^{-1}\,d
\bigl\langle\bar M^i\bigr\rangle_s,
\end{equation}
where $\bar M^{i,Q}$ is an $(\CF_{s_i+t})$-local martingale under
$Q_i$. Therefore $\bar X$ is a bounded nonnegative submartingale under
$Q_i$, and by the weak $L^1$ inequality
%
%e5.18 #&#
\begin{eqnarray}\label{meanbnd}
Q_i\bigl(\btau_i\le t\wedge(T_R-s_i)^+
\bigr)&=&Q_i \Bigl(\sup_{s\le
t\wedge
(T_R-s_i)^+} \bar
X_s\ge1 \Bigr)
\nonumber
\\[-8pt]
\\[-8pt]
\nonumber
&\le&\int\bar X_{t\wedge(T_R-s_i)^+} \,dQ_i.
\end{eqnarray}
It is not hard to show that $\bar M^{i,Q}$ is actually a martingale
under $Q_i$, but even without this we can localize and use Fatou's
lemma to see that the right-hand side of \eqref{meanbnd} is at most
%
%e5.19 #&#
\begin{equation}
\label{sqfnbnd} \ep+ E_{Q_i} \biggl[\int_0^t
\1\bigl(s\le(T_R-s_i)^+\wedge\btau_i\bigr)
\bar X_s^{-1}\,d\bigl\langle\bar M^i\bigr
\rangle_s \biggr]\equiv\ep+I.
\end{equation}

Next we use \eqref{eq25} and then the mean value theorem to see that
\begin{eqnarray*}
I&=&E_{Q_i} \biggl[\int_{s_i}^{s_i+t}\mathbf{1}
\bigl(s\le T_R\wedge (s_i+\btau _i)\bigr)
\\
&&\hspace*{46pt}{} \times\int \bigl(U(s,x)^{2\gamma
-1}U^i(s,x)+
\bar U(s,x)^{2\gamma}-U(s,x)^{2\gamma} \bigr)
\\
&&\hspace*{121pt}{}\times\widetilde U^i(s,x)\widetilde U(s,x)^{-1} \,dx\, \bar
U^i_s(1)^{-1} \,ds \biggr]
\\
&\le& \int_{s_i}^{s_i+t}E_{Q_i} \biggl[\1
\bigl(s\le T_R\wedge(s_i+\btau _i)\bigr)
\\
&&\hspace*{48pt}{} \times\int \bigl(U(s,x)^{2\gamma-1}U^i(s,x)+2\gamma\bar
U(s,x)^{2\gamma-1}\widetilde U^i(s,x) \bigr) \,dx\\
&&\hspace*{250pt}{}\times \bar
U^i_s(1)^{-1} \biggr] \,ds
\\
&\le& 2\gamma R^{2\gamma-1}\int_{s_i}^{s_i+t}E_{Q_i}
\biggl[\1(s\le s_i+\btau_i)\int\bar U^i(s,x)
\,dx \bar U^i_s(1)^{-1} \biggr] \,ds
\\
&\le& 2\gamma R^{2\gamma-1} t.
\end{eqnarray*}
We put the above bound into \eqref{sqfnbnd} and then use \eqref
{meanbnd} to conclude that
\[
Q_i\bigl(\btau_i\le t\wedge(T_R-s_i)^+
\bigr)\le\ep+2\gamma R^{2\gamma-1} t
\]
as required.
\end{pf*}

\begin{pf*}{Proof of Proposition~\ref{prop2}} Fix $i\le N_\ep$
and set
\[
X_t=U^i_{s_i+(t\wedge\bar\tau_i)}(1),\qquad D_t=
\tU^i_{s_i+(t\wedge\bar
\tau_i)}(1).
\]
If $f(x,d)=d/(x+d)$, then
%
%e5.20 #&#
\begin{equation}
\label{Rdefn} R_t\equiv\frac{\tU^i_{s_i+(t\wedge
\bar\tau
_i)}(1)}{\bar U^i_{s_i+(t\wedge\bar\tau_i)}(1)}=f(X_t,D_t)
\in[0,1].
\end{equation}
Proposition~\ref{thm11} shows that $X$ and $D$ are right-continuous
semimartingales with left limits.
We will work under $Q_i$ so that the denominator of $R$ is strictly
positive for all $t\ge0$ $Q_i$-a.s.
Our goal will be to show that $R$ remains small on $[0,t\wedge v_i]$
for $t$ small with high probability, uniformly in $\ep$. Then
$U^i_{s_i+s}(1)$ will be bounded below by a constant times $\bar
U_{s_i+s}(1)$ on this interval with high probability, and the latter
satisfies a uniform escape rate on the interval by the definition of $v_i$.

From Proposition~\ref{thm11}, and in particular \eqref{tUVdefn} and
\eqref{0early}, we have
\[
\tU^i_{s_i+(t\wedge\bar\tau_i)}(1)=\tM^i_t+K^{i,U}_{s_i+(t\wedge
\bar\tau_i)}(1),
\]
where $\tM^i$ is the continuous $(\CF_{s_i+t})$-local martingale (under
$P$) given by
\[
\tM^i_t=\int_{s_i}^{s_i+(t\wedge\bar\tau_i)}
\bigl(\bar U(s,x)^{2\gamma
}-U(s,x)^{2\gamma} \bigr)^{1/2}\sqrt{
\frac{\tU^i(s,x)}{\tU
(s,x)}}\tW ^{i,U}(ds,dx),
\]
and $K^{i,U}_{s_i+\cdot}$ is a right-continuous nondecreasing process.
By Girsanov's theorem [Theorem VIII.1.4 in \citet{ry99}] there is a
continuous $(\CF_{s_i+t})$-local martingale under $Q_i$, $\tM^{i,Q}$,
so that
%
%e5.21 #&#
\begin{eqnarray}
\label{tMQi}
 \tM^i_t&=&\tM^{i,Q}_t+
\int_{s_i}^{s_i+(t\wedge\bar\tau_i)}\bar U^i_s(1)^{-1}\,d
\bigl\langle\tM^i,\bar M^i\bigr\rangle_s
\nonumber
\\
&=&\tM^{i,Q}_t\\
&&{}+\int_{s_i}^{s_i+(t\wedge\bar\tau_i)}
\int \bigl(\bar U(s,x)^{2\gamma}-U(s,x)^{2\gamma} \bigr)
\frac{\tU^i(s,x)\tU
(s,x)^{-1}}{\bar
U^i_s(1)}\,dx \,ds.\nonumber
\end{eqnarray}

From \eqref{UVdefn} we have
\[
U^i_{s_i+(t\wedge\bar\tau_i)}(1)=\ep+M^i_t-K^{i,U}_{s_i+(t\wedge
\bar\tau_i)}(1),
\]
where $M^i$ is the continuous $(\CF_{s_i+t})$-local martingale (under $P$),
\[
M^i_t=\int_{s_i}^{s_i+(t\wedge\bar\tau_i)}\int
U(s,x)^{\gamma
-(1/2)}U^i(s,x)^{1/2} W^{i,U}(ds,dx).
\]
Another application of Girsanov's theorem implies there is a continuous
$(\CF_{s_i+t})$-local martingale under $Q_i$, $M_t^{i,Q}$, such that
%
%e5.22 #&#
\begin{equation}
\label{MQi} M_t^i=M_t^{i,Q}+\int
_{s_i}^{s_i+(t\wedge\bar\tau_i)}\int\frac
{U(s,x)^{2\gamma-1}U^i(s,x)}{\bar U^i_s(1)} \,dx \,ds.
\end{equation}
Note that $\langle M^i,\tM^i\rangle=0$ and so $M^{i,Q}$ and $\tM^{i,Q}$
are also orthogonal under~$Q_i$.

If
\begin{eqnarray*}
J_t=\sum_{s\le t} f(X_s,D_s)-f(X_{s-},D_{s-})&-f_x(X_{s-},D_{s-})
\Delta X_s
-f_d(X_{s-},D_{s-})\Delta D_s,
\end{eqnarray*}
then It\^o's lemma [e.g., Theorem VI.39.1 in \citet{rw87}] shows that
under $Q_i$,
%
%e5.23 #&#
\begin{eqnarray}
\label{ItoR}\quad R_t&=&R_0+\int_0^t
f_x(X_{s-},D_{s-})\,dX_s+\int
_0^t f_d(X_{s-},D_{s-})\,dD_s
\nonumber
\\
&&{}+\int_0^{t\wedge\bar\tau_i} \frac{1}{2}f_{xx}(X_{s-},D_{s-})
\int U(s_i+s,x)^{2\gamma-1}U^i(s_i+s,x) \,dx
\,ds
\nonumber
\\[-8pt]
\\[-8pt]
\nonumber
&&{}+\int_0^{t\wedge\bar\tau_i} \frac
{1}{2}f_{dd}(X_{s-},D_{s-})
\int\bigl[\bar U(s_i+s,x)^{2\gamma
}-U(s_i+s,x)^{2\gamma}
\bigr]
\\
\nonumber
&& \hspace*{126pt}{}
\times\tU^i(s_i+s,x)\tU (s_i+s,x)^{-1}
\,dx \,ds+J_t.
\end{eqnarray}
Since
\[
\Delta X_t=-\Delta K^{i,U}_{s_i+(t\wedge\bar\tau_i)}(1)=-\Delta
D_t,
\]
and $f_x=-d(x+d)^{-2}$, $f_d=x(x+d)^{-2}$, we conclude that
\begin{eqnarray*}
J_t&=&\sum_{s\le t} \bigl[
f(X_{s-}-\Delta D_s,D_{s-}+\Delta
D_s)-f(X_{s-},D_{s-})
\\
&&\hspace*{91pt}{}+[f_x-f_d](X_{s-},D_{s-})
\Delta D_s \bigr]
\\
&=&\sum_{s\le t}\frac{\Delta D_s}{X_{s-}+D_{s-}}-
\frac{\Delta
D_s}{X_{s-}+D_{s-}}=0.
\end{eqnarray*}
We use $f_{xx}=2d(x+d)^{-3}$, $f_{dd}=-2x(x+d)^{-3}$, \eqref{tMQi} and
\eqref{MQi} in \eqref{ItoR} to conclude that if $\bar X_t=\bar
U^i_{s_i+(t\wedge\bar\tau_i)}(1)$ and
\[
N_t=\int_0^t-D_{s-}
\bar X_s^{-2}\,dM^{i,Q}_s+\int
_0^t X_{s-}\bar X_s^{-2}
\,d\tM_s^{i,Q},
\]
then
%
%e5.24 #&#
\begin{eqnarray}
\label{Rdecomp}
\nonumber
R_t&=& R_0+N_t+
\int_0^{t\wedge\bar\tau_i}\bigl(-D_{s-}\bar
X_s^{-3}\bigr)\int U(s_i+s,x)^{2\gamma-1}
U^i(s_i+s,x) \,dx \,ds
\\
\nonumber
&&{}+\int_0^{t\wedge\bar\tau_i}D_{s-}\bar
X_s^{-2}\,dK_{s_i+s}^{i,U}(1)
\\
\nonumber
&&{}+\int_0^{t\wedge\bar\tau_i}X_s\bar
X_s^{-3}\int \bigl[\bar U(s_i+s,x)^{2\gamma}-U(s_i+s,x)^{2\gamma}
\bigr]
\\
&&\hspace*{81pt}{} \times\bar U^i(s_i+s,x)\tU(s_i+s,x)^{-1}\,dx
\,ds
\nonumber\\
\nonumber
&&{}+\int_0^{t\wedge\bar\tau_i}X_{s-}\bar
X_s^{-2}\,dK_{s_i+s}^{i,U}(1)
\\
&&{}+\int_0^{t\wedge\bar\tau_i}D_{s-}\bar
X_s^{-3}\int U(s_i+s,x)^{2\gamma-1}U^i(s_i+s,x)\,dx
\,ds
\\
\nonumber
&&{}-\int_0^{t\wedge\bar\tau_i} X_s\bar
X_s^{-3}\int\bigl[\bar U(s_i+s,x)^{2\gamma}-U(s_i+s,x)^{2\gamma}
\bigr]
\\
\nonumber
&&\hspace*{81pt}{} \times\tU^i(s_i+s,x)\tU(s_i+s,x)^{-1}\,dx
\,ds
\\
\nonumber
&=& R_0+N_t+\int_0^{t\wedge\bar\tau_i}
\bar X_s^{-1}\,dK^{i,U}_{s_i+s}(1).
\end{eqnarray}
Under $Q_i$, $N$ is a continuous $(\CF_{s_i+t})$-local martingale, and
the last term in \eqref{Rdecomp} is nondecreasing. It follows from
this and $R\in[0,1]$ that
%
%e5.25 #&#
\begin{equation}
\label{Rsub} R\mbox{ is an $(\CF _{s_i+t})$-submartingale under
}Q_i.
\end{equation}

As $R_0=K_{s_i}^{i,U}(1)/\ep$, integration by parts shows that
%
%e5.26 #&#
\begin{eqnarray}
\label{Rdecomp2}
R_t&=&R_0+N_t+
\frac{K^{i,U}_{s_i+(t\wedge\bar\tau
_i)}(1)}{\bar X_t}-\frac{K_{s_i}^{i,U}(1)}{\ep}-\int_0^tK_{s_i+s}^{i,U}(1)\,d
\biggl(\frac{1}{\bar X_s} \biggr)
\nonumber
\\[-8pt]
\\[-8pt]
\nonumber
&=&N_t-\int_0^tK_{s_i+s}^{i,U}(1)\,d
\biggl(\frac{1}{\bar X_s} \biggr)+\frac
{K^{i,U}_{s_i+(t\wedge\bar\tau_i)}(1)}{\bar U^i_{s_i+(t\wedge\bar
\tau_i)}(1)}.
\end{eqnarray}
Another application of It\^o's lemma using \eqref{barmassmart} and
\eqref{girsdec} shows that
\begin{eqnarray*}
\bar X_t^{-1}&=&\ep^{-1}-\int
_0^t\bar X_s^{-2}\,d\bar
X_s+\int_0^{t\wedge
\bar\tau_i}\bar
X_s^{-3}\,d\bigl\langle\bar M^i\bigr
\rangle_s
\\
&=&\ep^{-1}-\int_0^t\bar
X_s^{-2}\,d\bar M_s^{i,Q}-\int
_0^{t\wedge
\bar\tau
_i}\bar X_s^{-3}\,d
\bigl\langle\bar M^i\bigr\rangle_s+\int
_0^{t\wedge\bar\tau
_i}\bar X_s^{-3}\,d
\bigl\langle\bar M^i\bigr\rangle_s
\\
&=&\ep^{-1}-\int_0^t\bar
X_s^{-2}\,d\bar M_s^{i,Q}.
\end{eqnarray*}
Therefore $\bar X_t^{-1}$ is a continuous $(\CF_{s_i+t})$-local
martingale under $Q_i$ and hence the same is true of
$N^R_t=N_t-\int_0^tK^{i,U}_{s_i+s}(1)\,d (\frac{1}{\bar X_s} )$. From \eqref
{Rdecomp2} we have
%
%e5.27 #&#
\begin{equation}
\label{Rdecomp3} R_t=N_t^R+\frac{K^{i,U}_{s_i+(t\wedge\bar\tau_i)}(1)}{\bar
U^i_{s_i+(t\wedge\bar\tau_i)}(1)}.
\end{equation}

Recall from \eqref{Kjumps} and \eqref{eq22} that
%
%e5.28 #&#
\begin{equation}
\label{Kjump} \Delta K^{i,U}_{s_i+t}(1)\le\ep\qquad \mbox{for all }t
\ge0.
\end{equation}
Assume that (recall $\beta<3/2$)
%
%e5.29 #&#
\begin{equation}
\label{deltacond1} 0<2\delta_0\le\delta_1\le
\tfrac{1}{4} \bigl(\tfrac{3}{2}-\beta \bigr)\equiv
\delta_{\scriptsize{\ref{prop2}}}(\gamma).
\end{equation}
These last two inequalities (which give $\frac{3}{2}-\beta-\delta
_1-2\delta_0>0$) together with the continuity of $\bar U^i_{s_i+\cdot
}(1)$ [recall Proposition~\ref{thm11}(a)],\vadjust{\goodbreak} and the definitions of
$\theta_i\ge v_i$ and $H_i\ge v_i$ imply that
\[
\sup_{s\le v_i\wedge t}\frac{K^{i,U}_{s_i+s}(1)}{\bar
U^i_{s_i+s}(1)}\le\sup_{s\le v_i\wedge t}
\frac{(s+\ep
)^{(3/2)-2\delta
_0}+\ep}{(s+\ep)^{\beta+\delta_1}}\le(t+\ep)^{(3/2)-\beta
-2\delta
_0-\delta_1}+\ep^{1-\beta-\delta_1},
\]
and so from \eqref{Rdecomp3}
%
%e5.30 #&#
\begin{equation}
\label{NRbnd} \sup_{s\le v_i\wedge t} \bigl|N^R_s\bigr|
\le1+(t+\ep)^{(3/2)-\beta-2\delta
_0-\delta_1}+\ep^{1-\beta-\delta_1}<\infty.
\end{equation}

We now apply the weak $L^1$ inequality to the nonnegative
submartingale $R$
[recall \eqref{Rsub}] to conclude that ($\sup\varnothing=0$)
%
%e5.31 #&#
\begin{eqnarray}\label{maxineq}
&& Q_i \Bigl(\sup_{\ep^{2/3}\le s\le v_i\wedge
t}R_s
\ge1/2 \Bigr)\nonumber
\\
\nonumber
&&\qquad=E_{Q_i} \biggl[Q_i \biggl(\sup
_{\ep^{2/3}\le s\le
v_i\wedge
t}R_s\ge\frac{1}{2} \Big|
\CF_{\ep^{2/3}} \biggr)\mathbf {1}\bigl(v_i\wedge t\ge
\ep^{2/3}\bigr) \biggr]
\nonumber
\\[-8pt]
\\[-8pt]
\nonumber
&&\qquad\le 2E_{Q_i} \bigl[R_{v_i\wedge t}\mathbf{1}
\bigl(v_i\wedge t\ge \ep^{2/3}\bigr) \bigr]\1\bigl(t\ge
\ep^{2/3}\bigr)
\\
&&\qquad \le2E_{Q_i} \biggl[R_{(v_i\wedge t)-}+
\frac{\Delta
K^{i,U}_{s_i+(v_i\wedge t)}(1)}{\bar U^i_{s_i+(v_i\wedge
t)}(1)}\mathbf{1}\bigl(v_i\wedge t\ge\ep^{2/3}
\bigr) \biggr]\mathbf{1}\bigl(t\ge \ep ^{2/3}\bigr).
\nonumber
\end{eqnarray}
By \eqref{Kjump} and the definition of $H_i\ge v_i$ we have
%
%e5.32 #&#
\begin{eqnarray}\label{Rjumpbnd}
&&\frac{\Delta K^{i,U}_{s_i+(v_i\wedge t)}(1)}{\bar U^i_{s_i+(v_i\wedge
t)}(1)}\1\bigl(v_i\wedge t\ge\ep^{2/3}\bigr) \nonumber\\
&&\qquad
\le \frac{\ep}{(\ep+v_i\wedge t)^{\beta+\delta_1}}\1\bigl(v_i\wedge t\ge \ep^{2/3}\bigr)
\\
&&\qquad\le \ep^{1-(2/3)(\beta+\delta_1)}.\nonumber
\end{eqnarray}
From \eqref{Rdecomp3} and the definitions of $H_i\ge v_i$ and $\theta
_i\ge v_i$ we have
%
%e5.33 #&#
\begin{eqnarray}\label{Rminusbnd}
E_{Q_i} [R_{(v_i\wedge t)-} ]&=&E_{Q_i}
\bigl[N^R_{v_i\wedge
t}\bigr]+E_{Q_i}
\bigl[{K^{i,U}_{s_i+(v_i\wedge t)-}(1)}/{\bar U^i_{s_i+(v_i\wedge
t)}(1)}
\bigr]
\nonumber
\\
&\le& E_{Q_i} \bigl[\bigl(\ep+(v_i\wedge t)
\bigr)^{(3/2)-\beta
-2\delta_0-\delta_1} \bigr]\\
&\le&(\ep+t)^{(3/2)-\beta-2\delta
_0-\delta_1},\nonumber
\end{eqnarray}
where we used \eqref{NRbnd} to see that $N^R_{v_i\wedge t}$ is a mean
zero martingale and also applied~\eqref{deltacond1} to see the exponent
is positive. Inserting \eqref{Rjumpbnd} and \eqref{Rminusbnd} into
\eqref{maxineq} and using \eqref{deltacond1}, we get for $t\le1$,
%
%e5.34 #&#
\begin{eqnarray}
\label{Rbndaway} && Q_i \biggl(\sup_{\ep^{2/3}\le s\le v_i\wedge t}R_s
\ge\frac
{1}{2} \biggr)
\nonumber\\
&&\qquad\le \bigl[(\ep+t)^{(3/2)-\beta-2\delta_0-\delta
_1}+\ep ^{1-(2/3)(\beta+\delta_1)} \bigr]\1
\bigl(t\ge\ep^{2/3}\bigr)
\\
\nonumber
&&\qquad\le 2^{3/2}t^{(3/2)-\beta-2\delta_1}+t^{(3/2)-(\beta
+\delta
_1)}
\le5t^{(3/2)-\beta-2\delta_1}.
\end{eqnarray}
%
%Since $\beta<3/2$, we may take $\delta_1<\delta(\gamma)$ small enough
%so that
%$$\frac{3}{2}-\beta-2\delta_1\ge\frac{1}{2}\left(\frac{3}{2}-\beta
%The first inequality implies \eqref{deltacond1}, and
%the above inequalities, $\beta\ge1$, and \eqref{Rbndaway} show that
%for $\delta_1<\delta(\gamma)$ and $t\le1$,
Equation~\eqref{deltacond1} implies $(3/2)-\beta-2\delta_1\ge
(1/2)((3/2)-\beta
)$, and so for $t\le1$ we conclude
\[
Q_i \biggl(\sup_{\ep^{2/3}\le s\le v_i\wedge t}R_s\ge
\frac
{1}{2} \biggr)\le5t^{(1/2)((3/2)-\beta)}.
\]
The above is trivial for $t> 1$.
On $\{\sup_{\ep^{2/3}\le s\le v_i\wedge t}R_s<1/2\}$ we have for all
$s\in[\ep^{2/3},t\wedge v_i]$,
\[
U^i_{s_i+s}(1)\ge\tfrac{1}{2}\bar U^i_{s_i+s}(1)
\ge\tfrac
{1}{2}s^{\beta
+\delta_1},
\]
and so $B_i(t\wedge v_i)$ occurs. The result follows with $p_{\scriptsize{\ref
{prop2}}}=\frac{1}{2} (\frac{3}{2}-\beta )\in(0,\frac{1}{4}]$
(as~$\gamma\ge1/2$).
\end{pf*}%\qed

%s6 #&#
\section{Propagation speed of the supports and a comparison principle:
Proofs of Lemmas~\texorpdfstring{\protect\ref{rhobnd}}{4.4} and \texorpdfstring{\protect\ref
{TRbnd}}{4.5}}\label{seclem44}
\setcounter{equation}{0}

If $a>0$, $1>\gamma\ge1/2$ and $X_0\in C^+_{\mathrm{rap}}$, then Theorems~2.5 and
2.6 of \citet{shi94} show the existence of continuous $C^+_{\mathrm{rap}}$-valued
solutions to
%
%e6.1 #&#
\begin{equation}
\label{gammaspde} \frac{\partial X}{\partial t}=\frac{1}{2} \Delta X+a X^\gamma
\dot W,
\end{equation}
where as usual $\dot W$ is a space--time white noise on $\R_+\times\R$.
Theorem~1.1 of \citet{myt98w} then shows the laws
$\{P_{X_0}\dvtx X_0\in C^+_{\mathrm{rap}}\}$ of these processes on
$C(\R_+,C^+_{\mathrm{rap}})$ are unique.

We start with a quantified version of Theorem~3.5 of \citet{mp92}
applied to the particular equation \eqref{gammaspde}.

%le6.1 #&#
\begin{lemma}\label{qMP} Assume $X$ satisfies \eqref{gammaspde} with
$X_0=J_\ep^{x_0}$ for $x_0\in\R$ and $\ep\in(0,1]$. If $\gamma\in
(1/2,3/4)$ choose $\delta=\delta(\gamma)\in(0,1/5)$ sufficiently small
so that
$\beta_0=\beta_0(\gamma)=\frac{2\gamma-\delta}{1-\delta}\in(1,3/2)$
and for $N>1$, define
\[
T_N=\inf \biggl\{t\ge0\dvtx \int X(t,x)^\delta \,dx\ge N
\biggr\}.
\]
If $\gamma=1/2$, set $\beta_0=1$ and $T_N=\infty$. For $\delta_0\in
(0,1/4]$, define
%
%e6.2 #&#
\begin{equation}\quad
\label{rhodef2}\rho=\inf \bigl\{t\ge0\dvtx S(X_t)\not \subset
\bigl[x_0-\ep^{1/2}-t^{(1/2)-\delta_0}, x_0+
\ep^{1/2}+t^{(1/2)-\delta
_0}\bigr] \bigr\}.
\end{equation}
There is a $c_{\scriptsize{\ref{qMP}}}>0$ (depending on $\gamma$) so that
\[
P(\rho\le t\wedge T_N)\le c_{\scriptsize{\ref{qMP}}}a^{-1}N^{\beta_0-1}
\vep\exp \bigl(-t^{-\delta_0}/c_{\scriptsize{\ref{qMP}}}\bigr)\qquad\mbox{for all }\ep,t
\in(0,1].
\]
\end{lemma}
%
%pa6.subsection.subsubsection.1 #&#
\begin{pf}
Since $X$ is unique in law, the construction in
Section~4 of
\citet{mp92} allows us to assume the existence of a historical process
$H_t$, a~continuous $M_F(C)$-valued process, associated with $X$. Here $C$ is the
space of continuous $\R$-valued paths. $H$ will satisfy the martingale
problem $(M_{X_0})$ in \citet{mp92}, and the relationship with $X$ is that
%
%e6.3 #&#
\begin{equation}
\label{HX} H_t\bigl(\{y\in C\dvtx y_t\in B\}
\bigr)=X_t(B)\quad \mbox{for all $t\ge0$ and Borel subsets $B$ of }\R.\hspace*{-35pt}
\end{equation}
Hence the hypotheses of Theorem~3.5 of \citet{mp92} are satisfied with
$a_k\equiv a$ for all $k$. If $I_t=[x_0-\sqrt\ep-t^{(1/2)-\delta
_0},x_0+\sqrt\ep+t^{(1/2)-\delta_0}]$, that result implies
$S(X_t)\subset I_t$ for small enough $t$ a.s., but we need to quantify
this inclusion and so will follow the proof given there, pointing out
some minor changes and simplifications as we go.

If $\gamma=1/2$, $X$ is the density of one-dimensional super-Brownian
motion, and the argument in \citet{mp92} and its quantification\vadjust{\goodbreak} are both
much easier. As a result we will assume $3/4>\gamma>1/2$ in what
follows and leave the simpler case $\gamma=1/2$ for the reader. The
fact that $a_k=a$ for all~$k$ [i.e., for us $a(u)=au^\gamma$ for all
$u$ in the notation of \citet{mp92}], means that in the localization in
\citet{mp92}, the times $\{T_N\}$ may be chosen to agree with our
definition of $T_N$. We will work with the cruder modulus of
continuity, $\psi(t)=\frac{1}{2} t^{(1/2)-\delta_0}$, in place of the
more delicate $ch(t)=c(t\log^+(1/t))^{1/2}$ in \citet{mp92}, leading to
better bounds.

If
\begin{eqnarray}
G_{n,j,k}=\bigl\{y\in C\dvtx \bigl|y\bigl(k2^{-n}\bigr)-y
\bigl(j2^{-n}\bigr)\bigr|>\psi \bigl((k-j)2^{-n}\bigr)\bigr\},
\nonumber\\
\eqntext{0\le j<k; j,k,n\in\Z_+,}
\end{eqnarray}
and $B$ is a standard one-dimensional Brownian motion,
then for $k-j\le2^{n/2}$, (3.16) of \citet{mp92} becomes
\begin{eqnarray*}
&& Q_{X_0}\bigl(H_{(k+1)2^{-n}}(G_{n,j,k})>0, T_N\ge(k+1)2^{-n}\bigr)
\\
&&\qquad\le c_1 N^{\beta_0-1}a^{-1}2^n
X_0(1) P \bigl(\bigl|B\bigl(k2^{-n}\bigr)-B
\bigl(j2^{-n}\bigr)\bigr|>\psi\bigl((k-j)2^{-n}\bigr)
\bigr)^{2-\beta_0}
\\
&&\qquad\le c_2 N^{\beta_0 -1}a^{-1}2^n\ep\exp
\bigl(-\tfrac{1}{16} 2^{n\delta
_0} \bigr) \qquad\mbox{(recall $
\beta_0<3/2$).}
\end{eqnarray*}

Now we sum the above bound over $0\le j<k\le2^n$, $k-j\le2^{n/2}$,
$n\ge m$ and argue as in the proof of Theorem~3.5 in \citet{mp92} to see
that if
\[
\eta_m=c_3N^{\beta_0-1}a^{-1}\ep\exp
\bigl(-2^{(m\delta
_0/2)-4} \bigr),
\]
then with probability at least $1-\eta_m$,
\begin{eqnarray}
&&H_t(G_{n,j,k})=0\qquad\mbox{for all }0\le j<k\le2^n,
 k-j\le2^{n/2}, (k+1)2^{-n}\le T_N,
\nonumber\\
\eqntext{ t\ge(k+1)2^{-n},\mbox{ and }n\ge m.}
\end{eqnarray}
Rearranging this as in the proof of Theorem~3.5 of \citet{mp92}, we have
with probability at least $1-\eta_m$,
%
%e6.4 #&#
\begin{eqnarray}
\label{levymod}\qquad
\bigl|y\bigl(k2^{-n}\bigr)-y\bigl(j2^{-n}
\bigr)\bigr|\le\psi\bigl((k-j)2^{-n}\bigr)\qquad\mbox{for all }0\le j<k, k-j
\le2^{n/2},
\nonumber
\\[-6pt]
\\[-10pt]
\eqntext{(k+1)2^{-n}\le t \mbox{ and }n\ge m\mbox{ for }H_t\mbox{-a.a. }
 y\mbox{ for all }t\le T_N\wedge1.}\vadjust{\goodbreak}
\end{eqnarray}
Next, we can argue as in the last part of the proof of \citet{mp92},
which was
a slightly modified version of L\'evy's classical derivation of the exact
Brownian modulus of continuity, to see that \eqref{levymod} implies
\begin{eqnarray}
\bigl|y(v)-y(u)\bigr|\le2\psi\bigl(|v-u|\bigr)\qquad\mbox{for all }0\le u<v\le t\mbox{ satisfying }|v-u|
\le2^{-m/2}
\nonumber\\
\eqntext{\mbox{for }H_t\mbox{-a.a. } y  \mbox{ for all }t\le T_N
\wedge1.}
\end{eqnarray}
In particular, the above implies that
\[
P\bigl(\bigl|y(t)-y(0)\bigr|\le2\psi(t)\ H_t\mbox{-a.a. } y \mbox{ for all }t\le
2^{-m/2}\wedge T_N\bigr)\ge1-\eta_m.
\]

Now $H_t(|y(0)-x_0|>\sqrt\ep)$ is a nonnegative martingale starting at
$0$ by the martingale problem for $H$ [just as in the proof of
Corollary~3.9 in \citet{mp92}] and so is identically $0$ for all $t$
a.s. Therefore, the above and \eqref{HX} imply that
\[
P\bigl(\rho<2^{-m/2}\wedge T_N\bigr)\le
\eta_m.
\]
A simple interpolation argument now gives the required bound.
\end{pf}

%co6.2 #&#
\begin{cor}\label{qMP2} Assume $X$, $\delta_0$ and $\rho$ are as in
Lemma~\ref{qMP}. There is a \mbox{$c_{\scriptsize{\ref{qMP2}}}>0$}, depending on $a$,
$\delta_0$ and $\gamma$, so that
\[
P(\rho\le t)\le c_{\scriptsize{\ref{qMP2}}} \ep(t\vee\ep)\qquad\mbox{for all }t,\ep\in(0,1].
\]
\end{cor}
%
%pa6.subsection.subsubsection.2 #&#
\begin{pf}
We clearly may assume $x_0=0$ by translation
invariance. By Lem\-ma~\ref{qMP} with $N=N_0\equiv8$ and $\beta_0$,
$T_{N_0}$ as in that result, we have
%
%e6.5 #&#
\begin{equation}
\label{rhobndI} P(\rho\le t)\le c_{\scriptsize{\ref{qMP}}}a^{-1}8^{\beta_0-1}
\ep\exp \bigl(-t^{-\delta
_0}/c_{\scriptsize{\ref{qMP}}}\bigr)+P(t\wedge
T_{N_0}<\rho\le t).
\end{equation}
The result is now immediate if $\gamma=1/2$, so we assume $\gamma\in
(1/2,3/4)$. If $\delta\in(0,\frac{1}{5})$ is as in Lemma~\ref{qMP},
$I_s=[-\sqrt\ep-s^{(1/2)-\delta_0},\sqrt\ep+s^{(1/2)-\delta_0}]$, and
$0<t\le1$, then
%
%e6.6 #&#
\begin{eqnarray}\label{maxX}
&&
P(t\wedge T_{N_0}<\rho\le t)\nonumber\\
&&\qquad\le P(T_{N_0}<t
\wedge\rho)
\nonumber\\
&&\qquad\le P \biggl(\int_{I_s} X(s,x)^\delta
\,dx>8\mbox{ for some }s\le t\wedge\rho \biggr)\\
&&\qquad\le P \biggl( \biggl(\int X(s,x) \,dx \biggr)^\delta
|I_s|^{1-\delta
}>8\mbox{ for some }s\le t \biggr)
\nonumber\\
&&\qquad\le P\Bigl(\sup_{s\le t} X_s(1)>
\lambda\Bigr),\nonumber
\end{eqnarray}
where
$\lambda=8^{1/\delta} [[2 (\sqrt\ep+t^{(1/2)-\delta
_0}
)]^{(1-\delta)/\delta} ]^{-1}$. Recall that $X_t(1)$ is a
continuous nonnegative local martingale starting
at $\ep$, and so by the weak $L^1$ inequality and Fatou's lemma the right-hand
side of \eqref{maxX} is at most
\begin{eqnarray*}
\lambda^{-1}E\bigl[X_0(1)\bigr]&\le&\ep2^{-1-(2/\delta)}
\bigl(\sqrt\ep +t^{1/4} \bigr)^{(1-\delta)/\delta}\qquad \mbox{(by $
\delta_0\le 1/4$)}
\\
&\le& \ep \bigl[\max\bigl(t,\ep^2\bigr) \bigr]^{(1-\delta)/(4\delta)}
\\
&\le& \ep\max(t,\ep)\qquad \mbox{(since $\delta<1/5$)}.
\end{eqnarray*}
We use the above bound in \eqref{rhobndI} to conclude that
\begin{eqnarray*}
P(\rho\le t)&\le& \bigl[c_{\scriptsize{\ref{qMP}}}a^{-1}8^{\beta_0-1}\exp
\bigl(-t^{-\delta
_0}/c_{\scriptsize{\ref{qMP}}}\bigr)+(t\vee\ep) \bigr]\ep
\\
&\le& c_{\scriptsize{\ref{qMP2}}}(t\vee\ep)\ep.
\end{eqnarray*}
\upqed\end{pf}

The next proposition will allow us to extend the above bound to a
larger class of SPDEs. It will be proved at the end of this section.

%pr6.3 #&#
\begin{prop}
\label{prop1}
Let $a>0$, $1>\gamma\ge1/2$ and $Z$ be a continuous $C^+_{\mathrm{rap}}$-valued solution to
the following SPDE:
%
%e6.7 #&#
\begin{equation}
\frac{\partial Z}{\partial t}= \frac{1}{2}\Delta Z + \sigma(Z_s, s,
\omega)\dot W^1,
\end{equation}
where $\dot W^1$ is a space time white noise, $\sigma$ is Borel${}\times{}$previsable, and
\[
\sigma(y,s,\omega)\geq ay^{\gamma}\qquad \forall s,y, P\mbox{-a.s. } \omega.
\]
Assume also for each $t>0$ we have
%
%e6.8 #&#
\begin{equation}
\label{ZL2bound}\sup_{s\le t,x\in\R
}E\bigl[Z(s,x)^2\bigr]<
\infty.
\end{equation}
Let $X$ be a continuous $C^+_{\mathrm{rap}}$-valued solution to the
following SPDE, perhaps on a different space,
%
%e6.9 #&#
\begin{equation}
\label{gammaspde2} \frac{\partial X}{\partial t}= \frac{1}{2}\Delta X+
aX^{\gamma
}\dot W,
\end{equation}
with $Z(0,\cdot)=X(0,\cdot)\in C^+_{\mathrm{rap}}$.
Let $A$ be a Borel set in $ \R_+\times\R$. Then
\[
P\bigl(\operatorname{supp}(Z)\cap A=\varnothing\bigr)\geq P_{X_0}
\bigl(\operatorname{supp}(X)\cap A=\varnothing\bigr).
\]
\end{prop}

We will apply this result with $Z(t,x)=\bar U^i(s_i+t,x)$. To ensure
\eqref{ZL2bound} we will need the following moment bound which will
also give Lemma~\ref{TRbnd}. It will be proved in Appendix~\ref{secspdegr}.
%
%le6.4 #&#
\begin{lemma}
\label{pmom}
For any $q,T>0$, there exists $C_{q,T}$ such that
%
%e6.10 #&#
\begin{equation}
\sup_{0<\ep\le1}E \Bigl[\sup_{s\le T,x\in\R}\bigl(\bar
U(s,x)^{q}+\bar V(s,x)^{q}\bigr) \Bigr]\le
C_{q,T}.
\end{equation}
\end{lemma}
%
%Recall from \citet{myt98w} and \citet{shi94} that $\bar U$ is the unique
%in law continuous $C^+_{\mathrm{rap}}$-valued solution of \eqref{eq28}.
The proof of the above lemma is based on a simple adaptation of the
methods used for the proof of Proposition~1.8(a) of~\citet{mps06}, and
in particular Lemma \ref{mp44} of that paper.

%pa6.subsection.subsubsection.3 #&#
\begin{pf*}{Proof of Lemma~\protect\ref{TRbnd}}
This result with
$\delta_{\scriptsize{\ref
{TRbnd}}}(t)=C_{{1}/{2}, 2}t^{\ep_0/2}$ is an immediate corollary of
Markov's lemma and the above lemma with $q=1/2$. %
\end{pf*}

\begin{pf*}{Proof of Lemma~\ref{rhobnd}} We first fix $1\le i\le
N_\ep
$ and
argue conditionally on~$\CF_{s_i}$. Note that the inequalities in
\eqref{barsqfnbnd} hold pointwise, that is, without integrating over space.
These inequalities together with \eqref{eq25}, Lemma~\ref{pmom} and
Proposition~\ref{thm11} show the hypotheses of Proposition~\ref
{prop1} hold
with $Z(t,x)=\bar U^i(s_i+t,x)$, $Z_0=J_\ep^{x_i}$ and
$a=2^{{1}/{2}-\gamma}$. We apply this result to the open set
\[
A=A_t=\bigl\{(s,y)\dvtx |y-x_i|>\ep^{1/2}+s^{(1/2)-\delta_0},
 0<s<t\bigr\}
\]
and conclude that if $\rho$ is as in Lemma~\ref{qMP}, then
\[
P(\rho_i<t)=P_{J_\ep^{x_i}}\bigl(\mbox{supp}(Z)\cap A\neq
\varnothing\bigr)\le P(\rho<t).
\]
Corollary~\ref{qMP2} now shows there is a $c_{\scriptsize{\ref{rhobnd}}}=
c_{\scriptsize{\ref
{rhobnd}}}(\gamma,\delta_0)$ so that for $\ep,t\in(0,1]$,
\[
P(\rho_i\le t)\le c_{\scriptsize{\ref{rhobnd}}}\ep(t\vee\ep).
\]
It follows that for $p,\ep,t\in(0,1]$,
\[
P \Biggl(\bigcup_{i=1}^{pN_\ep}\{
\rho_i\le t\} \Biggr)\le\sum_{i=1}^{\lfloor
pN_\ep\rfloor}P(
\rho_i\le t)\le c_{\scriptsize{\ref{rhobnd}}}\lfloor pN_\ep \rfloor\ep
(t\vee\ep)\le c_{\scriptsize{\ref{rhobnd}}}p(t\vee\ep)\1(p\ge\ep).
\]
This finishes the proof of Lemma~\ref{rhobnd}.
\end{pf*}

We next turn to the proof of Proposition~\ref{prop1}.
Recall from the discussion at the beginning of this section that for
each $X_0\in C^+_{\mathrm{rap}}$ there is a unique law $P_{X_0}$ on $C(\R
_+,C^+_{\mathrm{rap}})$ of the solution to \eqref{gammaspde2}. We assume the
hypotheses of Proposition~\ref{prop1} for the rest of this section.

%le6.5 #&#
\begin{lemma}
\label{lem1} Let $\gamma\in[1/2,1)$.
For any nonnegative $\phi\in L^1(\R)$, and $t,s\geq0$, there exists a
sequence of $M_F(\R)$-valued processes $\{Y^{n}\}_{n\geq0}$ such that
$Y^n_0(dx)=\phi(x)\,dx$ and
%
%e6.11 #&#
%e6.12 #&#
\begin{eqnarray}
\label{eqdual1} E \bigl[ e^{-\langle\phi, Z_t\rangle}|\CF^{Z}_s
\bigr]&\geq& E \bigl[ e^{-\langle\phi, X_{t-s}\rangle}|X_0=Z_s \bigr]
\\
\label{eqdual2} &=& \lim_{n\rightarrow\infty}E^{Y^n}_{\phi}
\bigl[ e^{-\langle
Y^n_{t-s}, Z_s\rangle} \bigr],
\end{eqnarray}
where $P^{Y^n}_{\phi}$ is the probability law of $Y^n$.
\end{lemma}
%
%pa6.subsection.subsubsection.4 #&#
\begin{pf}
We may assume without loss of generality that $a=1$,
as only trivial adjustments are needed to the handle general $a>0$.
First we will prove the lemma for $\gamma>1/2$ and then explain the
modifications for the $\gamma=1/2$ case. For $\gamma\in(1/2,1)$,
(\ref{eqdual2}) follows from Proposition~2.3 of~\citet{myt98w}. To
simplify the exposition let us take $s=0$. For $s>0$ the proof goes
along the same lines as it depends only on the martingale properties of $Z$.

By the proof of Lemma~3.3 in~\citet{myt98w} we get that
for each $n$ there exists a stopping time $\tgamma_{k}(t)\le t$ and an
$M_F(\R)$-valued process
$Y^n$ such that, for $\eta=\frac{2\gamma(2\gamma-1)}{\Gamma
(2-2\gamma
)}$, and
\[
g(u,y)=\int_0^u \bigl(e^{-\lambda y} -1 +
\lambda y\bigr) \lambda^{-2\gamma-1} \,d\lambda,\qquad u,y\geq0,
\]
we have
%
%e6.13 #&#
\begin{eqnarray}
\label{eq0508}
&& E \bigl[ e^{-\langle Y^n_{\tgamma_k(t)},
Z_{t-\tgamma
_k(t)}\rangle}|Y^n_0=
\phi \bigr]\nonumber\\
&&\qquad = E_{\phi} \bigl[ e^{-\langle\phi, Z_t\rangle} \bigr]\nonumber\\
&&\qquad\quad{}- \frac{1}{2}E \biggl[ \int_0^{\tgamma_k(t)}
e^{-\langle Y^n_{s},
Z_{t-s}\rangle} \biggl\{ \eta\int_{\R}
\bigl(Y^n_s(x)\bigr)^2 g\bigl(1/n,
Z_{t-s}(x)\bigr)\,dx
\nonumber
\\[-8pt]
\\[-8pt]
\nonumber
&&\hspace*{150pt}{}+ \bigl\langle\sigma(Z_{t-s})^2-(Z_{t-s})^{2\gamma},
\bigl(Y^n_s\bigr)^2\bigr\rangle \biggr\}
\,ds \biggr]
\\
&&\qquad\leq
\nonumber
E_{\phi} \bigl[ e^{-\langle\phi, Z_t\rangle} \bigr]
\\
\nonumber
&&\qquad\quad\mbox{}- \frac{1}{2}E \biggl[ \int_0^{\tgamma_k(t)}
e^{-\langle Y^n_{s},
Z_{t-s}\rangle} \eta\int_{\R} \bigl(Y^n_s(x)
\bigr)^2 g\bigl(1/n, Z_{t-s}(x)\bigr)\,dx \,ds \biggr].\nonumber
\end{eqnarray}
If $k=k_n=\ln(n)$, we can easily get [as in the proof of Lemma~3.4
of~\citet{myt98w}] that
%
%e6.14 #&#
\begin{eqnarray}
\label{dualconvgt}&& E \biggl[ \int_0^{\tgamma_{k_n}(t)}
e^{-\langle Y^n_{s}, Z_{t-s}\rangle} \eta\int_{\R} \bigl(Y^n_s(x)
\bigr)^2 g\bigl(1/n, Z_{t-s}(x)\bigr)\,dx \,ds \biggr]
\nonumber\\
&&\qquad\leq C \sup_{x, s\leq t} E \bigl[ Z_s(x)^2
\bigr]k_n n^{2\gamma
-2}
\\
\nonumber
&&\qquad\rightarrow 0\qquad \mbox{as } n\rightarrow\infty.
\end{eqnarray}
Here we used \eqref{ZL2bound} in the last line.
Moreover, as is shown in the proof of Lemma~3.5 of~\citet{myt98w}, we have
\[
P\bigl(\tgamma_{k_n}(t)<t\bigr)\rightarrow 0 \qquad\mbox{as } n\rightarrow
\infty,
\]
or equivalently,
\[
P\bigl(\tgamma_{k_n}(t)=t\bigr)\rightarrow 1 \qquad\mbox{as } n\rightarrow
\infty.
\]
Hence we get from \eqref{eq0508}, \eqref{dualconvgt} and the above
\begin{eqnarray*}
&&
\lim_{n\rightarrow\infty} E \bigl[ e^{-\langle Y^n_{t},
Z_{0}\rangle
}|Y^n_0=
\phi \bigr] \nonumber\\
&&\qquad= \lim_{n\rightarrow\infty} E \bigl[ e^{-\langle Y^n_{\tgamma
_{k_n}(t)}, Z_{t-\tgamma_{k_n}(t)}\rangle}|Y^n_0=
\phi \bigr]
\\
\nonumber
&&\qquad\leq E \bigl[ e^{-\langle\phi, Z_{t}\rangle} \bigr]\qquad \forall t\geq0.
\end{eqnarray*}
But by Lemma~3.5 of \citet{myt98w} we have
%
%e6.15 #&#
\begin{eqnarray}
\label{dualconv2} \lim_{n\rightarrow\infty} E \bigl[ e^{-\langle Y^n_{t},
Z_{0}\rangle
}|Y^n_0=
\phi \bigr] = E \bigl[ e^{-\langle\phi, X_{t}\rangle} \bigr]\qquad\forall t\geq0
\end{eqnarray}
and we are done for $\gamma\in(1/2,1)$.

The case $\gamma=1/2$ is even easier. Now $X$ is just a super-Brownian
motion. Now take $Y^n=Y$ for all $n$, where $Y$ is a solution to the
log-Laplace equation
\[
\frac{\partial Y_t}{\partial t}= \frac{1}{2}\Delta Y_t -
\frac{1}{2}(Y_t)^2,
\]
so that \eqref{dualconv2} is the standard exponential duality for
super-Brownian motion.
Then~\eqref{eq0508} follows with $\tilde{\gamma}_k(t)=t$, and
$\eta
=0$, and so the result follows
immediately for $\gamma=1/2$.
\end{pf}

%le6.6 #&#
\begin{lemma}
\label{lem2}
For any $k\geq1$ and $0\leq t_1<t_2<\cdots< t_k$ and $\phi_1,\ldots
,\phi_k\geq0$,
%
%e6.16 #&#
\begin{equation}
\label{eqdual3} E \bigl[ e^{-\sum_{i=1}^k \langle\phi_i,
Z_{t_i}\rangle} \bigr]\geq E \bigl[
e^{-\sum_{i=1}^k\langle\phi_i, X_{t_i}\rangle} \bigr].
\end{equation}
\end{lemma}
%
%pa6.subsection.subsubsection.5 #&#
\begin{pf}
The proof goes by induction. For $k=1$ it follows
from the previous lemma. Suppose the equality holds for $k-1$. Let us
check it for $k$:
%
%e6.17 #&#
\begin{eqnarray}\label{eq26}
&&
E \bigl[ e^{-\sum_{i=1}^k \langle\phi_i, Z_{t_i}\rangle} \bigr]
\nonumber\\
&&\qquad= E \bigl[ e^{-\sum_{i=1}^{k-1}\langle\phi_i, Z_{t_i}\rangle} E \bigl[ e^{-\langle\phi_k, Z_{t_k}\rangle} |\CF^Z_{t_{k-1}}
\bigr] \bigr]
\nonumber\\
&&\qquad\geq E \Bigl[ e^{-\sum_{i=1}^{k-1}\langle\phi_i,
Z_{t_i}\rangle} \lim
_{n\rightarrow\infty} E^{Y^n}_{\phi_k} \bigl[
e^{-\langle Y^n_{t_k-t_{k-1}},
Z_{t_{k-1}}\rangle} \bigr] \Bigr]
\\
\nonumber
&&\qquad= \lim_{n\rightarrow\infty} E^{Y^n}_{\phi_k}
\times E^{Z} \bigl[ e^{-\sum_{i=1}^{k-2}\langle\phi_i, Z_{t_i}\rangle-\langle\phi
_{k-1}+ Y^n_{t_k-t_{k-1}}, Z_{t_{k-1}}\rangle} \bigr]
\\
\nonumber
&&\qquad\geq \lim_{n\rightarrow\infty} E^{Y^n}_{\phi_k}
\times E^{X} \bigl[ e^{-\sum_{i=1}^{k-2}\langle\phi_i, X_{t_i}\rangle-\langle\phi
_{k-1}+ Y^n_{t_k-t_{k-1}}, X_{t_{k-1}}\rangle} \bigr],
\end{eqnarray}
where the inequality in~(\ref{eq26}) follows by Lemma~\ref{lem1},
and the last inequality follows by the induction
hypothesis. Now, for $\gamma\in(1/2,1)$, we use conditioning and
Proposition~2.3 in~\citet{myt98w}
to get
%
%e6.18 #&#
\begin{eqnarray}\label{eq0608}
\nonumber
&&\lim_{n\rightarrow\infty} E^{Y^n}_{\phi
_k}
\times E^{X} \bigl[ e^{-\sum_{i=1}^{k-2}\langle\phi_i, X_{t_i}\rangle-\langle\phi
_{k-1}+ Y^n_{t_k-t_{k-1}}, X_{t_{k-1}}\rangle} \bigr]
\\
&&\qquad= E \Bigl[ e^{-\sum_{i=1}^{k-1}\langle\phi_i, X_{t_i}\rangle} \lim_{n\rightarrow\infty}
E^{Y^n}_{\phi_k} \bigl[ e^{-\langle Y^n_{t_k-t_{k-1}},
X_{t_{k-1}}\rangle} \bigr] \Bigr]
\\
 &&\qquad= E \bigl[ e^{-\sum_{i=1}^{k}\langle\phi_i, X_{t_i}\rangle
} \bigr],\nonumber
\end{eqnarray}\eject\noindent
and we are done for $\gamma\in(1/2,1)$. For $\gamma=1/2$, \eqref
{eq0608} follows immediately again by conditioning,
and the fact that $Y=Y^n$ is a solution to the log-Laplace equation
for super-Brownian motion.
\end{pf}

%le6.7 #&#
\begin{lemma}
\label{lem3}
For any nonnegative and Borel measurable function $\psi$ on \mbox
{$\R
_+\times\R$}
% \in\B_+(\R_+\times\R)$ (i.e., $\psi$ is nonnegative and Borel)
%
%e6.19 #&#
\begin{equation}
\label{eqdual4} \qquad E \bigl[ e^{-\int_0^t\int_{\R} \psi(s,x) Z(s,x)\,dx\,ds} \bigr]\geq E \bigl[ e^{-\int_0^t\int_{\R} \psi(s,x) X(s,x)\,dx\,ds}
\bigr]\qquad \forall t\geq0.
\end{equation}
\end{lemma}

Before starting the proof, we recall the following definition.

%de6.8 #&#
\begin{definition}
We say that a sequence $\psi_n(x)$ of functions converges
bounded-pointwise to
$\psi(x)$ provided $\lim_{n\to\infty}\psi_n(x)=\psi(x)$ for all $x$,
and there
exists a constant $K<\infty$ such that $\sup_{n,x}|\psi_n(x)|\leq K$.
\end{definition}

%pa6.subsection.subsubsection.6 #&#
\begin{pf*}{Proof of Lemma~\protect\ref{lem3}}
First suppose that $\psi\in C_+(\R_+\times\R)$ is bounded.
Then let us choose an approximating sequence of bounded functions
$\phi^n_1,\ldots,\break  \phi^n_{k_n}\in C_+(\R_+)$ such that
\[
\sum_{i=1}^{k_n} \langle
\phi_i, f_{t_i}\rangle\rightarrow\int_0^t
\int_{\R} \psi(s,x) f(s,x) \,ds\,dx\qquad \forall t\geq0
\]
for any $f\in D(\R_+, C_+(\R))$.
In this way for bounded $\psi\in C_+(\R_+\times\R)$ the result follows
immediately from Lemma~\ref{lem2}. Now pass to the bounded-pointwise closure
of this class of $\psi$'s, that is the smallest class containing the above
continuous $\psi$'s which is closed under bounded-pointwise limits. Finally
take monotone increasing limits to complete the proof. %
\end{pf*}

%pa6.subsection.subsubsection.7 #&#
\begin{pf*}{Proof of Proposition~\protect\ref{prop1}}
Take
\[
\psi_n(s,x)= n 1_{A}(s,x).
\]
Then by Lemma~\ref{lem3} we have
\[
E \bigl[ e^{-n Z(A)} \bigr]\geq E \bigl[ e^{-n X(A)} \bigr],
\]
where $Z(A)\equiv\int_A Z(s,x)\,dx\,ds$ and $X(A)\equiv\int_A X(s,x)\,dx\,ds$.
Take $n\rightarrow\infty$ on both sides to get
%
%e6.20 #&#
\begin{equation}
\label{eq27} P\bigl(Z(A)=0\bigr)\geq P\bigl(X(A)=0\bigr).
\end{equation}
The required result follows immediately for $A$ open because then
\[
\bigl\{\operatorname{supp}(Z)\cap A=\varnothing\bigr\}=\bigl\{Z(A)=0\bigr\}.
\]
It then follows for compact $A$ because
\[
\bigl\{\operatorname{supp}(X)\cap A=\varnothing\bigr\}=\bigcup
_{n}\bigl\{\operatorname {supp}(X)\cap A^{1/n}=
\varnothing\bigr\},
\]
where $A^{1/n}$ is the open set of points distance less than $1/n$ of $A$.
The general result now follows by
the inner regularity of the Choquet capacity $A\to P(\operatorname{supp}(Z)\cap
A\neq\varnothing)$; see page 39 of \citet{mey66}.
\end{pf*}

%s7 #&#
\section{Bounds on the killing measure: Proof of Lemma~\texorpdfstring{\protect\ref
{thetabnd}}{4.3}}\label{secKgrowth}
\setcounter{equation}{0}
Let
\[
G\bigl(\bar U^i\bigr)=\overline{\bigl\{(t,x)\dvtx \bar
U^i(t,x)>0\bigr\}}
\]
be the closed graph of $\bar U^i$, and let
\[
\Gamma^U_i(t)=\Gamma^U_i(t,
\delta_0)=\bigl\{(s,x)\dvtx s_i\le s\le
s_i+t, |x-x_i|\le(s-s_i)^{(1/2)-\delta_0}+
\ep^{1/2}\bigr\},\vadjust{\goodbreak}
\]
and let $\Gamma^V_j(t)$ be the corresponding set for $V$ with
$(t_j,y_j)$ in place of $(s_i,x_i)$. It is easy to check, using the
definition of $\rho_i$,
% and continuity of $\bar U^i_s$ for $s\ge s_i$,
that
%
%e7.1 #&#
\begin{equation}
\label{Umod} G\bigl(\bar U^i\bigr)\cap\bigl([s_i,s_i+
\rho_i]\times \R \bigr)\subset\Gamma_i^U(
\rho_i).
\end{equation}
%
%HERE IS A SPOT WHERE USING THE JOINT CONTINUITY OF $\bar U^i$ ON $[s_i,
%INTERVAL IN THE ABOVE--BUT IT IS MINOR AND CURRENTLY NOT NEEDED.
Of course an analogous inclusion holds for $\bar V^j$. If $K'(\cdot)$
is a nondecreasing right-continuous $M_F(\R)$-valued process, we let
$S(K')$ denote the closed support of the associated random measure on
space--time, $K'(ds,dx)$.

%le7.1 #&#
\begin{lemma}\label{KinU} $S(K^{i,U})\subset G(\bar U^i)$ and
$S(K^{j,V})\subset G(\bar V^j)$ for all $i,j\in\N_\ep$, \mbox{$P$-a.s.}
\end{lemma}
%
%pa7.subsection.subsubsection.1 #&#
\begin{pf}
It is easy to see from \eqref{UVdefn} that
$S(K^{i,U})\subset[s_i,\infty)\times\R$. Let $\CO$ be a bounded open
rectangle in $((s_i,\infty)\times\R)\cap G(\bar U^i)^c$ whose corners
have rational coordinates, and choose a smooth nonnegative function
$\phi$ on $\R$ so that $\CO=(r_1,r_2)\times\{\phi>0\}$. Then $\bar
U_r^i(\phi)=0$ for all $r\in(r_1,r_2)$ and hence for all $r\in
[r_1,r_2]$ a.s. by continuity. It then follows from \eqref{UVdefn} and
$U^i\le\bar U^i$ that a.s.\looseness=-1
\[
0=U^i_{r_2}(\phi)-U^i_{r_1}(\phi)=-
\bigl(K^{i,U}_{r_2}(\phi )-K^{i,U}_{r_1}(
\phi)\bigr).
\]\looseness=0
Therefore $K^{i,U}(\CO)=0$. Taking unions over such open ``rational''
rectangles, we conclude that
\[
K^{i,U}\bigl(G\bigl(\bar U^i\bigr)^c\cap
\bigl((s_i,\infty)\times\R\bigr)\bigr)=0 \qquad  \mbox{a.s}.
\]
On the other hand, from \eqref{eq25},
\begin{eqnarray*}
K^{i,U}\bigl(G\bigl(\bar U^i\bigr)^c\cap\bigl(
\{s_i\}\times\R\bigr)\bigr)&\le& K^{i,U}\bigl(
\{s_i\} \times [x_i-\sqrt\ep,x_i+\sqrt
\ep]^c\bigr)
\\
&=&0.
\end{eqnarray*}
In the last line we used \eqref{UVdefn} (recall from Section~\ref{secsetup} this implies $U^i_s=0$ for $s<s_i$) to see that
$K^{i,U}_{s_i}(\cdot)\le\langle J^{x_i},\cdot\rangle$.
The last two displays imply that $K^{i,U}(G(\bar U^i)^c)=0$ and hence
the result for $K^{i,U}$. The proof for $K^{j,V}$ is the same.
\end{pf}

Next we need a bound on the extinction times of nonnegative
martingales which is a slight generalization of Lemma~3.4 of \citet{mp92}.

%le7.2 #&#
\begin{lemma}\label{lemMhit0} Assume $\gamma'=\gamma''=\frac
{1}{2}$ or
$(\gamma',\gamma'')\in(1/2,1)\times[1/2,1]$. Let $M\ge0$ be a
continuous $(\CH_t)$-local martingale and $T$ be an $(\CH_t)$-stopping
time so that for some $\delta\ge0$ and $c_0>0$,
%
%e7.2 #&#
\begin{equation}
\label{sqfnineq} \frac{d\langle M\rangle_t}{dt}\ge c_0\1(t<T)M_t^{2\gamma'}(t+
\delta )^{(1/2)-\gamma''} \qquad\mbox{for }t>0.
\end{equation}
If $\tau_M(0)=\inf\{t\ge0\dvtx M_t=0\}$, then there is a $c_{\scriptsize{\ref
{lemMhit0}}}(\gamma')>0$ such that
\[
P\bigl(T\wedge\tau_M(0)\ge t|\CH_0\bigr)\le
c_{\scriptsize{\ref{lemMhit0}}}\bigl(\gamma '\bigr)c_0^{-1}M_0^{2-2\gamma'}
t^{\gamma''-(3/2)}\qquad
\mbox{for all }t\ge \delta/2.
\]
\end{lemma}
%
%pa7.subsection.subsubsection.2 #&#
\begin{pf}
If $\gamma'=\gamma''=\frac{1}{2}$, the lemma follows
from a
slight extension of the proof of Lemma~\ref{hittime}, so assume
$\gamma'\in(1/2,1)$. Let $V=T\wedge\tau_M(0)$. As usual there is a Brownian
motion $B(t)$ such that $M(t)=B(\langle M\rangle_t)$ for $t\le V$. By
\eqref{sqfnineq} we have
\begin{eqnarray*}
\int_0^V c_0(t+
\delta)^{(1/2)-\gamma''} \,dt&\le&\int_0^V
M_t^{-2\gamma'} \,d\langle M\rangle_t
\\[-2pt]
&\le&\int_0^{\langle M\rangle_V} B_u^{-2\gamma'}
\,du \le\int_0^{\tau_B(0)} B_u^{-2\gamma'}
\,du.
\end{eqnarray*}
If $L^x_t, x\in\R,  t\ge0$ is the semimartingale local time of $B$,
the Ray--Knight theorem [see Theorem~VI.52.1 in \citet{rw87}] and the
occupation time formula implies that the above gives
%
%e7.3 #&#
\begin{eqnarray}\label{Vlowerbnd}
&&
\nonumber
E \bigl[(V+\delta)^{(3/2)-\gamma''}-\delta^{(3/2)-\gamma
''} |
\CH_0 \bigr]
\\[-2pt]
\nonumber
&&\qquad\le \bigl((3/2)-\gamma''
\bigr)c_0^{-1}\int_0^\infty
x^{-2\gamma
'}E\bigl(L^x_{\tau_B(0)}|B_0\bigr)
\,dx
\nonumber
\\[-9pt]
\\[-9pt]
\nonumber
&&\qquad=\bigl((3/2)-\gamma''
\bigr)c_0^{-1}\int_0^\infty
x^{-2\gamma
'}2(M_0\wedge x) \,dx
\\[-2pt]
&&\qquad\le c_1\bigl(\gamma'
\bigr)c_0^{-1}M_0^{2-2\gamma'} \qquad\mbox{
(use $\gamma'>1/2$)}.\nonumber
\end{eqnarray}
A bit of calculus shows that
%
%e7.4 #&#
\begin{equation}
\label{calc}\qquad (t+\delta)^{(3/2)-\gamma''}-\delta^{(3/2)-\gamma''}\ge
\tfrac
{1}{2}(\sqrt3-\sqrt2)t^{(3/2)-\gamma''}\qquad \mbox{for all }t\ge \delta/2.
\end{equation}
Therefore by \eqref{Vlowerbnd} and \eqref{calc}, for $t\ge\delta/2$,
\begin{eqnarray*}
P(V\ge t|\CH_0)&\le& \frac{E [(V+\delta)^{(3/2)-\gamma
''}-\delta
^{(3/2)-\gamma''} |\CH_0 ]}{(t+\delta)^{(3/2)-\gamma
''}-\delta
^{(3/2)-\gamma''}}
\\[-2pt]
&\le& \frac{2c_1(\gamma')c_0^{-1}M_0^{2-2\gamma'}}{(\sqrt3-\sqrt
2)t^{(3/2)-\gamma''}}
\\[-2pt]
&\equiv& c_{\scriptsize{\ref{lemMhit0}}}c_0^{-1}M_0^{2-2\gamma'}t^{\gamma
''-(3/2)}.%\qed
\end{eqnarray*}
\upqed\end{pf}\eject

Define $\rho^V_j=\rho_j^{V,\delta_0,\ep}$ just as $\rho_i$ but with
$\bar V^j_{t_j+t}$ in place of $\bar U^i_{s_i+t}$ and $y_j$ in place of
$x_i$.\vspace*{-2pt}

%le7.3 #&#
\begin{lemma}\label{rhoVbnd} $Q_i (\bigcup_{j=1}^{pN_\ep}\{\rho
_j^V\le
t\} )\le c_{\scriptsize{\ref{rhobnd}}}(t\vee\ep)p\1(p\ge\ep)\mbox{ for
all }\ep
,p,t\in(0,1]$ and $i\in\N_\ep$.\vspace*{-2pt}
\end{lemma}
%
%pa7.subsection.subsubsection.3 #&#
\begin{pf}
All the $P$-local martingales and $P$-white noises
arising in the definition\vadjust{\goodbreak} of $\{\bar V^j,j\in\N_\ep\}$ remain such
under $Q_i$ because they are all orthogonal to
\[
\frac{dQ_i}{dP} \Big|_{\CF_t}=\1(t<s_i)+\1(t\ge
s_i)\frac{\bar
U^i_{t\wedge(s_i+\bar\tau_i)}(1)}{\ep}.
\]
The proof of Lemma~\ref{rhobnd} for $\{\rho_i\}$ under $P$ therefore
applies to $\{\rho_j^V\}$ under~$Q_i$.
\end{pf}

Recall we are trying to show that the killing measure $K^{i,U}_t$
associated with the~$i$ cluster of $U$ grows slowly enough for small
$t$. We will control the amount of killing here by controlling the
amount of killing by the $V^j$'s. The following result essentially
shows that with high probability for small $t$, there is no killing
during $[s_i,s_i+t]$ from the $V^j$'s which are born before time
$s_i$. Note it is particularly important that there is no $V$ mass on
the birth site of the $U^i$ cluster.

Recall from \eqref{bardef} that $\bar\delta=\bar\delta(\gamma
)=\frac
{1}{3} (\frac{3}{2}-2\gamma )$. We introduce
\[
\underline\rho_i^V=\min_{j:t_j\le s_i}
\rho_j^V.
\]

%le7.4 #&#
\begin{lemma}\label{presimass} There is a constant $c_{\scriptsize{\ref
{presimass}}}(\gamma)>0$ so that for $0<\delta_0\le\bar\delta
(\gamma)$,
\begin{eqnarray}
Q_i \biggl(\Gamma_i^U(t)\cap \biggl\{
\bigcup_{j:t_j\le s_i}G\bigl(\bar V^j\bigr) \biggr
\}\neq \varnothing, \underline\rho^V_i>2t \biggr)\le
c_{\scriptsize{\ref
{presimass}}}(\gamma ) (\ep\vee t)^{\bar\delta}
\nonumber
\\[-2pt]
\eqntext{\mbox{for all }\ep,t\in(0,1]\mbox{ and }s_i\le t.}
\end{eqnarray}
\end{lemma}
%
%pa7.subsection.subsubsection.4 #&#
\begin{pf}
Assume $\ep,t,s_i$ and $\delta_0$ are as above. Set $\alpha=\frac
{1}{2}-\delta_0(\ge\frac{1}{3})$ and choose $n_0\le n_1\in\Z_+$ so that
%
%e7.5 #&#
\begin{equation}
\label{n0cond} 2^{-n_0-1}<t\vee\ep\le2^{-n_0},\qquad
2^{-n_1-1}<\ep\le2^{-n_1}.
\end{equation}

Assume that
%
%e7.6 #&#
\begin{equation}
\label{rhoVbig}\underline\rho_i^V>2t,
\end{equation}
until otherwise indicated. Suppose $t_j\le s_i$ (hence $t_j<s_i$)
and
\[
(t_j,y_j)\notin[0,s_i)\times
\bigl[x_i-7\cdot2^{-n_0\alpha},x_i+7\cdot
2^{-n_0\alpha}\bigr].
\]
Then
\[
|y_j-x_i|>7\cdot2^{-n_0\alpha}\ge7(t\vee
\ep)^\alpha\ge t^\alpha +(t+s_i-t_j)^\alpha+2
\sqrt\ep,
\]
and so
\[
\Gamma_i^U(t)\cap\Gamma_j^V(s_i+t-t_j)=
\varnothing.
\]
By \eqref{rhoVbig} we have $\rho^V_j>s_i-t_j+t$, and so by \eqref
{Umod}, or more precisely its analogue for $\bar V^j$, we have
%
%e7.7 #&#
\begin{equation}
\label{possVj} \Gamma_i^U(t)\cap G\bigl(\bar
V^j\bigr)\subset\Gamma_i^U(t)\cap\Gamma
_j^V(s_i+t-t_j)=\varnothing.
\end{equation}
We therefore have shown that, assuming \eqref{rhoVbig},
%
%e7.8 #&#
\begin{eqnarray}
\label{possj} &&\bigl\{(t_j,y_j)\dvtx t_j\le
s_i, \Gamma_i^U(t)\cap G\bigl(\bar
V^j\bigr)\neq \varnothing \bigr\}
\nonumber
\\[-8pt]
\\[-8pt]
\nonumber
&&\qquad \subset[0,s_i)\times
\bigl[x_i-7\cdot2^{-n_0\alpha},x_i+7\cdot
2^{-n_0\alpha}\bigr].
\end{eqnarray}

Next we cover the rectangle on the right-hand side of the above by
rectangles as follows:
\begin{eqnarray*}
R_n^0&=&\bigl[s_i-2^{-n+1},s_i-2^{-n}
\bigr]\times\bigl[x_i-7\cdot2^{-n\alpha
},x_i+7
\cdot2^{-n\alpha}\bigr],
\\
R_n^r&=&\bigl[s_i-2^{-n},s_i
\bigr]\times\bigl[x_i+7\cdot2^{-(n+1)\alpha},x_i+7\cdot
2^{-n\alpha}\bigr],
\\
R_n^\ell&=&\bigl[s_i-2^{-n},s_i
\bigr]\times\bigl[x_i-7\cdot2^{-n\alpha},x_i-7\cdot
2^{-(n+1)\alpha}\bigr].
\end{eqnarray*}
Then it is easy to check that\vspace*{-1pt}
%
%e7.9 #&#
%e7.10 #&#
\begin{eqnarray}
\label{Rcontain1}
&&
\bigcup_{n=n_0}^\infty
\bigl(R^0_n\cup R^r_n\cup
R_n^\ell \bigr)
\nonumber
\\[-9pt]
\\[-9pt]
\nonumber
&&\qquad\supset \bigl[s_i-2^{-n_0+1},s_i\bigr)
\times\bigl[x_i-7\cdot2^{-n_0\alpha},x_i+7\cdot
2^{-n_0\alpha}\bigr]
\\[-1pt]
\label{Rcontain}&&\qquad\supset [0,s_i)\times\bigl[x_i-7
\cdot2^{-n_0\alpha
},x_i+7\cdot2^{n_0\alpha}\bigr].
\end{eqnarray}
We group together those $\bar V^j$'s which have their initial ``seeds''
in each of the above rectangles. That is, for $q=0,\ell,r$ consider
\begin{eqnarray*}
V^{n,q}(t,x)&=&\sum_j\1
\bigl((t_j,y_j)\in R_n^q
\bigr)V^j(t,x),
\\[-1pt]
{\widetilde V}^{n,q}(t,x)&=&\sum_j\1
\bigl((t_j,y_j)\in R_n^q\bigr){
\widetilde V}^j(t,x),
\\[-1pt]
{\bar V}^{n,q}(t,x)&=&\sum_j\1
\bigl((t_j,y_j)\in R_n^q\bigr){
\bar V}^j(t,x).
\end{eqnarray*}
We also let $V_t^{n,q}$, $\widetilde V_t^{n,q}$ and $\bar V_t^{n,q}$ denote
the corresponding measure-valued processes.

It follows from \eqref{possj} and \eqref{Rcontain} that\vspace*{-1pt}
%
%e7.11 #&#
\begin{eqnarray}
\label{HitBound0} && Q_i\biggl(\bigcup_{t_j\le s_i}
\bigl(G\bigl(\bar V^j\bigr)\cap\Gamma_i^U(t)
\neq\varnothing, \underline\rho_i^V>2t\bigr)\biggr)
\nonumber\\[-1pt]
&&\qquad\le \sum_{n=n_0}^{n_1}\sum
_{q=0,r,\ell}Q_i\bigl(G\bigl(\bar
V^{n,q}\bigr)\cap\Gamma_i^U(t)\neq
\varnothing, \underline\rho _i^V>2t\bigr)
\\[-1pt]
\nonumber
&&\qquad\quad{} + Q_i\Biggl(\bigcup_{n=n_1+1}^\infty
\bigcup_{q=0,r,l} \bigl(G\bigl(\bar V^{n,q}
\bigr)\cap\Gamma_i^U(t)\bigr)\neq\varnothing\Biggr).
\end{eqnarray}
We will use different arguments to show that each of the two terms on
the right-hand side of~\eqref{HitBound0} is small.
For the second term a very crude argument works. Namely, for
the supports of the $\bar V^j$ clusters with initial ``seeds'' in
$\bigcup_{n=n_1+1}^\infty (R_n^0\cup R_n^r\cup R_n^\ell)$ to intersect the
support of $U^i$, the $\bar V^j$ clusters must be born in $\bigcup_{n=n_1+1}^\infty (R_n^0\cup R_n^r\cup R_n^\ell)$, and the probability
of this event
is already small. More precisely,
%
%e7.12 #&#
\begin{eqnarray}
\label{HitBound01} && Q_i\Biggl(\bigcup
_{n=n_1+1}^\infty\bigcup_{q=0,r,l}
\bigl(G\bigl(\bar V^{n,q}\bigr)\cap\Gamma_i^U(t)
\bigr)\neq\varnothing\Biggr)
\nonumber
\\[-8pt]
\\[-8pt]
\nonumber
&&\qquad\le Q_i \Biggl(\eta_\ep^- \Biggl(\bigcup
_{n=n_1+1}^\infty \bigl(R_n^0
\cup R_n^r\cup R_n^\ell\bigr)
\Biggr)>0 \Biggr).
\end{eqnarray}
%
%LM It is easy to check that
By Proposition~\ref{girs} and
the decomposition for $\bar U^i(1)$ in \eqref{eq25} [see also~\eqref
{barmassmart}], we have
%
%e7.13 #&#
\begin{eqnarray}
\label{Qiunif} && Q_i\bigl((x_i,y_j)
\in A\bigr)
\nonumber
\\[-8pt]
\\[-8pt]
\nonumber
&&\qquad=E_P \biggl(\frac{\bar U^i_{s_i+[(t_j-s_i)^+\wedge\bar
\tau
_i]}(1)}{\ep}\1\bigl((x_i,y_j)
\in A\bigr) \biggr)=P\bigl((x_i,y_j)\in A\bigr).
%   A\subset[0,1].
\end{eqnarray}
%
% by using the decomposition for $\bar U^i(1)$ in \eqref{eq25}.
This and the analogue of \eqref{Rcontain1} with $n_1+1$ in place of
$n_0$, implies that the
right-hand side
of~\eqref{HitBound01} is at most
%
%e7.14 #&#
\begin{eqnarray}
\label{tailbnd}
&& Q_i \bigl(\eta_\ep^-
\bigl(\bigl[s_i-2^{-n_1},s_i\bigr)\times
\bigl[x_i-7\cdot 2^{-(n_1+1)\alpha},x_i+7
\cdot2^{-(n_1+1)\alpha}\bigr]\bigr)>0 \bigr)
\nonumber
\\[-8pt]
\\[-8pt]
\nonumber
&&\qquad\le 2\bigl(14\cdot2^{-(n_1+1)\alpha}\bigr)\le42\ep^\alpha.
\end{eqnarray}
Substitute this bound into \eqref{HitBound0} to get
%
%e7.15 #&#
\begin{eqnarray}
\label{HitBound1} && Q_i \biggl(\bigcup
_{t_j\le s_i}\bigl(G\bigl(\bar V^j\bigr)\cup
\Gamma_i^U(t)\bigr)\neq\varnothing, \underline
\rho_i^V>2t \biggr)
\nonumber
\\[-8pt]
\\[-8pt]
\nonumber
&&\qquad\le \sum_{n=n_0}^{n_1}\sum
_{q=0,r,\ell}Q_i\bigl(G\bigl(\bar
V^{n,q}\bigr)\cap\Gamma_i^U(t)\neq
\varnothing, \underline\rho_i^V>2t\bigr) %\\
% (R_n^0\cup R_n^r\cup R_n^\ell)\Bigr)>0\Bigr).
+ 42\ep^\alpha.
\end{eqnarray}

Now we are going to bound each term in the sum on the right-hand side
of~\eqref{HitBound1}.
To this end, in what follows, we assume that $n_0\leq n\leq n_1$, and,
for $q=0,r,l$, set
%
%e7.16 #&#
\begin{eqnarray}
\label{NDef}&& N_t^{n,q}=\sum
_j \mathbf{1}\bigl((t_j,y_j)\in
R_n^q\bigr)\nonumber\\
&&\hspace*{48pt}{}\times \int_0^t
\int_{\R
} \biggl(V(s,x)^{2\gamma-1}V^j(s,x)
\\
&&\hspace*{93pt}{}+\bigl(\bar V(s,x)^{2\gamma}-V(s,x)^{2\gamma}\bigr)
\frac{\widetilde
V^j(s,x)}{\widetilde V(s,x)} \biggr)^{1/2}\bar W^{j,V}(ds,dx).\nonumber
\end{eqnarray}
Note that $N^{n,q}$ is a continuous local martingale under $Q_i$.

The treatment of the cases $q=0$ and $q=r,l$ is different. First, let
$q=0$. Basically, in this case, we will show that, the on the event $\{
\underline\rho^V_i>2t\}$, the total
mass of $\bar V^{n,0}$ dies out with high probability before the time
$s_i$ (and, in fact, even before $s_i-2^{-n-1}$).
Hence, with this high probability, the support of $\bar V^{n,0}$ does
not intersect $\Gamma^U_i$. Let us make this precise.
We have from \eqref{eq25}\looseness=1
%
%e7.17 #&#
\begin{equation}
\label{barvdecomp} \bar V^{n,0}_{t+(s_i-2^{-n})^+}(1)=\bar
V^{n,0}_{(s_i-2^{-n})^+}(1)+\bar M^{n,0}_t,
\end{equation}\looseness=0
where
\[
\bar V^{n,0}_{(s_i-2^{-n})^+}(1)=\int\int\mathbf{1}\bigl((s,y)\in
R_n^0\bigr)\eta _\ep^-(ds,dy)+N^{n,0}_{(s_i-2^{-n})^+}
\]
and
%
%e7.18 #&#
\begin{equation}
\label{MDef} \bar M_t^{n,0}=N^{n,0}_{t+(s_i-2^{-n})^+}-N^{n,0}_{(s_i-2^{-n})^+}
\end{equation}
is a continuous $\CF_{t+(s_i-2^{-n})^+}$-local martingale under $Q_i$.

Assume for now that $s_i>2^{-n}$ since otherwise $\bar
V^{n,0}_{s_i}(1)=0$ and the bound \eqref{barVn0} below is trivial. An
easy localization argument shows that (recall that $n_0\le n\le n_1$)
%
%e7.19 #&#
\begin{eqnarray}
\label{initmassbnd} && Q_i \bigl(\bar V^{n,0}_{(s_i-2^{-n})}(1)
\ge2^{-n(1+\alpha-\bar\delta
)} \bigr)
\nonumber\\
&&\qquad\le 2^{n(1+\alpha-\bar\delta)}Q_i \biggl(\int\int\1 \bigl((s,y)\in
R_n^0\bigr)\eta^-_\ep(ds,dy) \biggr)
\nonumber
\\[-8pt]
\\[-8pt]
\nonumber
&&\qquad\le 2^{n(1+\alpha-\bar\delta)}\ep\bigl[\ep ^{-1}2^{-n}+1
\bigr]14\cdot 2^{-n\alpha}\qquad\bigl[\mbox{by }\eqref{Qiunif}\bigr]
\\
\nonumber
&&\qquad\le 14\bigl(2^{-n\bar\delta}\bigr) \bigl(2^n\ep+1\bigr)
\le28\cdot2^{-n\bar
\delta}.
\end{eqnarray}\eject

Now from \eqref{NDef} and \eqref{MDef}, if $t'\equiv s_i-2^{-n}+t<
T'\equiv\min_{j:t_j\le s_i}(\rho^V_j+t_j)$,
% and $n_0\le n\le n_1$,
then
%
%e7.20 #&#
\begin{eqnarray}
\label{sqfnM0bnd}
&&\frac{d}{dt}\bigl\langle\bar M^{n,0}\bigr
\rangle_t
\nonumber\\
&&\qquad=\int V\bigl(t',x\bigr)^{2\gamma-1}\bar
V^{n,0}\bigl(t',x\bigr)+\bigl(\bar V\bigl(t',x
\bigr)^{2\gamma}-V\bigl(t',x\bigr)^{2\gamma}\bigr)
\frac{\widetilde V^{n,0}(t',x)}{\tV
(t',x)} \,dx\nonumber
\\
&&\qquad\ge \int V^{n,0}\bigl(t',x
\bigr)^{2\gamma}+\tV ^{n,0}\bigl(t',x
\bigr)^{2\gamma} \,dx
\\
\nonumber
&&\qquad\ge 2^{-2\gamma}\int\bar V^{n,0}\bigl(t',x
\bigr)^{2\gamma}\1 \bigl(|x-x_i|\le7\cdot2^{-n\alpha}+
\bigl(2^{-n}+t\bigr)^\alpha+\sqrt\ep\bigr) \,dx
\\
\nonumber
&&\qquad\ge 2^{-2\gamma}\bar V^{n,0}_{t'}(1)^{2\gamma}
\bigl(2\bigl[7\cdot 2^{-n\alpha}+\bigl(2^{-n}+t
\bigr)^\alpha+\sqrt\ep\bigr]\bigr)^{1-2\gamma}.
\end{eqnarray}
In the last line we used Jensen's inequality and the fact that $T'> t'$
implies $\bar V^{n,0}(t',\cdot)$ is supported in the closed interval
with endpoints
$x_i\pm(7\cdot2^{-n\alpha}+(t+2^{-n})^\alpha+\sqrt\ep)$. A bit of
arithmetic\vadjust{\goodbreak} (recall $2^{-n}\ge\ep$ for $n\le n_1$) shows that \eqref
{sqfnM0bnd} implies for some
$c(\gamma)>0$,
%
%e7.21 #&#
\begin{eqnarray}
\label{sqfnM0bnd2} \frac{d}{dt}\bigl\langle\bar M^{n,0}\bigr
\rangle_t\ge c(\gamma) \bigl(\bar V^{n,0}_{t+(s_i-2^{-n})}(1)
\bigr)^{2\gamma}\bigl[2^{-n}+t\bigr]^{\alpha
(1-2\gamma
)}
\nonumber
\\[-8pt]
\\[-8pt]
\eqntext{\mbox{for }t<T\equiv \Bigl(\min_{j:t_j\le s_i}\bigl(\rho
^V_j+t_j\bigr)-\bigl(s_i-2^{-n}
\bigr) \Bigr)^+.} %\mbox{ and }n_0\le n\le n_1.
\end{eqnarray}
Note that $T$ is an $\CF_{(s_i-2^{-n})+t}$-stopping time. Therefore
\eqref{sqfnM0bnd2} allows us to apply Lemma~\ref{lemMhit0} to $t\to
\bar V^{n,0}_{(s_i-2^{-n})+t}(1)\equiv M_t$ with $\gamma'=\gamma$,
$\gamma''=\gamma-\delta_0(2\gamma-1)$ and $\delta=2^{-n}$. Here notice
that $\delta_0\le1/6$ implies $\gamma''\in[\frac{1}{2},\frac{3}{4}]$
and $\gamma''=1/2$ if $\gamma=1/2$. Therefore,
%for $n_0\le n\le n_1$, the Lemma
Lemma~\ref{lemMhit0}, the fact that $\underline\rho^V_i>2t$ implies
$T>t\ge s_i>2^{-n}$, and \eqref{initmassbnd} imply
%
%e7.22 #&#
\begin{eqnarray}
\label{barVn0}
&& Q_i\bigl(\bar V^{n,0}_{s_i-2^{-n-1}}(1)>0,
 \underline\rho^V_i>2t\bigr)
\nonumber\hspace*{-25pt}\\
&&\qquad\le Q_i\bigl(\bar V^{n,0}_{s_i-2^{-n}}(1)
\ge2^{-n(1+\alpha
-\bar
\delta)}\bigr)
\nonumber\hspace*{-25pt}\\
\nonumber
&&\qquad\quad{} +E_{Q_i} \bigl[Q_i\bigl(T\wedge
\tau_M(0)\ge2^{-n-1}|\CF _{s_i-2^{-n}}\bigr)\1\bigl(\bar
V^{n,0}_{s_i-2^{-n}}(1)<2^{-n(1+\alpha-\bar
\delta
)}\bigr) \bigr]\hspace*{-25pt}
\\
&&\qquad\le 28\cdot2^{-n\bar\delta}+c_{\scriptsize{\ref{lemMhit0}}}(\gamma )c(
\gamma)^{-1}2^{-n(1+\alpha-\bar\delta)(2-2\gamma
)}2^{-(n+1)(\gamma
-\delta_0(2\gamma-1)-(3/2))}\hspace*{-25pt}
\\
\nonumber
&&\qquad\le c'(\gamma) \bigl(2^{-n\bar\delta
}+2^{-n((3/2)-2\gamma
-2(1-\gamma)\bar\delta-\delta_0)}
\bigr) \qquad(\mbox{by the definition of }\alpha)\hspace*{-25pt}
\\
\nonumber
&&\qquad\le c'(\gamma) \bigl(2^{-n\bar\delta}+2^{-n(3\bar
\delta
-2(1-\gamma)\bar\delta-\delta_0)}
\bigr) \qquad(\mbox{by the definition of }\bar\delta)\hspace*{-25pt}
\\
\nonumber
&&\qquad\le c_0(\gamma)2^{-n\bar\delta},\hspace*{-25pt}
\end{eqnarray}
where $\delta_0\le\bar\delta$ and $\gamma\ge1/2$ are used in the
last line.

Next consider $\bar V^{n,r}$. The analogue of \eqref{barvdecomp} now is
\[
\bar V^{n,r}_{s_i+t}(1)=\bar V^{n,r}_{s_i}(1)+
\bar M^{n,r}_t,
\]
where
\[
\bar M^{n,r}_t=N^{n,r}_{s_i+t}-N^{n,r}_{s_i}.
\]
An argument similar to the derivation of \eqref{initmassbnd} shows that
% for $n_0\le n\le n_1$,
%
%e7.23 #&#
\begin{equation}
\label{initrmassbnd} Q_i\bigl(\bar V^{n,r}_{s_i}(1)
\ge2^{-n(1+\alpha-\bar\delta)}\bigr)\le28\cdot 2^{-n\bar\delta}.
\end{equation}
Next we argue as in \eqref{sqfnM0bnd} and \eqref{sqfnM0bnd2} to see
that for $s_i+t<T'\equiv\break \min_{j:t_j\le s_i}(\rho^V_j+t_j)$,
% and $n_0\le n\le n_1$,
\begin{eqnarray*}
\frac{d}{dt}\bigl\langle\bar M^{n,r}\bigr\rangle_t&
\ge& 2^{-2\gamma}\bar V^{n,r}_{s_i+t}(1)^{2\gamma}
\bigl(\bigl[7\cdot2^{-(n+1)\alpha
}+\bigl(2^{-n}+t\bigr)^{(1/2)-\delta_0}+
\sqrt\ep\bigr]2 \bigr)^{1-2\gamma}
\\
&\ge& c'(\gamma) \bigl(\bar V^{n,r}_{s_i+t}(1)
\bigr)^{2\gamma
}\bigl(2^{-n}+t\bigr)^{\alpha(1-2\gamma)},
\end{eqnarray*}
where we again used $n_0\le n\le n_1$.
Now we apply Lemma~\ref{lemMhit0} and \eqref{initrmassbnd}, as in the
derivation of \eqref{barVn0}, to conclude that
% for $n_0\le n\le n_1$,
%
%e7.24 #&#
\begin{equation}
\label{barVnr} Q_i\bigl(\bar V^{n,r}_{s_i+2^{-n}}(1)>0,
 \underline\rho^V_i>2t\bigr)\le c_1(
\gamma)2^{-n\bar\delta}.\vadjust{\goodbreak}
\end{equation}
If $\bar V_{s_i+2^{-n}}^{n,r}(1)=0$, then $\bar V^{n,r}_u(1)=0$ for all
$u\ge s_i+2^{-n}$, and so if in addition, $\underline\rho^V_i>2t$, then
by the definition of $\rho^V_j$,
%
%e7.25 #&#
\begin{eqnarray}
\label{Vnrinclus}
\nonumber
G\bigl(\bar V^{n,r}\bigr)&\subset& \bigl
\{(s,x)\dvtx s_i-2^{-n}\le s\le s_i+2^{-n},
\\
&& \hspace*{ 7pt}7\cdot2^{-(n+1)\alpha}-\bigl(s-s_i+2^{-n}
\bigr)^\alpha-\sqrt\ep
\\
&&\hspace*{7pt} \le x-x_i\le7\cdot2^{-n\alpha
}+
\bigl(s-s_i+2^{-n}\bigr)^\alpha+\sqrt\ep \bigr\}.\nonumber
\end{eqnarray}
A bit of algebra (using our choice of the factor $7$ and $n_0\le n\le
n_1$) shows that
\[
x_i+2^{-n\alpha}+\sqrt\ep<x_i+7
\cdot2^{-(n+1)\alpha
}-\bigl(2^{-n}+2^{-n}
\bigr)^\alpha-\sqrt\ep,
\]
and so the set on the right-hand side of \eqref{Vnrinclus} is disjoint
from $\Gamma_i^U(t)$. Therefore by~\eqref{barVnr} we may conclude that
%for $n_0\le n\le n_1$,
%
%e7.26 #&#
\begin{equation}
\label{barVnr2} Q_i\bigl(G\bigl(\bar V^{n,r}\bigr)\cap
\Gamma_i^U(t)\neq\varnothing, \underline\rho
^V_i>2t\bigr)\le c_1(\gamma)2^{-n\bar\delta}.
\end{equation}
Of course the same bound holds for $G(\bar V^{n,\ell})$.

Note that $\bar V^{n,0}_{s_i-2^{-n-1}}(1)=0$ implies $\bar
V^{n,0}_s(1)=0$ for all $s\ge s_i-2^{-n-1}$ and so $G(\bar V^{n,0})\cap
\Gamma_i^U(t)$ is empty. Therefore\vadjust{\goodbreak} \eqref{barVn0} and \eqref{barVnr2}
show that the summation on the right-hand side of \eqref{HitBound1} is
at most
\[
\sum_{n=n_0}^{n_1}\bigl(c_0(
\gamma)+2c_1(\gamma)\bigr)2^{-n\bar\delta}\le c_2(\gamma)
(t\vee\ep)^{\bar\delta}.
\]
We substitute the above
%and \eqref{tailbnd}
into \eqref{HitBound1} to see that
\begin{eqnarray*}
&& Q_i\biggl(\bigcup_{t_j\le s_i}\bigl(G\bigl(
\bar V^j\bigr)\cup\Gamma_i^U(t)\bigr)\neq
\varnothing, \underline\rho_i^V>2t\biggr)
\\
&&\qquad \le42\ep^\alpha+c_2(\gamma) (t\vee\ep)^{\bar\delta}\le
c_{\scriptsize{\ref
{presimass}}}(\gamma) (t\vee\ep)^{\bar\delta}.
\end{eqnarray*}
In the last line we used $\bar\delta\le1/6<1/4\le\alpha$.
\end{pf}
%pa7.subsection.subsubsection.5 #&#
\begin{pf*}{Proof of Lemma~\ref{thetabnd}} Fix $0<\delta
_0\le\bar\delta
$, $t\in(0,1]$ and assume $s_i,s\le t$. By~\eqref{Umod} and
Lemma~\ref
{KinU} on $\{\rho_i>s\}$ we have
\[
K^{i,U}_{s_i+s}(1)=K^{i,U}\bigl(\Gamma_i^U(s)
\bigr)\le\sum_j K^{j,V}\bigl(
\Gamma_i^U(s)\bigr),
\]
where \eqref{eq22} is used in the last inequality. Next use
$S(K^{j,V})\subset G(\bar V^j)$ (by Lemma~\ref{KinU}) and
$S(K^{j,V})\subset[t_j,\infty)\times\R$ to conclude that on
\[
\{\rho_i>s\}\cap \biggl\{ \biggl(\bigcup
_{t_j\le s_i}G\bigl(\bar V^j\bigr) \biggr)\cap
\Gamma_i^U(t)=\varnothing \biggr\}\equiv\{
\rho_i>s\}\cap D_i(t),
\]
we have\vspace*{-1pt}
%
%e7.27 #&#
\begin{equation}
\label{Kimassbnd1} K^{i,U}_{s_i+s}(1)\le\sum
_j \1(s_i<t_j\le
s_i+s)K^{j,V}\bigl(\Gamma_i^U(s)
\bigr).\vspace*{-1pt}
\end{equation}

Another application of \eqref{Umod} and Lemma~\ref{KinU}, this time to
$\bar V^j$, shows that for $t_j>s_i$,\vspace*{-1pt}
%
%e7.28 #&#
\begin{equation}
\label{SKjVinc} S\bigl(K^{j,V}\bigr)\cap\bigl([0,s_i+s]
\times\R\bigr)\subset\Gamma_j^V(s_i+s-t_j)
\qquad\mbox{on }\bigl\{\rho_j^V>s\bigr\}.\vspace*{-1pt}
\end{equation}
An elementary calculation shows that\vspace*{-1pt}
%
%e7.29 #&#
\begin{eqnarray}
\label{Gammaint} \Gamma_i^U(s)\cap\Gamma_j^V(s_i+s-t_j)=
\varnothing
\nonumber
\\[-10pt]
\\[-10pt]
\eqntext{\mbox{for }s_i< t_j\le s_i+s\mbox{ and }|y_j-x_i|>2\bigl(\sqrt\ep+s^{(1/2)-\delta_0}\bigr).}
\end{eqnarray}
If $F_i(t)=\bigcap_{j:t_j\le s_i+t}\{\rho_j^V>2t\}$, then use \eqref
{SKjVinc}\vspace*{1pt} and \eqref{Gammaint} in \eqref{Kimassbnd1} to see that on
$D_i(t)\cap F_i(t)$, for $s<t\wedge\rho_i$,\vspace*{-1pt}
%
%e7.30 #&#
\begin{eqnarray}
\label{KLi}
&&
K^{i,U}_{s_i+s}(1)\nonumber
\\[-1pt]
&&\qquad\le \sum
_j \1\bigl(s_i<t_j\le
s_i+s, |y_j-x_i|\le2\bigl(\sqrt
\ep+s^{(1/2)-\delta_0}\bigr)\bigr) K^{j,V}_{s_i+s}(1)
\\[-1pt]
&&\qquad\equiv L^i(s).\nonumber\vspace*{-1pt}
\end{eqnarray}
Note that $L^i$ is a nondecreasing process. If we sum the second
equation in \eqref{UVdefn} over $j$ satisfying $s_i<t_j\le s_i+s$,
$|y_j-x_i|\le2(\sqrt\ep+s^{(1/2)\delta_0})$, and denote this
summation by $\sum_j^{(i)}$, then\vspace*{-1pt}
%
%e7.31 #&#
\begin{eqnarray}\label{Libnd}
\nonumber
L^i(s)&\le&\sum_j^{(i)}K_{s_i+s}^{j,V}(1)+V_{s_i+s}^j(1)
\nonumber\\[-1pt]
&=&\int\int\1\bigl(s_i<t'\le
s_i+s, \bigl|y'-x_i\bigr|\le2\bigl(\sqrt\ep
+s^{(1/2)-\delta_0}\bigr)\bigr)\eta_\ep^-\bigl(dt',dy'
\bigr)
\\[-1pt]
&&{} +\sum_j^{(i)}
\int_0^{s_i+s}\int_{\R} V
\bigl(s',x\bigr)^{\gamma-(1/2)}V^j\bigl(s',x
\bigr)^{1/2}W^{j,V}\bigl(ds',dx\bigr).\nonumber
\end{eqnarray}
Now take means in \eqref{Libnd}, use \eqref{Qiunif} and use a standard
localization argument to handle the $Q_i$ martingale term, to conclude that\vspace*{-1pt}
\begin{eqnarray*}
&& E_{Q_i}\bigl(L^i(s)\bigr)
\\[-1pt]
&&\quad\le E_{Q_i} \biggl(\int\int\1\bigl(s_i<t'
\le s_i+s, \bigl|y'-x_i\bigr|\le 2\bigl(\sqrt\ep
+s^{(1/2)-\delta_0}\bigr)\bigr)\eta^-_\ep\bigl(dt',dy'
\bigr) \biggr)
\\[-1pt]
&&\quad=\sum_j \1(s_i<j\ep\le
s_i+s)\ep\\
&&\quad\qquad{}\times \int_0^1\int
_0^1\int_{y_j-\sqrt
\ep
}^{y_j+\sqrt\ep}J
\bigl(\bigl(y_j-y'\bigr)\ep^{-1/2}\bigr)
\ep^{-1/2}
\\
&&\hspace*{90pt}\qquad{}\times\1\bigl(\bigl|y'-x_i\bigr|\le\bigl(2\sqrt
\ep+2s^{(1/2)-\delta
_0}\bigr)\bigr)\,dy'\,dy_j\,dx_i
\\
&&\quad\le \sum_j \1(s_i<j\ep\le
s_i+s)\ep\int_0^1\int
_0^1 \1 \bigl(|y_j-x_i|
\le \bigl(3\sqrt\ep+2s^{(1/2)-\delta_0}\bigr)\bigr)\,dy_j\,dx_i
\\
&&\quad\le 2\bigl(3\sqrt\ep+2s^{(1/2)-\delta_0}\bigr) \biggl(\sum
_j \1(s_i<j\ep \le s_i+s)\ep \biggr)
\\
&&\quad\le 6\bigl(\sqrt\ep+s^{(1/2)-\delta_0}\bigr) (s+\ep)\le12(s+\ep
)^{(3/2)-\delta_0}.
\end{eqnarray*}
We take $s= 2^{-n}$ in the above, use Markov's inequality, and sum over
$n$ to conclude that
for some $c(\delta_0)>0$ independent of $\ep$,
\[
Q_i \biggl(L^i\bigl(2^{-n}\bigr)\le \biggl(
\frac{2^{-n-1}+\ep}{2} \biggr)^{(3/2)-2\delta
_0}\mbox{ for }N\le n\le
\log_2(1/\ep) \biggr)\ge1-c(\delta _0)2^{-N\delta_0}.
\]\eject\noindent
Recall that $L^i(\cdot)$ is nondecreasing and consider $s\in
[2^{-n-1},2^{-n}]$ to see that above implies that for $2^{-N}\ge\ep$,
\[
Q_i\bigl(L^i(s)\le(s+\ep)^{(3/2)-2\delta_0}\mbox{ for all
}s\in \bigl[0,2^{-N}\bigr]\bigr)\ge1-c(\delta_0)2^{-N\delta_0}.
\]
An easy interpolation argument in $N$ now shows that for some
$c_0(\delta_0)$, independent of $\ep$,
%
%e7.32 #&#
\begin{eqnarray}
\label{LiBC} Q_i\bigl(L^i(s)\le(s+\ep)^{(3/2)-2\delta_0}
\mbox{ for }0\le s\le u\bigr)\ge 1-c_0(\delta_0) (u\vee
\ep)^{\delta_0}
\nonumber
\\[-8pt]
\\[-8pt]
 \eqntext{\forall u\ge0.}
\end{eqnarray}

Apply \eqref{LiBC} in \eqref{KLi} and conclude
%
%e7.33 #&#
\begin{eqnarray}
\label{thetaboundI}
\nonumber
Q_i(\theta_i<\rho_i
\wedge t)&\le& Q_i\bigl(K^{i,U}_{s_i+s}(1)>(s+
\ep)^{(3/2)-2\delta_0}\ \exists s<\rho _i\wedge t\bigr)
\\
&\le& Q_i\bigl(F_i(t)^c
\bigr)+Q_i\bigl(D_i(t)^c\cap
F_i(t)\bigr)\nonumber\\
&&{} +Q_i\bigl(L^i(s)>(s+\ep)\ \exists s<
\rho_i\wedge t\bigr)
\\
&\le& Q_i \biggl(\bigcup_{j\le(2t/\ep)\wedge N_\ep}
\bigl\{\rho ^V_j\le 2t\bigr\} \biggr)+Q_i
\bigl(D_i(t)^c\cap\bigl\{\underline\rho^V_i>2t
\bigr\}\bigr)
 \nonumber\\
 &&{}+c_0(\delta_0) (t\vee\ep)^{\delta_0}.\nonumber
\end{eqnarray}
Recall from Section~\ref{intro} that $N_\ep=\lfloor\ep^{-1}\rfloor$.
The second term is at most $c_{\scriptsize{\ref{presimass}}}(\ep\vee t)^{\bar
\delta
}$ by Lemma~\ref{presimass}, and by Lemma~\ref{rhoVbnd}, if $4t\le1$
and $\ep\le1/2$, the first term is at most
\[
Q_i \biggl(\bigcup_{j\le4tN_\ep}\bigl\{
\rho_j^V\le2t\bigr\} \biggr)\le8c_{\scriptsize{\ref
{rhobnd}}}(t\vee
\ep)t\le8c_{\scriptsize{\ref{rhobnd}}}(t\vee\ep).
\]
If $4t>1$ or $\ep>1/2$, the above bound is trivial as $c_{\scriptsize{\ref
{rhobnd}}}\ge1$.
We conclude from~\eqref{thetaboundI} that
\[
Q_i(\theta_i<\rho_i\wedge t)
\le8c_{\scriptsize{\ref{rhobnd}}}(t\vee\ep )+c_{\scriptsize{\ref
{presimass}}}(\ep\vee t)^{\bar\delta}+c_0(
\delta_0) (t\vee\ep )^{\delta_0}.
\]
The result follows because $\delta_0\le\bar\delta\le1$.
\end{pf*}

\begin{appendix}\label{app}
%s8 #&#
\section{Moment bounds, tightness and proof of Proposition~\texorpdfstring{\lowercase{\protect\ref
{prop21}}}{2.2}}\label{secspdegr}
\setcounter{equation}{0}
We start with a moment bound obtained by a modification of the proof of
Lemma~4.2 in \citet{mp92}. Let $p(t,x)=p_t(x)$ denote that Gaussian
kernel, that is,
%
%e8.1 #&#
\begin{equation}
\label{Gauss} p_t(x)= \frac{1}{\sqrt{2\pi
t}}e^{-{x^2}/{(2t)}},\qquad t>0, x
\in\R.
\end{equation}
Let $S_t$ denote the corresponding semigroup, so $S_tf=p_t*f$ for
appropriate functions $f$.

%le8.1 #&#
\begin{lemma}\label{qmom} For any $q\ge1$ and $\lambda,T>0$ there is a
$C_{T,\lambda,q}$ such that for all $\ep\in(0,1]$:
\begin{longlist}[(a)]
\item[(a)] $\sup_{t\le T}\int e^{\lambda|x|}E(\bar U(t,x)^q+\bar
V(t,x)^q) \,dx\le C_{T,\lambda,q}$,

\item[(b)] $\sup_{t\le T,x\in\R}e^{\lambda|x|}E(\bar
U(t,x)^q+\bar
V(t,x)^q)\le C_{T,\lambda,q}$.
\end{longlist}
\end{lemma}
%
%re8.2 #&#
\begin{remark}
\label{rem09081} Lemma~\ref{qmom} and Theorem~1.1 of~\citet{myt98w}
easily imply uniqueness in law of each of
$\bU$ and $\bV$ separately for a pair
$(\bU,\bV)$ solving~\eqref{eq28}. To show the uniqueness in law for
the pair $(\bU,\bV)$, one should follow the proof of Theorem~1.1
of~\citet{myt98w} and derive the counterpart of Proposition~2.3
from~\citet{myt98w}, which is the main ingredient of the proof. More
specifically, suppose $t\in[s_i, t_i)$ for some $i\in\N_{\ep}$.
Following the argument from~\citet{myt98w}, for any
nonnegative $\phi_1,\phi_2\in L^1(\R)$, one can easily construct a
sequence of $M_F(\R)^2$-valued processes
$\{(Y^{1,n}, Y^{2,n})\}_{n\geq0}$ such that $\{Y^{1,n}\}_{n\geq1}$
and $\{Y^{2,n}\}_{n\geq1}$ are independent, and for any $(\bU,\bV)$
solving~\eqref{eq28} we have
%
%e8.2 #&#
\begin{eqnarray}
\label{eq2912}&& E \bigl[ e^{-\langle\phi_1, \bU_{t}\rangle+\langle\phi_2, \bV
_{t}\rangle} \bigr]
\nonumber
\\[-8pt]
\\[-8pt]
\nonumber
&&\qquad= \lim_{n\rightarrow\infty}
E \bigl[ e^{-\langle Y^{1,n}_{t-s_i},
\bU
_{s_i}\rangle+\langle Y^{2,n}_{t-s_i}, \bV_{s_i}\rangle
}|Y^{1,n}_0=\phi _1,
Y^{2,n}_0=\phi_2 \bigr].
\end{eqnarray}
A similar expression can be derived for $t\in[t_i, s_{i+1}), i\in\N
_{\ep}$, and then uniqueness in law for
the pair $(\bU,\bV)$ follows by standard argument: see again~\citet
{myt98w} where the single process without immigration is treated.
\end{remark}

%pa8.subsection.subsubsection.1 #&#
\begin{pf*}{Proof of Lemma~\protect\ref{qmom}} It suffices to
consider $\bar U$.
We let $C$ denote a constant which may depend on $q$, $\lambda$ and
$T$, and which may change from line to line. Note that equation~(\ref
{eq28}) for $\bU$ can be rewritten in the so-called mild form [see
Theorem~2.1 of \citet{shi94}]
%
%e8.3 #&#
\begin{eqnarray}
\label{eqmild1}
\bU_t(x) &=& \int_0^t
\int_{\R} p_{t-s}(x-y)\eta^+_{\ep
}(ds,dy)
\nonumber
\\[-8pt]
\\[-8pt]
\nonumber
&&{} +\int_0^t\int_{\R}
p_{t-s}(x-y) \bU(s,y)^{\gamma} \bW^{U}(ds,dy),\qquad t\geq0,
x\in\R.
\end{eqnarray}
Let $N(t,x)$ denote the stochastic integral term in the above.
The first term on the right-hand side of~(\ref{eqmild1}) can be
rewritten as
%
%e8.4 #&#
\begin{equation}
I_1(t,x)=I_{1,\ep}(t,x)=\sum_{s_i\in\CG_{\ep}^{\mathrm{odd}}, s_i\leq t}
\int_{\R} p_{t-s_i}(x-y) J^{x_i}_{\ep}(y)
\,dy
\end{equation}
(the meaning of the above if $t=s_i$ some $i$ is obvious). Recall that
$x_i\in[0,1]$ and so $y$ in\vadjust{\goodbreak} the above integral may be restricted to
$|y|\le2$. Therefore for $s_i\le t\le T$,
%
%e8.5 #&#
\begin{equation}
\label{expabs}e^{\lambda|x|}p_{t-s_i}(x-y)\le Cp_{2(t-s_i)}(x-y).
\end{equation}
It follows that
%
%e8.6 #&#
\begin{eqnarray}
\label{eqpmom2} &&\sup_{t\le T,x\in\R}e^{\lambda|x|}I_1(t,x)\nonumber
\\
\nonumber
&&\qquad\le \sum_{s_i\leq t- 2\ep} C (t-s_i)^{-1/2}
\ep\\
&&\qquad\quad{}+ \sum_{t-2\ep
<s_i<t}\sqrt {\ep} \int_{\R}
p_{2(t-s_i)}(x-y) \,dy
 +\1(s_i=t)e^{\lambda|x|}J_\ep^{x_i}(x)
\\
\nonumber
 &&\qquad\leq C \biggl[\int_{0}^{t}(t-s)^{-1/2}
\,ds + \ep^{1/2} \biggr]
\\
\nonumber
 &&\qquad\leq C,
\end{eqnarray}
uniformly on $\ep\in(0,1]$.
By \eqref{eqmild1} and \eqref{eqpmom2} we have for $t\le T$ and all $x$,
%
%e8.7 #&#
\begin{eqnarray}
\label{Uqmom1} E\bigl(\bar U(t,x)^q\bigr)&\le& C \bigl[E
\bigl(I_1(t,x)^q\bigr)+E\bigl(\bigl|N(t,x)\bigr|^q
\bigr) \bigr]
\nonumber
\\[-8pt]
\\[-8pt]
\nonumber
&\le& C \bigl[e^{-\lambda|x|}+E\bigl(\bigl|N(t,x)\bigr|^q\bigr)
\bigr].
\end{eqnarray}

For $q\ge1$ and $\lambda,t>0$ let
\[
\nu(q,\lambda,t)=\sup_{0\leq s\leq t}\int e^{\lambda|x|}E \bigl[\bar
U(s,x)^q \bigr]\,dx,
\]
and note that $\nu$ implicitly depends on $\ep$. Using the
Burkholder--Davis--Gundy
inequality and Jensen's inequality, we get for $q\ge2$,
%
%e8.8 #&#
\begin{eqnarray}
\label{Nqbound} && E \bigl[\bigl|N(t,x)\bigr|^q \bigr]\nonumber
\\
&&\qquad\leq CE \biggl[ \biggl(\int_{0}^{t}
\int p_{t-s}(x-y)^2\bar U(s,y)^{2\gamma}\,dy\,ds
\biggr)^{q/2} \biggr]
\\
\nonumber&&\qquad\leq CE \biggl[\int_{0}^{t}\int
p_{t-s}(x-y)^2\bar U(s,y)^{\gamma q}\,dy\,ds \biggr]
\\
\nonumber
&&\qquad\quad{}\times \biggl(\int_{0}^{t}\int
p_{t-s}(x-y)^2\,dy\,ds \biggr)^{(q/2)-1}
\\
\nonumber
&&\qquad\leq Ct^{(q-2)/4}E \biggl[\int_{0}^{t}
\int p_{t-s}(x-y)^2\bigl[\bar U(s,y)^{q/2}+\bar
U(s,y)^{q}\bigr]\,dy\,ds \biggr].
\end{eqnarray}
The final inequality follows because $p_{t-s}(x-y)^2\leq
(t-s)^{-1/2}p_{t-s}(x-y)$ and $a^{\gamma q}\le a^{q/2}+a^q$. A short
calculation using the above bound, just as in the bottom display on
page 349 of \citet{mp92} shows that
\begin{eqnarray}
\nonumber
\nu(q,\lambda,t)&\le& C \biggl[1+\sup_{s\le t}\int
e^{\lambda|x|}E\bigl(\bigl|N(t,x)\bigr|^q\bigr) \,dx \biggr]
\\
\eqntext{\mbox{[by \eqref{Uqmom1}\mbox{ with $2\lambda$ in place of $
\lambda$}]}}
\\
\nonumber
&\le& C+C\int_0^t
(t-s)^{-1/2}\bigl[\nu(q/2,\lambda,s)+\nu (q,\lambda,s)\bigr] \,ds
\\
\nonumber
&\le& C \biggl[1+\nu(q/2,\lambda,t)+\int_0^t(t-s)^{-1/2}
\nu (q,\lambda,s) \,ds \biggr].
\end{eqnarray}

A generalized Gronwall inequality [e.g., see Lemma~4.1 of \citet{mp92}]
shows that the above implies that for $q\ge2$,
%
%e8.9 #&#
\begin{equation}
\label{nuind} \nu(q,\lambda,t)\le\bigl(1+\nu(q/2,\lambda,t)\bigr)\exp
\bigl(4Ct^{1/2}\bigr)\qquad \mbox{for all }t\le T.
\end{equation}
The obvious induction on $q=2^n$ will now give (a) providing we can show
%
%e8.10 #&#
\begin{equation}
\label{init} \nu(1,\lambda,T)\le C.
\end{equation}
It follows from \eqref{eqmild1} and an argument using localization and
Fubini's theorem that
\begin{eqnarray*}
\sup_{t\le T}\sup_x e^{\lambda|x|}E
\bigl[\bar U(t,x)\bigr]
\le\sup_{t\le
T}\sup_x
e^{\lambda|x|}E\bigl[I_1(t,x)\bigr]\le C,
\end{eqnarray*}
the last inequality by \eqref{eqpmom2}. By optimizing over $\lambda$
we get
\eqref{init}. Therefore we have proved Lemma~\ref{qmom} part (a) except
for one detail. To use
Lemma~4.1 in \citet{mp92} to derive \eqref{nuind} we need to know that
$\nu(q,\lambda,T)<\infty$ (the bound can now depend on $\ep$). To
handle this
issue one can localize just as in \citet{mp92} using the facts that
$t\to\bar U_t$ is in
$D(\R_+,C^+_{\mathrm{rap}})$, and (from Proposition~\ref{thm11} and
$\bar U=\sum_i\bar U^i$) that the jumps of $\bar U$ occur
at $\{s_i\}$ with the $i$th jump equaling $J^{x_i}\le\sqrt\ep$.

Turning to Lemma~\ref{qmom} part (b), it suffices to consider $q>2$. By
\eqref{eqmild1}, \eqref{eqpmom2} and the first line of \eqref
{Nqbound} for $t\le T$, $p=q/(q-2)$ and $p'=q/2$, we have by H\"older's
inequality
\begin{eqnarray*}
&&\hspace*{-4pt}\sup_x e^{\lambda|x|}E \bigl[\bar
U(t,x)^q \bigr]
\\[-1pt]
&&\hspace*{-10pt}\qquad\le C \biggl(1+\sup_x E \biggl[ \biggl(\int
_0^t\int \bigl[p_{t-s}(x-y)^{1/p}e^{2\lambda|x|/q-2\lambda|y|/q}
\bigr]
\\[-1pt]
&&\hspace*{93pt}\hspace*{-10pt}\qquad\quad{}\times
\bigl[e^{2\lambda
|y|/q}\bar U(s,y)^{2\gamma}\bigr]p_{t-s}(x-y)^{2-(1/p)}\,dy\,ds
\biggr)^{q/2} \biggr] \biggr)
\\[-1pt]
&&\hspace*{-10pt}\qquad\le C \biggl(1+\sup_x \biggl(\int_0^t
\hspace*{-1pt}\int p_{t-s}(x-y)e^{2\lambda
p|x|/q-2\lambda p|y|/q}\,dy(t-s)^{-1+(1/2p)}\,ds
\biggr)^{q/2p}
\\[-1pt]
&&\hspace*{-10pt}\qquad\quad{}\times E \biggl[\int
_0^t\int e^{2\lambda
p'|y|/q}\bar
U(s,y)^{2\gamma p'}\,dy(t-s)^{-1+(1/2p)}\,ds \biggr] \biggr)
\\[-1pt]
&&\hspace*{-10pt}\qquad\le C \biggl(1+\int_0^t(t-s)^{-(q+2)/(2q)}\,ds
\nu(\gamma q,\lambda,t) \biggr)\\[-1pt]
&&\hspace*{-10pt}\qquad\le C.
\end{eqnarray*}
In the next to last line we have used Lemma~6.2 of \citet{shi94} and in the
last line we have used Lemma~\ref{qmom} part (a).
\end{pf*}

%pa8.subsection.subsubsection.2 #&#
\begin{pf*}{Proof of Lemma~\protect\ref{pmom}} It suffices to
consider $\bar U$.
Let $C$ denote a constant depending on $q$ and $T$ which may change
from line to line. We adapt the proof of Lemma A.3 of~\citet{mps06} to
the white noise setting and with
$\lambda=0$.

By \eqref{eqmild1}, \eqref{eqpmom2} and the continuity properties of
$\bar U$, we have\vspace*{-1pt}
\begin{eqnarray*}
&&E \Bigl[\sup_{t\le T,x\in\R}\bar U(t,x)^{q}
\Bigr]
\\[-1pt]
&&\qquad\le C_{q,T} \biggl(1 + E \biggl[\sup_{t\le T,x\in\Q, t\in\Q_+}\biggl
\llvert \int_0^t\int_{\R}
p_{t-s}(x-y) \bU(s,y)^{\gamma} \bW^{U}(ds,dy)\biggr
\rrvert ^q \biggr] \biggr).
\end{eqnarray*}
To handle the above expectation we carry out the argument in the proof
of Lemma~A.3 of \citet{mps06} with $\lambda=0$ and $W$ a white noise. We
take $a\in(0,1/4)$ and $q>\frac{3}{2a}$ in that work. With this choice
of $q$, the arguments in Lemma~A.3 of \citet{mps06} then go through to
show that
the expectation in the above is at most
\begin{eqnarray*}
&& C\int_0^T\int E \biggl[ \biggl|\int
_0^t\int(t-s)^{-a}p_{t-s}(x-y)
\bar U(s,y)^{\gamma} \,d\bar W^U(s,y)\biggr |^q
\biggr]\,dx\,dt
\\[-1pt]
&&\qquad\le C\int_0^T\int E \biggl[ \biggl|\int
_0^t\int (t-s)^{-2a}p_{t-s}(x-y)^2
\bar U(s,y)^{2\gamma}\,dy\,ds \biggr|\biggr]^{q/2}\,dx\,dt
\\[-1pt]
&&\qquad\le C\int_0^T\int \biggl[\int
_0^t\int (t-s)^{-2a-(1/2)}p_{t-s}(x-y)E
\bigl(\bar U(s,y)^{q\gamma}\bigr)\,dy\,ds \biggr]\,dx\,dt
\\[-1pt]
&&\qquad\le C,
\end{eqnarray*}
by Fubini, Lemma~\ref{qmom} part (a) and the choice of $a$. This gives
the result for $q>3/2a$ and hence for all $q>0$.
\end{pf*}%\qed

We turn next to the proof of Proposition~\ref{prop21} which is fairly
standard. We
follow the proof in Section~4 of \citet{mp92}, where a
similar existence proof is given. The main difference is the
immigration term
in the present situation.

By the mild form of \eqref{eqspde-u-ep} we have
%
%e8.11 #&#
\begin{eqnarray}
\label{eqspde-u-ep1} u_\ep(t,x) &=& \sum_{i}
\int p(t-s_i,y-x)J^{x_i}(y)\1(t\geq s_i)\,dy
\nonumber\\
&&{}- \sum_{j}\int p(t-t_j,y-x)J^{y_j}(y)
\1(t\geq t_j)\,dy
\nonumber
\\[-8pt]
\\[-8pt]
\nonumber
&&{}+ \int_{0}^{t}\int p(t-s,y-x)\bigl|u_\ep(s,y)\bigr|^\gamma
W(ds,dy)
\nonumber
\\
&\equiv&I_{1,\ep}(t,x)-I_{2,\ep}(t,x)+N_\ep(t,x).
\nonumber
\end{eqnarray}

Now we give a modified version of Lemma~4.4 of \citet{mp92}. The only
difference is that Lemma~4.4 of \citet{mp92} \,deals with $C^+_{\mathrm{rap}}$
instead of
$C_{\mathrm{rap}}$, but the proof carries over with almost no change.

%le8.3 #&#
\begin{lemma}
\label{mp44}
Let $\{X_n(t,\cdot)\dvtx t\geq0, n\in\N\}$ be a sequence of continuous
$C_{\mathrm{rap}}$-valued
processes. Suppose $\exists q>0,\gamma>2$ and $\forall T,\lambda>0$
$\exists C=C(T,\lambda)>0$ such that
%
%e8.12 #&#
\begin{eqnarray}
\label{lem-X-replace-by-N} E \bigl[\bigl|X_n(t,x)-X_n
\bigl(t',x'\bigr)\bigr|^q \bigr] \leq C
\bigl(\bigl|x-x'\bigr|^\gamma+\bigl|t-t'\bigr|^\gamma
\bigr)e^{-\lambda|x|}
\nonumber
\\[-8pt]
\\[-8pt]
\eqntext{\forall t,t'\in[0,T], \bigl|x-x'\bigr|\leq1, n\in\N.}
\end{eqnarray}
If $\{P_{X_n(0)\dvtx n\in\N}\}$ is tight on $C_{\mathrm{rap}}$, then $\{
P_{X_n}\dvtx n\in\N\}$
is tight on $C(\R_+,C_{\mathrm{rap}})$.
\end{lemma}

We also need Lemma~4.3 of \citet{mp92}:
%
%le8.4 #&#
\begin{lemma}
\label{lemmp-43}
If $T,\lambda>0$ there is a constant $C(T,\lambda)<\infty$ such that
\begin{eqnarray}
&& \int_{0}^{t}\int
\bigl(p_{t-s}(y-x)-p_{t'-s}\bigl(y-x'\bigr)
\bigr)^2e^{-\lambda
|y|}\,dy\,ds
\nonumber\\
&&\qquad\leq C(T,\lambda) \bigl(\bigl|x-x'\bigr|+\bigl(t-t'
\bigr)^{1/2} \bigr)e^{-\lambda|x|}
\nonumber\\
 \eqntext{\forall 0<t'<t\leq T,\bigl |x-x'\bigr|\leq1, \lambda>0,}
\end{eqnarray}
where $p_u(z)$ is defined to be 0 if $u<0$.
\end{lemma}

Clearly $t\to I_{\ell,\ep}(t,\cdot)$ is in $D(\R_+,C_{\mathrm{rap}})$ with
jumps only at $\{s_i\}$ for $\ell=1$ and at $\{t_j\}$ if $\ell=2$.
It is fairly easy to see that for $t,x$ fixed $I_{\ell,\ep}(t,x)$
converges in probability
to
\[
I(t,x)=\int_{0}^{t\wedge1}\int_{0}^{1}p(t-s,x-y)\,dy\,ds
\]
by the weak law of large numbers.
We need convergence in path space. It is easy to check that $t\to
I(t,\cdot)$ is in $C(\R_+,C_{\mathrm{rap}})$.

%le8.5 #&#
\begin{lemma}\label{Iconv} For $\ell=1,2$, $I_{\ell,\ep}$ converges in
probability in $D(\R_+,C_{\mathrm{rap}})$ to $I$ as $\ep\downarrow0$.
\end{lemma}
%
%pa8.subsection.subsubsection.3 #&#
\begin{pf}
The argument is routine if a bit tedious. We sketch
the proof for $\ell=2$ where $t_j=j\ep$. If $\delta=\ep^{3/4}$, write
\begin{eqnarray*}
I_{2,\ep}(t,x)&=&\sum_{t_j\le t-\delta}\ep\int
\bigl[p_{t-t_j}(y_j-x+\sqrt\ep w)-p_{t-t_j}(y_j-x)
\bigr]J(w) \,dw
\\
&&{}+\sum_{t-\delta<t_j\le t}S_{t-t_j}J^{y_j}_\ep(x)+
\sum_{t_j\le
t-\delta
}\ep p_{t-t_j}(y_j-x)
\\
&=&T_{1,\ep}+T_{2,\ep}+T_{3,\ep}.
\end{eqnarray*}
It is easy to check that for any $\lambda,T>0$,
\[
\sup_{t\le T,x\in\R}e^{\lambda|x|}\bigl|T_{2,\ep}(t,x)\bigr|\le
C_{T,\lambda
}\delta/\sqrt\ep\to0
\]
and
\[
\sup_{t\le T,x\in\R}e^{\lambda|x|}\bigl|T_{1,\ep}(t,x)\bigr|\le
C_{T,\lambda
}\sqrt \ep\bigl(1+\ln(1/\ep)\bigr)\to0.
\]
So it suffices to show that $T_{3,\ep}$ converges in probability in
$D(\R_+,C^+_{\mathrm{rap}})$ to $I$.

We next write
\begin{eqnarray*}
T_{3,\ep}(t,x)&=&\sum_{t_j\le t-\delta} \biggl(\ep
p_{t-t_j}(y_j-x)-\ep \int_0^1p_{t-t_j}(y-x)
\,dy \biggr)
\\
&&{} +\sum_{t_j\le t-\delta}\ep\int_0^1p_{t-t_j}(y-x)
\,dy
\\
&\equiv& T_{4,\ep}+T_{5,\ep}.
\end{eqnarray*}
$T_{5,\ep}$ is a Riemman sum for $\int_0^{t\wedge1}\int_0^1
p_{t-s}(y-x) \,dy \,ds$ (note that $t_j\le1$, whence the truncation by
$1$), and using the $t-\delta$ cut-off, the Gaussian tail and $y\in
[0,1]$, it is easy to see that for any $\lambda,T>0$,
\[
\lim_{\ep\to0}\sup_{t\le T, x\in\R}e^{\lambda|x|}
\biggl|T_{5,\ep
}-\int_0^{t\wedge1}\int
_0^1 p_{t-s}(y-x) \,dy \,ds \biggr|=0.
\]

Therefore it remains to show that $T_{4,\ep}\to0$ in probability in
$D(\R_+,C_{\mathrm{rap}})$. $T_{4,\ep}$ is a sum of mean $0$ independent
random variables, and so one easily sees that
\[
E\bigl(T_{4,\ep}(t,x)^2\bigr)\le\ep^2\sum
_{t_j\le t-\delta
}p_{2(t-t_j)}(0)\to 0\qquad\mbox{as }\ep
\downarrow0.
\]
If we could show for any $\ep_n\downarrow0$,
%
%e8.13 #&#
\[
\{T_{4,\ep_n}\dvtx n\}\mbox{ is $C$-tight in $D(
\R_+,C_{\mathrm{rap}})$}
\]
the result would follow as the only possible weak limit point is $0$ by
the above.

Let $\hat p_{t-t_j}(y_j-x)=p_{t-t_j}(y_j-x)-\int_0^1 p_{t-t_j}(y-x)\,dy$ and
\[
[t-\delta]_\ep=\max\{j\ep\dvtx j\ep\le t-\delta, j\in\Z_+\}.
\]
To work in the space of continuous $C_{\mathrm{rap}}$-valued paths, we
interpolate $T_{4,\ep}$ linearly and define
\begin{eqnarray*}
\widetilde T_{4,\ep_n}(t,x)&=&\sum_{t_j\le[t-\delta
_n]_{\ep_n}} \ep\hat
p_{t-t_j}(y_j-x)
\\
&& {} + \bigl((t-\delta_n)-[t-\delta_n]_{\ep_n}
\bigr)\hat p_{t-[t-\delta
_n]_{\ep_n}-\ep_n}(y_{1+([t-\delta_n]_{\ep_n}/\ep_n)}-x),
\end{eqnarray*}
so that $t\to\widetilde T_{4,\ep_n}(t,\cdot)\in C(\R_+,C_{\mathrm{rap}})$. If
$d$ is the metric on $C_{\mathrm{rap}}$, then it is clear that
\[
\lim_{n\to\infty} \sup_{t\le T}\,d\bigl(\widetilde
T_{4,\ep_n}(t),T_{4,\ep
_n}(t)\bigr)=0\qquad\mbox{for all }T>0.
\]
Therefore it remains to show that
%
%e8.14 #&#
\begin{eqnarray}
\label{tightcrit} \{\widetilde {T}_{4,\ep_n}\dvtx n\}\mbox{ is tight in }
C(\R_+,C_{\mathrm{rap}}).
\end{eqnarray}
This is proved by a straightforward application of Lemma~\ref{mp44},
as we illustrate below.\vadjust{\goodbreak}

To illustrate the method of the aforementioned proof let us bound the
spatial moments and work with $T_{4,\vep}$, hence dropping the trivial
continuity correction and dependence on $n$. Assume $0\le t\le T$,
$\lambda>0$ and $|x-x'|\le1$. For $q\ge2$ we use a predictable square
function inequality of Burkholder [see Theorem~21.1 of \citet{bur73}]
as follows:
%
%e8.15 #&#
\begin{eqnarray}
\label{sqfnin}\nonumber
&&
e^{\lambda|x|}E \bigl[\bigl|T_{4,\ep}(t,x)-T_{4,\ep
}
\bigl(t,x'\bigr)\bigr|^q \bigr]
\\
&&\qquad\le e^{\lambda|x|}c_q \biggl[\biggl |\sum
_{t_j\le
[t-\delta]_\ep
}\ep^2E\bigl(\bigl(\hat p_{t-t_j}(y_j-x)-
\hat p_{t-t_j}\bigl(y_j-x'\bigr)
\bigr)^2\bigr) \biggr|^{q/2}
\\
&&\hspace*{34pt}\qquad\quad{}+\sum
_{t_j\le[t-\delta
]_\ep}\ep ^qE\bigl(\bigl|\hat p_{t-t_j}(y_j-x)-
\hat p_{t-t_j}\bigl(y_j-x'\bigr)
\bigr|\bigr)^q \biggr].\nonumber
\end{eqnarray}

Now for $q\ge2$ and for, say $x>x'$,
\begin{eqnarray*}
&&e^{\lambda|x|}E \bigl[\bigl|\hat p_{t-t_j}(y_j-x)-
\hat p_{t-t_j}\bigl(y_j-x'\bigr)\bigr|^q
\bigr]
\\
&&\qquad\le ce^{\lambda|x|}\int_0^1\bigl|p_{t-t_j}(y-x)-p_{t-t_j}
\bigl(y-x'\bigr)\bigr|^q \,dy
\\
&&\qquad\le C_{\lambda,T}(t-t_j)^{-1/2}\int
_0^1\bigl|p_{t-t_j}(y-x)-p_{t-t_j}
\bigl(y-x'\bigr)\bigr|^{q-1} \,dy.
\end{eqnarray*}
In the last line we used the bound on $|x-x'|$ and the fact that $y\in
[0,1]$ to use the Gaussian tail of $(p_{t-t_j}(y-x)+p_{t-t_j}(y-x'))$
to absorb the $e^{\lambda|x|}$ as in \eqref{expabs}. By using the
spatial derivative of $p_t(z)$ and then carrying out a change of
variables, we may bound the above by
\begin{eqnarray*}
&& C_{\lambda,T}(t-t_j)^{-1/2}\int
_0^1(t-t_j)^{-(q-1)/2}
\\
&&\hspace*{77pt}\quad{}\times \biggl|\int\1 \biggl(\frac{y-x}{\sqrt{t-t_j}}\le z
\le\frac{y-x'} {
\sqrt{t-t_j}}
\biggr)zp_1(z) \,dz \biggr|^{q-1} \,dy
\\
&&\qquad\le C_{\lambda,T}(t-t_j)^{-q+0.5}\bigl|x-x'\bigr|^{q-1}.
\end{eqnarray*}
We use the above in \eqref{sqfnin} with $q=2$ and general $q$ to
conclude that
\begin{eqnarray*}
&&e^{|x|}E\bigl[\bigl|T_{4,\ep}(t,x)-T_{4,\ep}
\bigl(t,x'\bigr)\bigr|^q\bigr]
\\
&&\qquad\le C_{\lambda,T} \biggl(\sum_{t_j\le[t-\delta]_\ep}\ep
^2(t-t_j)^{-3/2} \biggr)^{q/2}\bigl|x-x'\bigr|^{q/2}
\\
&&\qquad\quad{}+C_{\lambda,T}\sum
_{t_j\le
[t-\delta
]_{\ep}}\ep^q(t-t_j)^{-q+0.5}
\bigl|x-x'\bigr|^{q-1}
\\
&&\qquad\le C_{\lambda,T}\bigl|x-x'\bigr|^{q/2},
\end{eqnarray*}
where we used $\delta=\ep^{3/4}$, $q\ge2$ and an elementary
calculation in the last line. So taking $q>4$ gives the required
spatial increment bound in Lemma~\ref{mp44}.

A similar, but slightly more involved, argument verifies the
hypotheses of Lemma~\ref{mp44} for the time increments. Here when
$0\le t'-t\le\ep$ the linear interpolation term must be used and the
cases $[t'-\delta]_\ep=[t-\delta]_\ep$ and $[t'-\delta]_\ep
=[t-\delta
]_\ep+\ep$ are treated separately. The details are left for the reader.
This establishes \eqref{tightcrit} and so completes the proof.
\end{pf}

Next we apply Lemma~\ref{mp44} to $X_n(t,x)=N_{\ep_n}(t,x)$ for any
$\ep_n\downarrow0$ by showing that
(\ref{lem-X-replace-by-N}) holds for $X_n=N_{\ep_n}$.

%le8.6 #&#
\begin{lemma}
\label{lemN-satisfies}
$\exists q>0,\gamma>2$ and $\forall T,\lambda>0$
$\exists C=C(T,\lambda)>0$ such that
%
%e8.16 #&#
\begin{eqnarray}
\label{mp47}
E \bigl[\bigl|N_\ep(t,x)-N_\ep
\bigl(t',x'\bigr)\bigr|^q \bigr] &\leq& C
\bigl(\bigl|x-x'\bigr|^\gamma+\bigl|t-t'\bigr|^\gamma
\bigr)e^{-\lambda|x|}
\nonumber
\\[-8pt]
\\[-8pt]
\eqntext{\forall t,t'\in[0,T], \bigl|x-x'\bigr|\leq1, 0<\ep<1.}
\end{eqnarray}
\end{lemma}

%pa8.subsection.subsubsection.4 #&#
\begin{pf}
Here we follow the proof of Proposition~4.5 of \citet{mp92}.
%For convenience we will omit the dependence on $\ep$ and simply write
%$N(t,x)$, while noting
%that it will be crucial that our constants \,do not depend on the
%invisible $\ep$.
Let $q\ge1$, $\lambda>0$, $0\leq t'<t\leq T$ and $|x-x'|\leq1$. First,
Jensen's inequality shows that for nonnegative functions $f,g$, we have
\[
\biggl(\int fg \biggr)^q \leq \biggl(\int f^qg \biggr)
\biggl(\int g \biggr)^{q-1}.
\]
Now using the Burkholder--Davis--Gundy inequality and Jensen's
inequality and allowing $c_q$ to vary
from line to line, we find
\begin{eqnarray}
&& E \bigl[\bigl|N_\ep(t,x)-N_\ep\bigl(t',x'
\bigr)\bigr|^{2q} \bigr]
\nonumber\\
&&\qquad\leq c_qE \biggl[ \biggl(\int_{0}^{t}
\int \bigl(p_{t-s}(y-x)-p_{t'-s}\bigl(y-x'\bigr)
\bigr)^2e^{-\lambda|y|}\nonumber\\
 &&\hspace*{142pt}{}\times e^{\lambda|y|}\bigl|u_\ep(s,y)\bigr|^{2\gamma}\,dy \,ds
\biggr)^q \biggr]
\nonumber\\
&&\qquad\leq c_qE \biggl[\int_{0}^{t}
\int\bigl|u_\ep(s,y)\bigr|^{2\gamma
q}e^{\lambda|y|(q-1)}
\bigl(p_{t-s}(y-x)-p_{t'-s}\bigl(y-x'\bigr)
\bigr)^2\,dy\,ds \biggr]
\nonumber\\
&&\qquad\quad{} \times \biggl( \int_{0}^{t}\int
\bigl(p_{t-s}(y-x)-p_{t'-s}\bigl(y-x'\bigr)
\bigr)^2 e^{-\lambda|y|}\,dy\,ds \biggr)^{q-1}
\nonumber\\
&&\qquad\leq c_qE \biggl[\int_0^t
\int\bigl|u_\ep(s,y)\bigr|^{8\gamma q}e^{4\lambda
|y|(q-1)}\,dy \,ds
\biggr]^{1/4}
\nonumber\\
&&\qquad\quad{} \times \biggl(\int_{0}^{t}\int
\bigl|p_{t-s}(y-x)-p_{t'-s}\bigl(y-x'
\bigr)\bigr|^{8/3}\,dy\,ds \biggr)^{3/4}
\nonumber\\
&&\qquad\quad{} \times C'(T,\lambda,q) \bigl(\bigl|x-x'\bigr|^{q-1}+
\bigl|t-t'\bigr|^{(q-1)/2}
\bigr)e^{-\lambda(q-1)|x|}\nonumber
\\
 \eqntext{\mbox{(H\"older's inequality and Lemma~\ref
{lemmp-43}}
)}
\\
&&\qquad\leq C'(T,\lambda,q) \bigl(\bigl|x-x'\bigr|^{q-1}+\bigl|t-t'\bigr|^{(q-1)/2}
\bigr)e^{-\lambda(q-1)|x|}\nonumber
\end{eqnarray}
by Lemma~\ref{qmom}(a) (recall that $|u_\ep|=|U_\ep-V_\ep|\le\bar
U_\ep
+\bar V_\ep$) and an elementary calculation.
The result follows.%\qed
\end{pf}

\begin{pf*}{Proof of Proposition~\ref{prop21}} Recall that $\ep
_n=\frac{1}{n}$. Lemma~\ref{lemN-satisfies} allows us to conclude that
$N_{\ep_n}(t,x)$ is tight in $C(\R_+,C_{\mathrm{rap}})$
as $n\to\infty$. Hence by Lemma~\ref{Iconv} and \eqref{eqspde-u-ep1},
$\{u_{\ep_n}\}$ is $C$-tight
in $D(\R_+,C_{\mathrm{rap}})$.

It remains to show that any limit point satisfies equation \eqref{spde}
(it will then necessarily be a $C_{\mathrm{rap}}$-valued solution). Recall
from \eqref{eqspde-u-ep} we have
%
%e8.17 #&#
\begin{eqnarray}
\label{eqspde-u-ep2}
\bigl\langle u_{\ep}(t),\phi\bigr\rangle&=&
\sum_{i}\1(s_i\le t)\bigl\langle
J_\ep^{x_i},\phi\bigr\rangle- \sum
_{j}\1(t_j\le t)\bigl\langle J_\ep^{y_j},
\phi\bigr\rangle
\nonumber
\\[-8pt]
\\[-8pt]
\nonumber
&&{} +\int_0^t \frac{1}{2}\bigl\langle
u_{\ep}(s),\Delta\phi\bigr\rangle \,ds + \int_0^t
\int\bigl|u_{\ep}(s,x)\bigr|^{\gamma} \phi(x) W(ds,dx)
\end{eqnarray}
for $\phi\in C^\infty_c$.

If $\phi\in C_c(\R)$, then a simple calculation using the strong law
of large
numbers shows that with probability 1,
%
%e8.18 #&#
\begin{eqnarray}
\label{immcvgt} \lim_{n\to\infty}\sum_{i}
\1(s_i\le t)\bigl\langle J_{\ep_n}^{x_i},\phi \bigr
\rangle &=&(t\wedge1)\int_{0}^{1}\phi(x)\,dx,
\nonumber
\\[-8pt]
\\[-8pt]
\nonumber
\lim_{n\to\infty}\sum_{j}
\1(t_j\le t)\bigl\langle J_{\ep
_n}^{y_j},\phi\bigr
\rangle &=&(t\wedge1)\int_{0}^{1}\phi(x)\,dx.
\end{eqnarray}
It is easy to interpolate in $t$ and conclude that the above
convergence is
uniform in $t$ with probability 1. By considering a countable dense set
of $\phi$ in
$C_c(\R)$, we may conclude that with probability 1 for all $\phi\in
C_c(\R)$ the
convergence in \eqref{immcvgt} holds uniformly in $t$.

Choose a subsequence
$\{n_k\}$ so that $u_{\ep_{n_k}}$ converges weakly to $u$ in
$D(\R_+,\break C_{\mathrm{rap}})$ where $u$ has continuous paths. To ease eye
strain, we
write $u_k$ for $u_{\ep_{n_k}}$. By Skorokhod's theorem we may change spaces
so that (recall convergence in cadlag space $D$ to a continuous path means
uniform convergence on compacts)
\[
\lim_{k\to\infty}\sup_{t\le T}\,d\bigl(u_{k}(t),u(t)
\bigr)=0 \qquad\mbox{for all }T>0 \mbox{ a.s.}
\]
This fact and the above convergence in \eqref{immcvgt} show that with
probability 1 for all $\phi\in C_c^\infty$, the left-hand side of
\eqref{eqspde-u-ep2} and first three terms on the right-hand side of
the same
equation converge uniformly in $t$ to the same terms but with $u$ in
place of~$u_\ep$, or in the case of \eqref{immcvgt}, to the right-hand side of\vadjust{\goodbreak}
\eqref{immcvgt}.
Hence the last term on the right-hand side of \eqref{eqspde-u-ep2}
must also converge uniformly in $t$ a.s. to a continuous limit
$M_t(\phi
)$. So for all $\phi\in C_c^\infty$ we have
%
%e8.19 #&#
\begin{eqnarray}
\label{martp} \langle u_t,\phi\rangle=\int_0^t
\frac{1}{2}\bigl\langle u(s),\Delta \phi \bigr\rangle \,ds+M_t(
\phi).
\end{eqnarray}
We see that $M_t(\phi)$ is the a.s. limit of the stochastic integral in
\eqref{eqspde-u-ep2}. Using the boundedness of the moments uniformly
in $\ep$
from Lemma~\ref{qmom}, it is now standard to deduce that $M_t(\phi)$
is a
continuous $\CF_t$-martingale with square function
$\int_0^t\int|u(s,x)|^{2\gamma}\phi(x)^2  \,dx \,ds$. Here $\CF_t$
is the
right continuous filtration generated by $t\to u_t$. It is also routine to
construct a white noise $W$, perhaps an enlarged space, so that
$M_t(\phi)=\int_0^t\int u(s,x)^{\gamma}\phi(x)\,dW(s,x)$ for all
$t\ge0$ a.s.
for all $\phi\in C_c^\infty$. Put this into \eqref{martp} to see that
$u$ is
a $C_{\mathrm{rap}}$-valued solution of \eqref{spde} and we are done. %
\end{pf*}

%s9 #&#
\section{Construction of approximate solutions and proof of
Proposition~\texorpdfstring{\lowercase{\protect\ref{thm11}}}{2.1}}\label{secconstr}
\setcounter{equation}{0}

Let us fix $\ep\in(0,1]$. For this $\ep$ we construct the sequence of
processes mentioned in Proposition~\ref{thm11}, approximating them by
a system of processes with ``soft-killing.''
Fix $n>0$, and define the sequence of processes $(U^{i,n}, V^{i,n}, \tU^{i,n},
\tV^{i,n})$ as follows. For any $\phi\in C^2_{b}(\R)$, let
%
%e9.1 #&#
{\fontsize{10.2}{12.2}\selectfont
\begin{eqnarray}\qquad
\cases{ %
\displaystyle U^{i,n}_t(
\phi) = \bigl\langle J^{x_i},\phi\bigr\rangle\1(t\geq s_i)
\vspace*{2pt}\cr
 \hspace*{45pt}\displaystyle{}+ \int_0^t \int_{\R}U^n(s,x)^{\gamma-1/2}
U^{i,n}(s,x)^{1/2} \phi(x) W^{i,n,U}(ds,dx)
\vspace*{2pt}\cr
\hspace*{45pt}\displaystyle{}+ \int_0^t
U^{i,n}_s\biggl(\frac{1}{2}\Delta\phi\biggr) \,ds - n
\int_0^t\bigl\langle U^{i,n}_s
V^{n}_s,\phi\bigr\rangle \,ds,\qquad  t\geq0, i\in
\NN_{\ep},
\vspace*{2pt}\cr
\displaystyle V^{j,n}_t(\phi) = \bigl\langle J^{y_j},\phi
\bigr\rangle\1(t\geq t_j)
\vspace*{2pt}\cr
\hspace*{47pt}\displaystyle{} + \int_0^t\int_{\R}
V^n(s,x)^{\gamma-1/2} V^{j,n}(s,x)^{1/2}
\phi(x) W^{j,n,V}(ds,dx)
\vspace*{2pt}\cr
\hspace*{47pt}\displaystyle{}+ \int_0^t
V^{j,n}_s\biggl(\frac{1}{2}\Delta\phi\biggr) \,ds - n
\int_0^t \bigl\langle V^{j,n}_s
U^{n}_s,\phi\bigr\rangle \,ds,\qquad t\geq0, j\in
\NN_{\ep},
\vspace*{2pt}\cr
\displaystyle\tU^{i,n}_t(\phi) = \int_0^t
\int_{\R} \bigl[ \bigl(\tU ^n(s,x)+U^n(s,x)
\bigr)^{2\gamma} - U^n(s,x)^{2\gamma}
\bigr]^{1/2}
\vspace*{2pt}\cr
\hspace*{72pt}\displaystyle{}\times \sqrt{\frac{\tU^{i,n}(s,x)}{\tU^n(s,x)}} \phi(x) \tW
^{i,n,U}(ds,dx)
\vspace*{2pt}\cr
\hspace*{42pt}\displaystyle{}+ \int_0^t
\tU^{i,n}_s\biggl(\frac{1}{2}\Delta\phi\biggr) \,ds +
n \int_0^t\bigl\langle U^{i,n}_s
V^{n}_s,\phi\bigr\rangle \,ds,\qquad t\geq0, i\in
\NN_{\ep},
\vspace*{2pt}\cr
\displaystyle\tV^{j,n}_t(\phi) = \int_0^t
\int_{\R} \bigl[ \bigl(\tV ^n(s,x)+V^n(s,x)
\bigr)^{2\gamma} - V^n(s,x)^{2\gamma}
\bigr]^{1/2}
\vspace*{2pt}\cr
\hspace*{72pt}\displaystyle{}\times \sqrt{\frac{\tV^{j,n}(s,x)}{\tV^n(s,x)}} \phi(x) \tW
^{j,n,V}(ds,dx)
\vspace*{2pt}\cr
\hspace*{46pt}\displaystyle{}+ \int_0^t
\tV^{j,n}_s\biggl(\frac{1}{2}\Delta\phi\biggr) \,ds +
n \int_0^t\bigl\langle V^{j,n}_s
U^{n}_s,\phi\bigr\rangle \,ds,\qquad t\geq0, j\in
\NN_{\ep}, }\label{tUVndefn}
\end{eqnarray}}
\hspace*{-2pt}where
\begin{eqnarray*}
U^{n}_t &=& \sum_i
U^{i,n}_t,\qquad  V^{n}_t = \sum
_j V^{j,n}_t,
\\
\tU^{n}_t& =& \sum_i
\tU^{i,n}_t, \qquad \tV^{n}_t =\sum
_j \tV^{j,n}_t,
\end{eqnarray*}
and $\{ W^{i,n,U}, W^{j,n,V},\tW^{k,n,U},\tW^{l,n,V}\}_{i,j,k,l\in
\NN
_{\ep}}$ is a collection of mutually independent white noises.
For $\phi\in C^2_{b}(\R)$, let $\{M_t^{i,n,U}(\phi)\}_{t\geq0}, \{
M_t^{j,n,V}(\phi)\}_{t\geq0}$,
$\{\tM_t^{i,n,U}(\phi)\}_{t\geq0},
\{\tM_t^{j,n,V}(\phi)\}_{t\geq0}$ denote the stochastic integrals on
the right-hand side of the
equations for $U^{i,n}, V^{j,n}, \tU^{i,n}, \tV^{j,n}$, respectively,
in~(\ref{tUVndefn}). For each $n$, a solution taking values in
$(C^+_{\mathrm{rap}})^{4N_{\ep}}$ to the system of
above equations can be constructed via standard steps by extending the
procedure in~\citet{shi94}. We will comment further on this point below.

We also define the following nondecreasing $M_F(\R)$-valued processes:
\begin{eqnarray*}
K_t^{i,n,U}(\phi)&=& n \int_0^t
\bigl\langle U^{i,n}_s V^{n}_s,
\phi\bigr\rangle \,ds, \qquad t\geq0, \phi \in C_{b}(\R),
\\
K_t^{j,n,V}(\phi)&=& n \int_0^t
\bigl\langle V^{j,n}_s U^{n}_s
\phi,\bigr\rangle \,ds, \qquad t\geq0, \phi \in C_{b}(\R).
\end{eqnarray*}
Clearly,
\[
\sum_{i\in\NN_{\ep}} K_t^{i,n,U}= \sum
_{j\in\NN_{\ep}} K_t^{j,n,V}=:
K_t^n,
\]
and $(U^n,V^n,\tU^n,\tV^n)$ satisfies the following system of equations
for $\phi\in C_b^2(\R)$:
%
%e9.2 #&#
\begin{eqnarray}
\label{eq43} \cases{ %
\displaystyle U^{n}_t(
\phi) = \sum_{i\in\NN_{\ep}}\bigl\langle J^{x_i},
\phi\bigr\rangle \1 (t\geq s_i)
\vspace*{2pt}\cr
\hspace*{42pt}\displaystyle{} + \int_0^t \int_{\R}U^n(s,x)^{\gamma}
\phi(x) W^{n,U}(ds,dx)
\vspace*{2pt}\cr
\hspace*{42pt}\displaystyle{}+ \int_0^t U^{n}_s
\biggl(\frac{1}{2}\Delta\phi\biggr) \,ds - K^n_t(\phi
),\qquad t\geq0,
\vspace*{2pt}\cr
\displaystyle V^{n}_t(\phi) = \sum_{j\in\NN_{\ep}}
\bigl\langle J^{y_j},\phi\bigr\rangle \1 (t\geq t_j)
\vspace*{2pt}\cr
\hspace*{42pt}\displaystyle{}+ \int_0^t\int_{\R}
V^n(s,x)^{\gamma} \phi(x) W^{n,V}(ds,dx)
\vspace*{2pt}\cr
\hspace*{42pt}\displaystyle{}+ \int_0^t V^{n}_s
\biggl(\frac{1}{2}\Delta\phi\biggr) \,ds - K^n_t(\phi
),\qquad t\geq0,
\vspace*{2pt}\cr
\displaystyle\tU^{n}_t(\phi) = \int_0^t
\int_{\R} \bigl[ \bigl(\tU ^n(s,x)+U^n(s,x)
\bigr)^{2\gamma} - U^n(s,x)^{2\gamma}
\bigr]^{1/2}
\displaystyle\vspace*{2pt}\cr\hspace*{67pt}{}\times\phi(x) \tW^{n,U}(ds,dx)
\vspace*{2pt}\cr\hspace*{42pt}\displaystyle{}+ \int_0^t \tU^{n}_s
\biggl(\frac{1}{2}\Delta\phi\biggr) \,ds + K^{n}_t(
\phi), \qquad t\geq0,
\vspace*{2pt}\cr
\displaystyle\tV^{n}_t(\phi) = \int_0^t
\int_{\R} \bigl[ \bigl(\tV ^n(s,x)+V^n(s,x)
\bigr)^{2\gamma} - V^n(s,x)^{2\gamma}
\bigr]^{1/2}
\vspace*{2pt}\cr{}\hspace*{67pt}\displaystyle\times\phi(x) \tW^{n,V}(ds,dx)
\vspace*{2pt}\cr{}\hspace*{42pt}\displaystyle+ \int_0^t \tV^{n}_s
\biggl(\frac{1}{2}\Delta\phi\biggr) \,ds + K^{n}_t(
\phi), \qquad t\geq0,}
\end{eqnarray}
with
$W^{n,U}, W^{n,V}, \tW^{n,U}, \tW^{n,V}$ being a collection of
independent space--time white noises.
For $i\in\NN_{\ep}$,
define $\bU^{i,n}_t\equiv U^{i,n}_t+\tU^{i,n}_t, \bV^{i,n}_t
\equiv
V^{i,n}_t+\tV^{i,n}_t, t\geq0$ and
%
%e9.3 #&#
\begin{eqnarray}
\label{11081} \bU^{n}_t\equiv\sum
_i \bU^{i,n}_t,\qquad \bV^{n}_t
\equiv\sum_j \bV^{j,n}_t,\qquad  t
\geq0.
\end{eqnarray}
Since $\{ W^{i,n, U}, W^{j,n,V},\tW^{k,n, U}, \tW^{l,n,V}, i,j,k,l\in
\NN_{\ep}\}$ is a collection
of independent white noises, and by stochastic calculus, one can easily show
that the
processes $\bU^{n}, \bV^{n}$ satisfy equations~(\ref{eq28}), and
so by
\citet{myt98w} they have laws on $D([0,T],C^+_{\mathrm{rap}})$ which are
independent of $n$.

Here we comment further on the construction of
$(U^{i,n}, V^{i,n}, \tU^{i,n}, \tV^{i,n})_{i\in\N_{\ep}}$, the solution
to~(\ref{tUVndefn}). As we have mentioned above, one can follow the procedure
indicated in the proof of Theorem~2.6 in~\citet{shi94}
by extending it to systems of equations.
In the proof, one constructs an approximating sequence of processes
$\{(U^{i,n,k}, V^{i,n,k}, \tU^{i,n,k}, \tV^{i,n,k})_{i\in\N_{\ep
}}\}
_{k\geq1}$
with
globally Lipschitz coefficients, and shows that this sequence is tight in
\[
\prod_{i=1}^{\N_{\ep}} \bigl(C
\bigl([s_i,\infty), C^+_{\mathrm{rap}}\bigr)\times C
\bigl([t_i,\infty), C^+_{\mathrm{rap}}\bigr)\times C
\bigl([s_i,\infty), C^+_{\mathrm{rap}}\bigr)\times C
\bigl([t_i,\infty), C^+_{\mathrm{rap}}\bigr)\bigr ),
\]
and each limit point satisfies~(\ref{tUVndefn}).
The only subtle point is that the drift coefficients $U^{i,n}(\cdot
)V^n(\cdot)$ and
$V^{i,n}(\cdot)U^n(\cdot)$
in the system of limiting equations~(\ref{tUVndefn})
do not satisfy a linear growth condition. However, note that, by \eqref
{11081}, any solution to~(\ref{tUVndefn})
satisfies the following bounds:
%
%e9.4 #&#
\begin{eqnarray}
\label{11082} U^{i,n}, \tU^{i,n},U^{n},
\tU^n\leq\bU^n,\qquad V^{i,n}, \tV^{i,n},V^{n},
\tV^n\leq\bV^n,
\end{eqnarray}
where $\bU^n$ and $\bV^n$ have
good moment bounds by Lemma~\ref{pmom}. Hence, it is possible to construct
$\{(U^{i,n,k}, V^{i,n,k}, \tU^{i,n,k}, \tV^{i,n,k})_{i\in\N_{\ep
}}\}
_{k\geq1}$
so that
the bound in Lemma~\ref{pmom} holds
uniformly in $k$: for any $q,T>0$, there exists $C_{q,T}$ such that
\begin{eqnarray*}
&&\hspace*{-4pt}\sup_{k\ge1}\sup_{i\in\N_{\ep}}E \Bigl[\sup
_{s\le T,x\in\R
}\bigl(U^{i,n,k}(s,x)^{q} +
\tU^{i,n,k}(s,x)^{q} +V^{i,n,k}(s,x)^{q} +
\tU^{i,n,k}(s,x)^{q}\bigr) \Bigr]
\\
&&\hspace*{-4pt}\qquad\le C_{q,T}.
\end{eqnarray*}
With this uniform bound in hand, it is not difficult to check that the
moment bound~(6.5) from~\citet{shi94} [which is
in fact \eqref{lem-X-replace-by-N} with $\lambda=0$],
holds for $\{U^{i,n,k}\}_{k\geq1}, \{V^{i,n,k}\}_{k\geq1}$, $\{ \tU
^{i,n,k}\}_{k\geq1}$, $ \{\tV^{i,n,k}\}_{k\geq1}$,
for all $i\in\N_{\ep}$, on time intervals of the form $[\frac
{(i-1)\ep
}{2}, \frac{i\ep}{2}), i\in\N_{\ep}$ and $[N_{\ep}\ep, T]$.
This, in turn,
by Lemma~6.3 in~\citet{shi94} implies the tightness of the corresponding
processes in $D^{\ep}(\R_+, C^+_{\mathrm{tem}})$. Here
\[
C_{\mathrm{tem}}:=\bigl\{f \in C(\R)\dvtx \Vert f\Vert _{\lambda} <
\infty\mbox{ for any } \lambda<0\bigr\},
\]
endowed with the topology induced
by the norms $\Vert \cdot\Vert _{\lambda}$ for $\lambda<0$, and $C^+_{\mathrm{tem}}$ is the set of nonnegative functions in $C_{\mathrm{tem}}$.
Finally, since the limiting processes $U^{i,n}$, $\tU^{i,n}, i\in\N
_{\ep}$, (resp., $V^{i,n}$, $\tV^{i,n}, i\in\N_{\ep}$) are
dominated by
$\bU$ (resp., $\bV$) in $D^{\ep}(\R_+, C^+_{\mathrm{rap}})$,
it follows that $U^{i,n}$, $\tU^{i,n}$, $V^{i,n}$, $\tV^{i,n}, i \in
\N
_{\ep}$,
are in $D^{\ep}(\R_+, C^+_{\mathrm{rap}})$ as well. This, together with the
domination~\eqref{11082} and Lemma~\ref{qmom}, allows us to take
functions in $C^2_{\mathrm{tem}}$ as test functions in \eqref{tUVndefn};
however for our purposes it will be enough to use functions from
$C_b^2(\R)$ as test functions.

Fix an arbitrary $T>1$.
%
%re9.1 #&#
\begin{remark}
\label{rem04}
In what follows we are going to show the tightness of the sequence of
the processes constructed above on the
time interval $[0,T]$. We will prove that limit points have the
properties stated in Proposition~\ref{thm11} on\vadjust{\goodbreak}
$[0,T]$.
Since $T>1$ is arbitrary, this argument immediately yields the claim of
the theorem on the time interval $[0,\infty)$.
\end{remark}

Define $E=[0,T]\times\R$. We identify a finite measure $K$ on $E$ with
the nondecreasing path in $D([0,T],M_F(\R))$ given by $t\to K_t(\cdot
)=K([0,t]\times\{\cdot\})$.

%pr9.2 #&#
\begin{prop}
\label{prop4}
$\{(U^{i,n},\tU^{i,n},V^{i,n}, \tV^{i,n}, K^{i,n,U},K^{i,n,V})_{i\in
\NN
_{\ep}}\}_{n\geq1}$
is tight in $ (C([0,T]\setminus\CG_{\ep},M_F(\R))^4\times
M_F(E)^2 )^{N_{\ep}}$.
Moreover, any limit point $(U^{i},\tU^{i},V^{i},\break  \tV^{i}, K^{i,U},
K^{i,V})_{i\in\NN_{\ep}}$
has the following properties:
\begin{longlist}[(1)]
\item[(1)] $U^i, \tU^i, V^i, \tV^i\in C([0,T]\setminus\CG_{\ep},
C^+_{\mathrm{rap}})\cap D^{\ep}([0,T], L^1(\R)),
\forall i\in\NN_{\ep}$;
%with respect to the Lebesgue measure;
%$U^i, \tU^i, V^j, \tV^j\in D([0,2]\setminus\CG_{\ep}, L^2(\R)),
%
\item[(2)] $K^{i,U}, K^{i,V}\in D^{\ep}([0,T], M_F(\R)), \forall
i\in
\NN_{\ep}$;
\item[(3)] $(U^i, \tU^i, V^i, \tV^i,K^{i,U},K^{i,V})_{i\in\NN_{\ep}}$
satisfy \eqref{UVdefn}--\eqref{tUVdefn}.
\end{longlist}
\end{prop}

The above proposition is the key for proving Proposition~\ref{thm11}.
The proposition will be proved via a series of lemmas.

%le9.3 #&#
\begin{lemma}
\label{lem0171}
$\{K^n\}_{n\geq1}$ is tight in $M_F(E)$, and $\{K^n_T(1)\}_{n\ge1}$
is $L^1(dP)$-bounded.
\end{lemma}
%
%pa9.subsection.subsubsection.1 #&#
\begin{pf}
First note that by rewriting equation~(\ref{eq28}) for $\bU^n$ in the
mild form [see \eqref{eqmild1}] one can
easily get that for any $\phi\in C_{b}^+(\R)$,
\begin{eqnarray}
\nonumber
E \bigl[\bU^n_t(\phi) \bigr]&\leq& E \biggl[
\sum_{s_i\in\CG_{\ep}^{\mathrm{odd}}, s_i\leq t} \int_{\R}\int
_{\R} p_{t-s_i}(z-y) J^{x_i}_{\ep}(y)
\phi(z) \,dy \,dz \biggr]
\\
&=& \sum_{s_i\in\CG_{\ep}^{\mathrm{odd}}, s_i\leq t} \int_0^1
\int_{\R
}\int_{\R} p_{t-s_i}(z-y)
J^{x}_{\ep}(y)\phi(z) \,dy \,dz \,dx
\\
\nonumber
&=& \sum_{s_i\in\CG_{\ep}^{\mathrm{odd}}, s_i\leq t} \int_{\R}
\int_{\R}\int_{0}^1
p_{t-s_i}(z-y) J^{x}_{\ep}(y)\phi (z) \,dx \,dz \,dy.
\nonumber
\end{eqnarray}
Estimating the above integrals, we have
\begin{eqnarray*}
E \bigl[\bU^n_t(\phi) \bigr] &\leq& \ep\sum
_{s_i\in\CG_{\ep}^{\mathrm{odd}}, s_i\leq t} \int_{\R
}S_{t-s_i}\phi(y)
\1\bigl(|y|\leq2\bigr) \,dy
\\
\label{1101} &\leq& \sup_{s\leq t}\int_{\R}
S_{s}\phi(y) \1\bigl(|y|\leq2\bigr) \,dy,
\end{eqnarray*}
where $\{S_t\}_{t\geq0}$ is the Brownian semigroup corresponding to
the transition density function $\{ p_t(x), t\geq0,
x\in\R\}$.

For any nonnegative $\phi\in{C^2_{b}}(\R)$ we have from \eqref{eq43},
%
%e9.5 #&#
\begin{eqnarray}\label{eq0671}
\nonumber
E \bigl[K^n_t(\phi) \bigr] &\leq& E \biggl[
\sum_{i\in\CG_{\ep
}^{\mathrm{odd}}} \int_{\R}
J^{x_i}_{\ep}(y) \phi(y) \,dy \biggr] + E \biggl[ \int
_0^t U^n_s \biggl(
\biggl\llvert \frac
{\Delta
\phi}{2}\biggr\rrvert \biggr) \,ds \biggr]
\\
&\leq& \sum_{i\in\CG_{\ep}^{\mathrm{odd}}} \int_0^1
\int_{\R} J^{x}_{\ep}(y) \phi(y) \,dy \,dx +
E \biggl[ \int_0^t \bU^n_s
\biggl(\biggl\llvert \frac
{\Delta\phi
}{2}\biggr\rrvert \biggr) \,ds \biggr]
\\
\nonumber
&\leq& \int_{\R} \1\bigl(|y|\leq2\bigr) \phi(y) \,dy
+ \int_0^t \sup
_{r\leq s}\int_{\R} S_{r} \biggl(
\biggl\llvert \frac{\Delta\phi}{2}\biggr\rrvert \biggr) (y) \1\bigl(|y|\leq 2\bigr) \,dy \,ds.
\end{eqnarray}
Now by taking $\phi=1$ we get that the sequence of the total masses $\{
K^n_T(1)\}_{n\geq1}$ is
bounded in $L^1(dP)$. Moreover for any $\delta>0$ we can choose
$R>3$ sufficiently large and $\phi$ such that $\phi(z)=0$ for
$|z|\leq
R-1, \phi(z)=1$ for $|z|\ge R$ with the
property that
\[
S_{t} \biggl(\biggl\llvert \frac{\Delta\phi}{2}\biggr\rrvert \biggr)
(y)\leq \delta\qquad \forall t\in[0,T], y\in[-2, 2].
\]
This shows that
\[
E \biggl[ \int_{|z|\geq R} K^n_T(dz)
\biggr]\leq E \bigl[K^n_T(\phi ) \bigr] \leq4T\delta\qquad
\forall n\geq1,
\]
by~(\ref{eq0671}), and our choice of $\phi$ and $R$. This, in turn,
together with the
$L^1(dP)$-boundedness of
total masses $\{K^n_T(1)\}_{n\geq1}$, implies tightness of $\{K^n\}
_{n\geq1}$ in $M_{F}(E)$.
\end{pf}

%co9.4 #&#
\begin{cor}
\label{cor151}
$\{K^{i,n,U}\}_{n\geq1}$ and $\{K^{i,n,V}\}_{n\geq1}$ are tight in
$M_F(E)$ for any $i\in\NN_\ep$.
\end{cor}
%
%pa9.subsection.subsubsection.2 #&#
\begin{pf}
The assertion follows immediately from the bound
\[
K^{n,i,U}, K^{n,i,V} \leq K^n \qquad \forall n\geq1, i
\in\NN_{\ep}.
\]
\upqed\end{pf}\eject

Before we start dealing with tightness of $\{(U^n,V^n,\tU^n,\tV
^n,K^n)\}
_{n\geq1}$ we need to introduce a lemma
that will be frequently used.
%
%le9.5 #&#
\begin{lemma}
\label{lem3061}
We have:
\begin{longlist}[(a)]
\item[(a)] Let $\{W^n\}_{n\geq1}$ be a sequence of $\{\CF^n_t\}
_{t\geq0}$-adapted space--time white noises, and
$\{b^n(t,x,\omega)\}_{n\geq1}$ be a sequence of
$\{\CF^n_t\}_{t\geq0}$-predictable${}\times{}$Borel measurable processes
such that
%
%e9.6 #&#
\begin{eqnarray}
\sup_{n\geq1} \sup_{x\in\R}\sup
_{t\in[0,T]} E \bigl[ \bigl|b^{n}(t,x,\cdot )\bigr|^p
\bigr]<\infty\qquad \mbox{for some $p>4$}.
\end{eqnarray}
Then the sequence of processes $\{X^n(t,x),   t\in[0,T], x\in\R\}
_{n\geq1}$ defined by
\[
X^n(t,x) = \int_0^t \int
_{\R} p_{t-s}(x-y) b^n(s,y,\cdot)
W^n(ds,dy),\qquad t\in[0,T], x\in\R,
\]
have versions which are tight in $C([0,T], C_{\mathrm{tem}})$.
\item[(b)] Let $W$ be an $\{\CF_t\}_{t\geq0}$-adapted space--time
white noise, and
$b(t,x,\omega)$ be an $\{\CF_t\}_{t\geq0}$-predictable${}\times{}$Borel
measurable process such that
%
%e9.7 #&#
\begin{eqnarray}
\sup_{x\in\R}\sup_{t\in[0,T]} E \bigl[ \bigl|b(t,x,
\cdot)\bigr|^p \bigr]<\infty \qquad \mbox{for some $p>4$}.
\end{eqnarray}
Then the process $X$ defined by
\[
X(t,x) = \int_0^t \int_{\R}
p_{t-s}(x-y) b(s,y,\cdot) W^n(ds,dy),\qquad t\in[0,T], x\in\R,
\]
has a version in $C([0,T], C_{\mathrm{tem}})$. If moreover, $|X(t,x)|\leq
|\widetilde X(t,x)|$ for some $\widetilde X\in
D([0,T],C_{\mathrm{rap}})$, then $X\in C([0,T], C_{\mathrm{rap}})$.
\end{longlist}
\end{lemma}
%
%pa9.subsection.subsubsection.3 #&#
\begin{pf}
(a) This assertion follows immediately from the estimates on
increments of a stochastic integral [see, e.g., step 2
in the proof of Theorem~2.2 of~\citet{shi94}, page 432] and then an
application of Lemmas~6.2 and 6.3(ii) from~\citet{shi94}.

(b) This again follows by using the estimates on
increments of a stochastic integral [see again step 2
in the proof of Theorem~2.2 of~\citet{shi94}, page 432] and then
applying Lemmas~6.2 and 6.3(i) in~\citet{shi94},
to get that the process is in $C([0,T], C_{\mathrm{tem}})$. The last
assertion is obvious.
\end{pf}
%
%le9.6 #&#
\begin{lemma}
\label{lem0172}
Let
\[
w^{n}= U^n-V^n,\qquad n\geq1.
\]
Then $\{w^n\}_{n\geq1}$ is tight in $D([0,T], C_{\mathrm{rap}})$, and every
limit point is in
$D^{\ep}([0,T],\break   C_{\mathrm{rap}})$.
\end{lemma}
%
%pa9.subsection.subsubsection.4 #&#
\begin{pf}
By writing the equation for $w^n$ in mild form we get
\begin{eqnarray*}
w^n(t,x)&=& \int_0^t \int
_{\R} p_{t-s}(x-y) \bigl(\eta^+_{\ep
}(ds,dy)-
\eta ^-_{\ep}(ds,dy)\bigr)
\\
&&{} +\int_0^t\int_{\R}
p_{t-s}(x-y) U^n(s,y)^{\gamma} W^{n,U}(ds,dy)
\\
&&{}- \int_0^t\int
_{\R} p_{t-s}(x-y) V^n(s,y)^{\gamma}
W^{n,V}(ds,dy),\qquad t\geq0, x\in\R.
\end{eqnarray*}
Clearly, by the definition of $\eta^+_{\ep}, \eta^-_{\ep}$, the first
term, $I(t,x)$ (being independent of~$n$) is tight in $D([0,T], C_{\mathrm{rap}})$, and is in
$D^{\ep}([0,T], C_{\mathrm{rap}})$.
%is tight in $C([s_i, t_i), C_{\mathrm{rap}})$ and
%$C([t_i, s_{i+1}), C_{\mathrm{rap}})$, for any $i\in\NN_{\ep}$, and hence it
%is tight in $D^{\ep}([0,T], C_{\mathrm{rap}})$.
Using the domination
%
%e9.8 #&#
\begin{eqnarray}
\label{eq0171} U^n\leq\bU^n\in D\bigl([0,T],C^+_{\mathrm{rap}}
\bigr), \qquad V^n \leq\bV^n\in D\bigl([0,T],C^+_{\mathrm{rap}}
\bigr),
\end{eqnarray}
and Lemmas~\ref{pmom} and \ref{lem3061}(a), the stochastic integral
terms are tight in $C([0,T],\break   C_{\mathrm{tem}})$. If $S^n(t,x)$ is the
difference of the above stochastic integral terms, then
the domination
\[
\bigl|S^n(t,x)\bigr|\le\bar U^n(t,x)+\bar V^n(t,x)+\bigl|I(t,x)\bigr|
\in D^\vep \bigl([0,T],C^+_{\mathrm{rap}}\bigr),
\]
and the definition of the norms on $C_{\mathrm{tem}}$ and $C_{\mathrm{rap}}$
shows that
$\{S^n\}$ is tight in $C([0,T],C_{\mathrm{rap}})$.
%By tightness of the first term in $D^{\ep}([0,T], C_{\mathrm{rap}})$ and
%the domination
%$$ |w^n| \leq|\bU|+|\bV|\in D([0,T],C_{\mathrm{rap}}),$$
%the tightness of $\{w^n\}_{n\geq1}$ in $D^{\ep}([0,T], C_{\mathrm{rap}})$
%follows.
\end{pf}

Now we are ready to deal with the tightness of $\{(U^n,V^n,\tU^n,\tV
^n,K^n)\}_{n\geq1}$.
Let $L^p(E)$ denote the usual $L^p$ space with respect to Lebesgue
measure on~$E$.

%le9.7 #&#
\begin{lemma}
\label{lem151}
The following assertions hold:
\begin{longlist}[(a)]
\item[(a)]
$\{(U^n,V^n,\tU^n,\tV^n,K^n)\}_{n\geq1}$ is tight in
$L^p(E)^4\times M_F(E)$ for any
$p\geq1$.
Moreover any limit point has a version
\[
(U,V,\tU,\tV,K)\in D^{\ep}\bigl([0,T], C^+_{\mathrm{rap}}
\bigr)^4\times D^{\ep}\bigl([0,T], M_F(\R)
\bigr).
\]
\item[(b)]
\[
t\mapsto\int_0^{t}\int_{\R}
p_{t-s}(\cdot-y)K(ds,dy)\in D^{\ep
}\bigl([0,T],
C_{\mathrm{rap}}\bigr).
\]
\item[(c)] $\{K^n\}_{n\geq1}$ is also tight in
$C([0,T]\setminus
\CG_{\ep}, M_F(\R))$, and any of its limit points
satisfies
\[
{\bolds\Delta} K_t(1)\leq\ep \qquad\forall t\in[0,T].
\]
\end{longlist}
\end{lemma}
%
%pa9.subsection.subsubsection.5 #&#
\begin{pf}
(a) We will give the proof just for the tightness\vspace*{1pt} of $\{(U^n,V^n,\break K^n)\}
_{n\geq1}$ and the properties of its limit points,
since the corresponding results for
$\{(\tU^n, \tV^n)\}_{n\geq1}$ and its limit points will follow along
the same lines.

Recall the domination (\ref{eq0171}), where the laws of the upper
bounds are independent of $n$. By this domination we immediately get that
\[
\bigl\{\bigl(U^{n}(s,x)\,dx\,ds, V^{n}(s,x)\,dx\,ds\bigr)\bigr
\}_{n\geq1}
\]
is
tight in $(M_F(E)\times M_F(E))$. Recall also that by Lemma~\ref
{lem0171}, $\{K^n\}_{n\geq1}$ is tight in
$M_F(E)$.
%Also define
%$$w^{n}= U^n-V^n.$$
%Since the killing term is the same for both $U^n$ and $V^n$ the drift
%does not appear in the
%equation for $w^n$. Hence by standard argument one can show the
%tightness in $D_{C_{\mathrm{rap}}}$ for
%$\{w^n\}_{n\geq1}$ ($C$-tight in $D[s_i,t_{i+1})$ for all $i$?).
This, the fact that the laws of $\bar U_n$, $\bar V_n$ are independent
of $n$, and Lemma~\ref{lem0172} allows us to choose a convergent
subsequence of $(U^n,V^n,K^n,w^n, \bar U^n,\bar V^n)$ in
$M_F(E)^3\times D([0,T],C_{\mathrm{rap}})^3$.
For simplicity of notation, we will again index this subsequence by
$n$. Denote the corresponding limit point by $(U,V,K,w,\bar U,\bar V)$.

Now, for any $\phi\in C_{b}(\R)$, let
\begin{eqnarray*}
M^{n,U}_t(\phi)&\equiv& \int_0^t
\int_{\R}U^n(s,x)^{\gamma} \phi(x)
W^{n,U}(ds,dx),\qquad t\in[0,T],
\\
M^{n,V}_t(\phi)&\equiv& \int_0^t
\int_{\R}V^n(s,x)^{\gamma} \phi(x)
W^{n,V}(ds,dx),\qquad t\in[0,T],
\end{eqnarray*}
denote the martingales given by the stochastic integrals in the
semimartingale decomposition~(\ref{eq43})
for $U^n_t(\phi)$ and $V^n_t(\phi)$.
For any $\phi\in C_{b}(\R)$, use the Burkholder--Davis--Gundy
inequality, and
again the domination~(\ref{eq0171}),
to get, that for any $p\geq2, \lambda>0$,
%
%e9.9 #&#
\begin{eqnarray}
\label{eq0172}
&& E \bigl[ \bigl\llvert M^{n,U}_{t}(
\phi)-M^{n,U}_{u}(\phi)\bigr\rrvert ^p \bigr]
\nonumber\\
&&\qquad\leq C_p \sup_{s\leq T, x\in\R}e^{({\lambda p}/{2}) |x|}E
\bigl[ \bU (s,x)^{p\gamma} \bigr]
\nonumber
\\[-8pt]
\\[-8pt]
\nonumber
&&\hspace*{82pt}{}\times \biggl[ \int_{\R}
e^{-\lambda|x|}\bigl\llvert \phi(x) \bigr\rrvert ^2 \,dx
\biggr]^{p/2} (t-u)^{p/2},
\\
 \eqntext{\forall0\leq u\leq t\leq T.}
\end{eqnarray}
This, together with Lemma~\ref{qmom}(b) and Kolmogorov's tightness
criterion, implies that
%
%e9.10 #&#
\begin{eqnarray}
\bigl\{M^{n,U}_{\cdot}(\phi) \bigr\}_{n\geq1} \mbox{is
tight in } C\bigl([0,T],\R\bigr)
\end{eqnarray}
for any $\phi\in C_{b}(\R)$. Similarly,
%
%e9.11 #&#
\begin{eqnarray}
\bigl\{M^{n,V}_{\cdot}(\phi) \bigr\}_{n\geq1} \mbox{is
tight in } C\bigl([0,T],\R\bigr)
\end{eqnarray}
for any $\phi\in C_{b}(\R)$.
Let ${\cal D}$ be a countable subset of
$C^2_{b}(\R)$ which is bounded-pointwise dense in $C_b(\R)$. That is,
the smallest class containing ${\cal D}$ and closed under bounded
pointwise limits contains $C_b(\R)$. By the above, we can take a
further subsequence, which for simplicity we will index again by
$n$, so that
all the sequences of martingales $ \{M^{n,U}_{\cdot}(\phi)
\}
_{n\geq1}$,
$ \{M^{n,V}_{\cdot}(\phi) \}_{n\geq1}$
indexed by functions $\phi$ from ${\cal D}$, converge in $C([0,T],\R)$.
For $\phi\in{\cal D}$, we will denote the limiting processes
by $M^U_{\cdot}(\phi), M^V_{\cdot}(\phi)$, respectively.
Now let us switch to a probability space where
%
%e9.12 #&#
\begin{eqnarray}
\label{asconv}
\bigl(U^n,V^n,K^n,w^n,
\bar U^n,\bar V^n\bigr)&\rightarrow&(U,V,K,w,\bar U,\bar
V)\nonumber\\
 \eqntext{\mbox{in } M_F(E)^3\times D\bigl([0,T],C_{\mathrm{rap}}
\bigr)^3,}
\\[-8pt]
\\[-8pt]
\nonumber
\bigl(M^{n,U}(\phi_1),M^{n,V}(
\phi_2)\bigr)&\rightarrow&\bigl(M^{U}(\phi
_1),M^{V}(\phi_2)\bigr)\\
\eqntext{\mbox{in } C\bigl([0,T],
\R\bigr)^2\ \forall\phi_1,\phi_2\in{\cal D},}
\end{eqnarray}
as $n\rightarrow\infty$, a.s.

In our next step, we will verify convergence of $\{(U^n, V^n)\}_{n\geq
1}$ in $L^p(E)^2$, for any $p\geq1$. First,
by $L^1(dP)$-boundedness of the total mass of $K^n$ (Lem\-ma~\ref
{lem0171}), we have
%
%e9.13 #&#
\begin{eqnarray}
n E \biggl[ \int_0^T \int
_{\R} U^n_s(x)V^n_s(x)\,dx\,ds
\biggr]=E \bigl[K^n_T(1) \bigr] \leq C,
\end{eqnarray}
uniformly in $n$ for some constant $C$. Therefore we get
%
%e9.14 #&#
\begin{eqnarray}
E \biggl[ \int_0^T \int_{\R}
U^n_s(x)V^n_s(x)\,dx\,ds
\biggr] \rightarrow0\qquad\mbox{as } n\rightarrow\infty,
\end{eqnarray}
and hence
%
%e9.15 #&#
\begin{eqnarray}
\int_0^T \int_{\R}
\bigl(U^n_s(x)\wedge V^n_s(x)
\bigr)^2 \,dx\,ds \rightarrow0,
\end{eqnarray}
in $L^1(dP)$.
%Therefore, with almost sure convergence of $\{(U^n,V^n)\}_{n\geq1}$ in
%$M_F(E)^2$, this implies
By taking another subsequence if necessary, we may assume
\[
\bigl(U^n_s(x)\wedge V^n_s(x)
\bigr) \rightarrow0 \qquad\mbox{in } L^2(E), P\mbox{-a.s.}
\]
Now recall again the domination
\[
U^n \leq\bU^n\to\bU \qquad\mbox{in }D\bigl([0,T],
C^+_{\mathrm{rap}}\bigr), P\mbox{-a.s.},
\]
which implies that for any $p\geq1$,
\[
\bigl(U^n_s(x)\wedge V^n_s(x)
\bigr) \rightarrow0 \qquad\mbox{in } L^p(E), P\mbox{-a.s.}
\]
Also by
\[
U^{n}_t(x)= \bigl(U^n_s(x)
\wedge V^n_s(x)\bigr) + \bigl(w^n_t(x)
\bigr)^+,
\]
we get that in fact
%
%e9.16 #&#
\begin{eqnarray}
\label{genpc} U^n\rightarrow(w)^+\qquad \mbox{in } L^p(E),
\mbox{ for any $p \geq1$,} P\mbox{-a.s.},
\end{eqnarray}
and hence $U(dt,dx)=w_t(x)^+\, dt\,dx$. With some abuse of notation we
denote the density of $U(dt,dx)$
by $U_t(x)$. Similarly we get
\[
V(dt,dx) = w_t(x)^- \,dt\,dx,
\]
and we denote its density by $V_t(x)$. In what follows we will use the
continuous in space versions of the densities of
$U(dt,dx),V(dt,dx)$, that is, $U_t(x)=w_t(x)^+,   V_t(x)=w_t(x)^-$,
and hence, by Lemma~\ref{lem0172},
we get that $(U,V)\in D^{\ep}([0,T], C_{\mathrm{rap}})^2$. We delay the
proof of the assertion that $K\in D^{\ep}([0,T],M_F(\R))$ until
the proof of part {(b)}.

%pa9.subsection.subsubsection.6 #&#
(b) Fix an arbitrary $\phi\in{\cal D}$.
We will go to the limit in (\ref{eq43}) for $\{ U^n_{\cdot}(\phi)\}
_{n\geq1}$.
As $\{U^n\}_{n\geq1}$ converges a.s. to $w^+$ in $L^2(ds,dx)$, and
\[
U^n\leq\bU^n\to\bU\qquad \mbox{in } D\bigl([0,T],
C_{\mathrm{rap}}\bigr),
\]
it is easy to see that $\{U^n_{\cdot}(\phi)\}_{n\geq1}$
converges to
$w^+_{\cdot}(\phi)\equiv\int w^+_{\cdot}(x)\phi(x) \,dx$ in $L^2[0,T]$
a.s. As for the right-hand side, use~(\ref{genpc}) with $p=1$ to get
\[
\sup_{t\le T} \biggl|\int_0^tU_s^n
\biggl(\frac{1}{2}\Delta\phi\biggr) \,ds-\int_0^tU_s
\biggl(\frac{1}{2}\Delta\phi\biggr) \,ds \biggr|\le\bigl\| U^n-U
\bigr\|_{L^1(E)}\| \Delta \phi/2\|_{\infty}\to0.
\]
In particular this implies that
$\{
\int_0^{\cdot} U^{n}_s(\frac{1}{2}\Delta\phi)  \,ds\}_{n\geq1}$
converges to
$
\int_0^{\cdot} U_s(\frac{1}{2}\Delta\phi)  \,ds$
in $C([0,T],\R)$
(and hence in $L^2[0,T]$).
By (a) $\{K^n(\phi)(ds)\}_{n\geq1}$ converges to $K(\phi)(ds)$ as
finite signed measures on $[0,T]$ a.s., and therefore
$\{K^n_{\cdot}(\phi)\}_{n\geq1}$ converges in $L^2[0,T]$ to
$K_{\cdot
}(\phi)$ a.s. Since the immigration term
does not change with~$n$,
it also converges in $L^2[0,T]$.

Now we have to deal with convergence of the stochastic integral term,
that we denoted by $M^{n,U}(\phi)$.
We proved in~(a) that $\{M^{n,U}(\phi)\}_{n\geq1}$
converges a.s. in $C([0,T],\R)$.
Moreover, by~(\ref{eq0172}), the martingales $M_t^{n,U}(\phi)$ are
bounded in $L^p(dP)$ uniformly in $n$ and $t\in[0,T]$,
for all $p\geq2$, and hence the limiting process is a continuous
martingale that we will call $M^{U}(\phi)$. Turning to its quadratic variation,
%note that as $U^{i,n}, U^{n}\leq\bU\in C_{\mathrm{rap}}$, then the
it follows from \eqref{genpc}
that the sequence $\{ (U^n)^{2\gamma}\}_{n\geq1}$ converges
to $U^{2\gamma}$ in $L^2(E)$ a.s.
and this implies that
%for any $\phi\in C_{\mathrm{rap}}(\R)$,
%
%e9.17 #&#
\begin{eqnarray}
\label{eq0173}
\bigl\langle M^{n,U}_{\cdot}(\phi)
\bigr\rangle_t &=& \int_0^t \int
_{\R} U^n(s,x)^{2\gamma}
\phi(x)^2 \,dx\,ds
\nonumber
\\[-8pt]
\\[-8pt]
\nonumber
&\rightarrow& \int_0^t \int
_{\R}U(s,x)^{2\gamma} \phi(x)^2 \,dx\,ds\qquad
\mbox{as } n\rightarrow\infty, P\mbox{-a.s.}
\end{eqnarray}
Hence, again by boundedness of $M_t^{n,U}(\phi)$ in $L^p(dP), p\geq2$,
uniformly in $t\in[0,T], n\geq1$,
we get that the limiting continuous
martingale $M^{U}$ has quadratic variation
\[
\bigl\langle M^{U}_{\cdot}(\phi)\bigr\rangle_t=
\int_0^t \int_{\R
}U(s,x)^{2\gamma
}
\phi(x)^2 \,dx\,ds
\]
for any $\phi\in{\cal D}$. Since ${\cal D}$ is bounded-pointwise dense
in $C_{b}(\R)$, $M^U$ can be extended to a martingale
measure on $E$, and
one can show by\vadjust{\goodbreak} standard procedure
that there is a space--time white noise $W^U$ such that
\[
M^U_t(\phi)=\int_0^t
\int_{\R} U(s,x)^{\gamma} \phi(x) W^U(ds,dx),\qquad
t\in[0,T], \forall\phi\in C_{b}(\R).
\]

Now we are ready to take limits in \eqref{eq43} in $L^2([0,T])$.
%Recalling that $U=w^+$
We get
%
%e9.18 #&#
\begin{eqnarray}\label{eq44}
\nonumber
U_t(\phi) &=& \sum
_i\bigl\langle J^{x_i},\phi\bigr\rangle\1(t\geq
s_i)
\\
&& {}+ \int_0^t \int
_{\R}U(s,x)^{\gamma} \phi(x) W^{U}(ds,dx)
\\
&& {}+ \int_0^t U_s\biggl(
\frac{1}{2}\Delta\phi\biggr) \,ds - K_t(\phi), \qquad t\in[0,T].
%
%&&\mbox{}\\
%w^-_t(\phi) &=& \sum_j\langle J^{y_j},\phi\rangle\1(t\geq t_j) \\
%&& + \int_0^t\int_{\R} V(s,x)^{\gamma} \phi(x) W^{V}(ds,dx) \\
%&& \mbox{}+ \int_0^t V_s(\frac{1}{2}\Delta\phi)  \,ds - K_t(\phi),
% a.e. t\geq0,
%&&
\nonumber
%&&\mbox{}
\end{eqnarray}
Note that although some of the convergences leading to the above
equation hold in $L^2[0,T]$,
all terms are right continuous in $t$ and so the equality
holds for all $t$, and not just for a.e. $t$. By equation~\eqref
{eq44}\vadjust{\goodbreak} and the fact that $U\in D^{\ep}([0,T], C_{\mathrm{rap}})$ [from (a)] we see that $K_\cdot(\phi)\in D^\ep([0,T],\R)$. It then
follows from $K\in M_F(E)$ that $K_\cdot\in D^{\ep}([0,T], M_F(\R))$,
and this proves the last part of (a).

Now we will rewrite the above equation in the mild form. The derivation
is a
bit more complicated than, for example, (\ref{eqmild1}) for $\bU$,
due to the presence of the measure-valued term $K$.
%Now it is standard to rewrite the above equation in the mild form (see
%e.g. (\ref{eqmild1}) for $\bU$):
For any $\phi\in C^+_b(\R)$, $t\in[0,T]\setminus\CG_{\ep}$,
\begin{eqnarray*}
U_t(\phi) &=& \sum_{s_i\in\CG_{\ep}^{\mathrm{odd}}, s_i\leq t} \int
_{\R} S_{t-s_i}\phi(y) J^{x_i}_{\ep}(y)
\,dy
\\
&&{} +\int_0^t\int
_{\R} S_{t-s}\phi(y) U(s,y)^{\gamma}
W^{U}(ds,dy)
\\
&&{}-\int_0^t\int
_{\R} S_{t-s}\phi(y)K(ds,dy).
\end{eqnarray*}
Writing $S_t$ in terms of $p_t$, we have
\begin{eqnarray}
\label{eq231112}
\nonumber
U_t(\phi) &=& \int_{\R}
\phi(x)\sum_{s_i\in\CG_{\ep
}^{\mathrm{odd}}, s_i\leq t} \int_{\R}
p_{t-s_i}(y-x) J^{x_i}_{\ep}(y) \,dy \,dx
\\
&&{} +\int_{\R} \phi(x)\int_0^t
\int_{\R} p_{t-s}(x-y) U(s,y)^{\gamma}
W^{U}(ds,dy) \,dx
\\
\nonumber
&&{}-\int_{\R}\phi(x)\int
_0^t\int_{\R}
p_{t-s}(x-y)K(ds,dy) \,dx, \qquad P\mbox{-a.s.},
\nonumber
\end{eqnarray}
where the last equality follows by the Fubini and the stochastic Fubini
theorems. Note that we take the time $t$ outside the set $\CG_{\ep}$
since, for $t\in\CG_{\ep}$, $K(\{t\}, dx)$ could be strictly positive,
and with $p_0$ being a delta measure, this creates difficulties with
applying the Fubini theorem. Therefore the case of $t\in\CG_{\ep}$ will
be treated separately.

By {(a)}, we know that
%
%e9.19 #&#
\begin{eqnarray}
\label{eq23111} U\in D^{\ep}\bigl([0,T], C^+_{\mathrm{rap}}\bigr),
\qquad P\mbox{-a.s.}
\end{eqnarray}
By the domination
\[
U^{\gamma}\leq\bU^{\gamma} \in D^\ep\bigl([0,T],
C^+_{\mathrm{rap}}\bigr),
\]
Lemma~\ref{pmom}, and Lemma~\ref{lem3061}(b) we may choose a version
of the stochastic integral so that
%
%e9.20 #&#
\begin{eqnarray}
\label{eq23112}\qquad t \mapsto\int_0^t\int
_{\R} p_{t-s}(\cdot-y) U(s,y)^{\gamma}
W^{U}(ds,dy)\in C\bigl([0,T], C_{\mathrm{rap}}\bigr),
\nonumber
\\[-8pt]
\\[-8pt]
\eqntext{P\mbox{-a.s.},}
\end{eqnarray}
and in what follows we will always consider such a version. This, and
the fact that $K_\cdot\in D^{\ep}([0,T], M_F(\R))$, implies that the
equality in (\ref{eq231112}) holds
$P$-a.s. \textit{for all} $t\in[0,T]\setminus\CG_{\ep}$, and, hence,
we get
%(\ref{eq231112}) implies
%
%e9.21 #&#
\begin{eqnarray}
\label{eqmild2}
U_t(x) &=& \sum
_{s_i\in\CG_{\ep}^{\mathrm{odd}}, s_i\leq t} \int_{\R} p_{t-s_i}(x-y)
J^{x_i}_{\ep}(y) \,dy
\nonumber\\
&&{}+\int_0^t\int_{\R}
p_{t-s}(x-y) U(s,y)^{\gamma} W^{U}(ds,dy)
\nonumber
\\[-8pt]
\\[-8pt]
\nonumber
&&{}
-\int_0^t\int_{\R}
p_{t-s}(x-y)K(ds,dy),
\\
\eqntext{\mbox{Leb-a.e. } x \in\R, \mbox{for each $t\in\bigl([0,T]\setminus
\CG_{\ep}\bigr)$ }, P\mbox{-a.s.}}
\end{eqnarray}
Now let us check that the above equation holds for all $ (t,x) \in
([0,T]\setminus\CG_{\ep})\times\R$, $P$-a.s.
[recall again that Lemma~\ref{lem3061}(b) is used to select an
appropriate jointly continuous version of the stochastic integral].
First, note that the steps similar to those leading to~\eqref{eqmild2}
easily imply
\begin{eqnarray}
\label{eqmild3}
\nonumber
U_t(x) &=& S_{t-r}U_r(x)+
\sum_{s_i\in\CG_{\ep}^{\mathrm{odd}},
r<s_i\leq
t} \int_{\R}
p_{t-s_i}(x-y) J^{x_i}_{\ep}(y) \,dy
\\
&&{}+\int_r^t\int_{\R}
p_{t-s}(x-y) U(s,y)^{\gamma} W^{U}(ds,dy)
\nonumber
\\[-8pt]
\\[-8pt]
\nonumber
&&{}-\int_r^t\int
_{\R} p_{t-s}(x-y)K(ds,dy),
\\
\eqntext{\mbox{Leb-a.e. } x \in\R, \mbox{for all $r,t\in [0,T]\setminus
\CG_{\ep}, r\le t$, } P\mbox{-a.s.}}
\end{eqnarray}
%
%To this end, we first establish the continuity of the mapping
%(t,x)\mapsto\int_0^t\int_{\R} p_{t-s}(x-y)K(ds,dy)
% on $([0,T]\setminus\CG_{\ep})\times\R$.
%Moreover,
Lemma~\ref{lem3061}(b) could be easily strengthened to assure, that,
in fact, the process
%
%e9.22 #&#
\begin{eqnarray}\label{eq23114}
X(r,t,x) \equiv \int_r^t\int
_{\R} p_{t-s}(x-y) U(s,y)^{\gamma}
W^{U}(ds,dy),
\nonumber
\\[-8pt]
\\[-8pt]
\eqntext{0\leq r\leq t\leq T, x\in\R, \mbox{is $P$-a.s. continuous in $(r,t,x)$}}
\end{eqnarray}
and
%
%e9.23 #&#
\begin{eqnarray}
\label{eq23115} X(t,t,\cdot)=0\qquad \forall t\in[0,T].
\end{eqnarray}
Again, to be more precise, there exists just a version of the process
$X$ such that~\eqref{eq23114} holds, and, in what follows, we will
always consider such a version.

As was already noted following Lemma~\ref{lemmp-43},
%
%e9.24 #&#
\begin{eqnarray}
\label{eq23113} t \mapsto\sum_{s_i\in\CG_{\ep}^{\mathrm{odd}}, s_i\leq t} \int
_{\R} p_{t-s_i}(\cdot-y) J^{x_i}_{\ep}(y)
\,dy \in D^{\ep}\bigl([0,T], C^+_{\mathrm{rap}}\bigr),\qquad  P\mbox{-a.s.}\hspace*{-35pt}
\end{eqnarray}
Let us take $A\subset\Omega$ such that $P(A)=1$ and for each $\omega
\in A$, (\ref{eq23111}) and (\ref{eqmild2})--(\ref{eq23113}) hold.
Fix an arbitrary $\omega\in A$ and $(t,x)\in((0,T]\setminus\CG_{\ep
})\times\R$. Then choose
$\{(r_l,z_k)\}_{l,k\geq1}$ such that
the equality in (\ref{eqmild3}) holds with $(r_l,t,z_k)$ in place of
$(r,t,x)$, and
$(r_l,z_k)\rightarrow(t,x)\in([0,T]\setminus\CG_{\ep})\times\R$, as
$l,k\rightarrow\infty$. Also assume that
$r_l<t$, for all $l\geq1$. Note that both $\{(r_l,z_k)\}_{l,k\geq1}
,(t,x)$ may depend on $\omega$.
We would like to show
%
%e9.25 #&#
\begin{eqnarray}
\label{eq231113}&& \lim_{k\rightarrow\infty} \int_{0}^{t}
\int_{\R} p_{t-s}(z_k-y)K(ds,dy)
\nonumber
\\[-8pt]
\\[-8pt]
\nonumber
&&\qquad =
\int_{0}^{t}\int_{\R}
p_{t-s}(x-y)K(ds,dy).
\end{eqnarray}

Fix $\delta>0$.
By (\ref{eq23111}), (\ref{eq23114}) and (\ref{eq23115}) we can
choose $l^*$ sufficiently large so that,
with $r^*\equiv r_{l^*}$, we have
%
%e9.26 #&#
\begin{eqnarray}
\label{eq23118}&& \bigl\llvert U_{t}(z_k) -
S_{t-r^*}U_{r^*}(z_k) \bigr\rrvert
\nonumber
\\[-8pt]
\\[-8pt]
\nonumber
&&\qquad{}+ \biggl\llvert
\int_{r^*}^{t}\int_{\R}
p_{t-s}(z_k-y) U(s,y)^{\gamma} W^{U}(ds,dy)
\biggr\rrvert \leq\delta
\end{eqnarray}
for all $k\geq1$.
Note that we assume without loss of generality that
%$r^*< t_k$, for all $k\geq1$, and
%$\eta$ is sufficiently small and $t_1, t_2,\ldots,$ are close to $t$,
%such that,
\[
\bigl[r^*,t\bigr] %\cup\Bigl(\cup_{k=1}^\infty[r^*,t_k]\Bigl)
\subset[0,T]\setminus\CG_{\ep}.
\]

Now we are ready to show~(\ref{eq231113}).
First, by the bounded
convergence theorem and $K_\cdot\in D^\ep([0,T],M_F(\R))$, we get
%
%e9.27 #&#
\begin{eqnarray}
\label{eq231110} \int_0^{r^*}\int
_{\R} p_{t-s}(z_k-y)K(ds,dy)
\rightarrow\int_0^{r^*}\int_{\R}
p_{t-s}(x-y)K(ds,dy)
\end{eqnarray}
as $k\rightarrow\infty$.
Next consider \eqref{eqmild3} with $r=r^*$, $x=z_k$,
%apply $S_\eta$ to both sides and evaluate the result at $z_k$. Compare
%the result to \eqref{eqmild2} with $(t,x)=(t_k,z_k)$
to conclude that
\begin{eqnarray}
U_{t}(z_k) &=& S_{t-r^*}U_{r^*}(z_k)
\nonumber\\
&&{} +\int_{r^*}^{t}\int
_{\R} p_{t-s}(z_k-y)
U(s,y)^{\gamma} W^{U}(ds,dy)
\\
\nonumber
&&{}-\int_{r^*}^{t}\int
_{\R} p_{t-s}(z_k-y)K(ds,dy) \qquad\forall
k\geq1.
\nonumber
\end{eqnarray}
Therefore,
%
%e9.28 #&#
\begin{eqnarray}
\label{eq23117}
\nonumber\quad
&&\int_{r^*}^{t}\int
_{\R} p_{t-s}(z_k-y)K(ds,dy) \\
&&\qquad\leq
\bigl\llvert U_{t}(z_k) - S_{t-r^*}U_{r^*}(z_k)
\bigr\rrvert
+ \biggl\llvert \int_{r^*}^{t}\int
_{\R} p_{t-s}(z_k-y)
U(s,y)^{\gamma} W^{U}(ds,dy)\biggr\rrvert
\\
\nonumber
&&\qquad\leq \delta\qquad \forall k\geq1,
\end{eqnarray}
where the last bound follows from~(\ref{eq23118}).
This together with Fatou's lemma and $K\in D^\ep([0,T],M_F(\R))$ implies
%
%e9.29 #&#
\begin{eqnarray}
\label{eq23119}&&\int_{r^*}^{t}\int
_{\R} p_{t-s}(x-y)K(ds,dy)
\nonumber
\\[-8pt]
\\[-8pt]
\nonumber
&&\qquad\leq \liminf_{k\to\infty} \int_{r^*}^{t}
\int_{\R} p_{t-s}(z_k-y)K(ds,dy)\leq
\delta.
\end{eqnarray}
Equations~(\ref{eq23117}), (\ref{eq23119}) and (\ref{eq231110}) imply
\[
\limsup_{k\rightarrow\infty} \biggl\llvert \int_{0}^{t}
\int_{\R} p_{t-s}(x-y)K(ds,dy)- \int
_{0}^{t}\int_{\R}
p_{t-s}(z_k-y)K(ds,dy)\biggr\rrvert \leq3\delta,
\]
and since $\delta$ was arbitrary, (\ref{eq231113}) follows.

Equation~(\ref{eq231113}) together with (\ref{eq23111}), (\ref
{eq23112}), (\ref{eq23113})
implies that the equality in~(\ref{eqmild2}) holds for \textit{all}
$(t,x)\in([0,T]\setminus\CG_{\ep})
\times\R$ on a set of full probability measure. Moreover, since all
the other terms in~(\ref{eqmild2}) except
$\int_0^t\int_{\R} p_{t-s}(\cdot-y)K(ds,dy)$ are in $ D^{\ep}([0,T],
C^+_{\mathrm{rap}})$, we get that, in fact,
\[
t\mapsto\int_0^t\int_{\R}
p_{t-s}(\cdot-y)K(ds,dy)\in C\bigl([0,T]\setminus \CG_{\ep},
C^+_{\mathrm{rap}}\bigr),\qquad P\mbox{-a.s.}
\]

Now let $t\in\CG_{\ep}$, and let us show that, at $t$, the
$C^+_{\mathrm{rap}}$-valued mapping $r\mapsto\int_0^r\int_{\R}
p_{r-s}(\cdot
-y)K(ds,dy)$ is right continuous and with a left limit.
We will prove it for $t=s_j\in\CG^{\mathrm{odd}}_{\ep}$
for some $j$ (for $t\in\CG^{\mathrm{even}}_{\ep}$ the argument is the
same, even simpler).
Note that the measure $K(\{s_j\},dx)$ is absolutely continuous with
respect to
Lebesgue measure. This follows from \eqref{eq44} and the fact that
$U$ is in $D^\ep([0,T],C^+_{\mathrm{rap}})$.
We will denote the density of $K(\{s_j\},dx)$ by $K(\{s_j\},x), x\in\R
$. Take $\eta>0$ sufficiently small such that
$(s_j,s_j+\eta]\subset[0,T]\setminus\CG_{\ep}$. Then, since~(\ref
{eqmild2}) holds for all $(t,x)\in([0,T]\setminus\CG_{\ep})
\times\R$, we get
%
%e9.30 #&#
\begin{eqnarray}
\label{eq231111}
\nonumber
U_{s_j+\eta}(x) &=& \sum
_{s_i\in\CG_{\ep}^{\mathrm{odd}}, s_i< s_j} \int_{\R
} p_{s_j+\eta-s_i}(x-y)
J^{x_i}_{\ep}(y) \,dy
\\
&&{}+\int_{\R} p_{\eta}(x-y)
J^{x_j}_{\ep}(y) \,dy
\nonumber\\
&&{} +\int_0^{s_j+\eta}\int
_{\R} p_{s_j+\eta-s}(x-y) U(s,y)^{\gamma}
W^{U}(ds,dy)
\\
\nonumber
&&{}-\int_0^{s_j+\eta}\int
_{\R} p_{s_j+\eta
-s}(x-y) \bigl(K(ds,dy)-\delta
_{s_j}(ds)K\bigl(\{s_j\},dy\bigr)\bigr)
\\
&&{}-\int_{\R} p_{\eta}(x-y)K\bigl(
\{s_j\},y\bigr) \,dy \qquad\forall x\in\R.
\nonumber
\end{eqnarray}
Take $\eta\downarrow0$. Since the measure $(K(ds,dy)-\delta
_{s_j}(ds)K(\{s_j\},dy))$ gives zero mass to the set $\{s_j\}\times\R$,
by the
argument similar to the one used in the case of $t\in[0,T]\setminus
\CG
_{\ep}$,
we can easily derive that
\begin{eqnarray*}
&&\int_0^{s_j+\eta}\int_{\R}
p_{s_j+\eta-s}(\cdot -y) \bigl(K(ds,dy)-\delta _{s_j}(ds)K\bigl(
\{s_j\},dy\bigr)\bigr)
\\
&&\qquad\rightarrow \int_0^{s_j}\int
_{\R} p_{s_j-s}(\cdot-y) \bigl(K(ds,dy)-\delta
_{s_j}(ds)K\bigl(\{ s_j\},dy\bigr)\bigr),
\end{eqnarray*}
in $C_{\mathrm{rap}}$, as $\eta\downarrow0$. Moreover,
$U_{s_j+\eta}(\cdot)$ and the first three terms on the right-hand side
of~(\ref{eq231111}) converge in $C_{\mathrm{rap}}$. This immediately
implies that the last term
$\int_{\R} p_{\eta}(\cdot-y)K(\{s_j\},y) \,dy$ also converges in
$C_{\mathrm{rap}}$, and clearly the limit is
%
%e9.31 #&#
\begin{eqnarray}
\label{eq24112} K\bigl(\{s_j\},\cdot\bigr)\in C_{\mathrm{rap}},
\end{eqnarray}
or more precisely a $C_{\mathrm{rap}}$-valued version of this density.
%Hence, all the terms in~(\ref{eq231111}) converge in $C_{\mathrm{rap}}$,
%the last one also converges %in $C_{\mathrm{rap}}$ to $K(\{t\},x),  x\in
All together we get that~(\ref{eqmild2}) holds also for $t\in\CG
^{\mathrm{odd}} _{\ep}$ with $p_0$
being the Dirac measure; moreover the $C_{\mathrm{rap}}$-valued mapping
$r\mapsto\int_{0}^{r}\int_{\R} p_{r-s}(\cdot-y)K(ds,dy)$
is right
continuous at $t\in\CG^{\mathrm{odd}} _{\ep}$. The existence of left-hand
limits for $r\mapsto\int_{0}^{r}\int_{\R} p_{r-s}(\cdot
-y)K(ds,dy)$ at
$t\in\CG^{\mathrm{odd}} _{\ep}$ follows by a similar argument. As we noted
above, the same proof works for $t\in\CG^{\mathrm{even}} _{\ep}$,
and this finishes the proof of {(b)}.

%pa9.subsection.subsubsection.7 #&#
(c)
%ver on any interval $[s_n,t_{n+1})$ the
% the left-hand side is continuous and hence the right-hand side is
%continuous and hence $K$ does not have
By the above $t\mapsto K_t$ is continuous on $[0,T]\setminus\CG_{\ep}$.
Since $\{K^n\}$ is a sequence of continuous, nondecreasing
measure-valued processes, its tightness in $M_F(E)$ immediately
implies tightness on all the open intervals between the jumps of the
limiting process,
in the space of continuous measure-valued paths,
that is,
in $C([0,T]\setminus\CG_{\ep}, M_F(\R))$.

So, the only jumps $K$ may possibly have are at the points $s_i, t_i\in
\CG_{\ep}$. We recall that a jump of measure-valued process $K$ at any
$t\in[0,T]$ equals $K(\{t\},dx)=K(\{t\},x) \,dx$, where by \eqref
{eq24112} $K(\{t\},\cdot)\in C_{\mathrm{rap}}$ for all $t\in\CG_\ep$.
%%recall that as above, we can consider $K$ as a measure on $E$.
We now calculate the sizes of those jumps. Consider the possible jump
at $s_i$.
Assume $\phi$ is a nonnegative function in $C^2_{c}(\R)$.
By~\eqref{eq44} (and it's analogue for $V$), $U=w^+$ and $V=w^-$,
we have the following conditions on $w^{\pm}_{s_i}$:
%
%e9.32 #&#
%e9.33 #&#
\begin{eqnarray}
\label{eq411} {\bolds\Delta}\bigl\langle w^+,\phi\bigr\rangle(s_i)
&=& \bigl\langle J^{x_i},\phi \bigr\rangle- \bigl\langle K\bigl(
\{s_i\},\cdot\bigr), \phi\bigr\rangle,
\\
\label{eq412} {\bolds\Delta}\bigl\langle w^-,\phi\bigr\rangle(s_i)
&=& - \bigl\langle K\bigl(\{s_i\},\cdot\bigr), \phi\bigr\rangle
\leq0. %\\
%w^+_{s_i}(x)w^-_{s_i}(x)&=&0, w^+_{s_i-}(x)w^-_{s_i-}(x)=0,
\end{eqnarray}
The above are preserved under bounded pointwise limits in $\phi$ and
so continue
to hold for any bounded Borel $\phi\ge0$.

We consider two cases. First assume $\phi$ is such that
\[
\operatorname{supp}(\phi)\subset % \{x: w^+_{s_i-}(x)>0\}\subset
\bigl\{x\dvtx w^-_{s_i-}(x)=0\bigr\}.
\]
%
%$x$ be such that $w^-_{s_i-}(x)=0, w^+_{s_i-}(x)>0.$
Then ${\bolds\Delta}\langle w^-,\phi\rangle(s_i)=\langle
w^-_{s_i},\phi\rangle\ge0$
and so (\ref{eq412}) immediately implies that
$\langle K(\{s_i\},\cdot), \phi\rangle=0$.

Now let $\phi$ be such that
\[
\operatorname{supp}(\phi)\subset % \{x: w^+_{s_i-}(x)>0\}\subset
\bigl\{x\dvtx w^+_{s_i-}(x)=0\bigr
\}.
\]
%
%$x$ be such that $w^-_{s_i-}(x)=0, w^+_{s_i-}(x)>0.$
Then ${\bolds\Delta}\langle w^+_{s_i},\phi\rangle=\langle
w^+_{s_i},\phi\rangle\geq0$
and so (\ref{eq411}) immediately implies that\break
$\langle K(\{s_i\}, \cdot), \phi\rangle\leq\langle J^{x_i},\phi
\rangle$.

We may write $1=\phi_1+\phi_2$, where $\phi_i$ is as in case $i$
($i=1,2$) [because $w^+_{s_i-}(x)w^-_{s_i-}(x)\equiv0$]. It therefore
follows that
\[
{\bolds\Delta}\langle K_{s_i},1\rangle= \bigl\langle K\bigl(
\{s_i\},\cdot\bigr), 1\bigr\rangle\leq\bigl\langle J^{x_i},1
\bigr\rangle=\ep,
\]
and we are done.
\end{pf}

%le9.8 #&#
\begin{lemma}
\label{lem152}
The following assertions hold.
\begin{longlist}[(a)]
\item[(a)] For any $i\in\NN_{\ep}$,
$\{U^{i,n}\}_{n\geq1}$, $\{\tU^{i,n}\}_{n\geq1} $,
$\{V^{i,n}\}_{n\geq1},\{\tV^{i,n}\}_{n\geq1}$ are
tight in $C([0,T]\setminus\CG_{\ep}, M_F(\R))$.
% and $\{K^{i,n,U}\}_{n\ge1}$ and $\{K^{i,n,V}\}_{n\ge1}$ are tight
%in $M_F(E)$.
%
\item[(b)] For any $i,j\in\NN_{\ep}$, and $\phi_l\in
C_{b}(\R),
l=1,\ldots, 4$,
\[
\bigl\{\bigl(M^{i,n,U}(\phi_1),M_t^{j,n,V}(
\phi_2),\tM_t^{i,n,U}(\phi_3),\tM
_t^{j,n,V}(\phi_4)\bigr)\bigr\}_{n\geq1}
\]
is
tight in $C([0,T],\R)^4$.
\end{longlist}
\end{lemma}
%
%pa9.subsection.subsubsection.8 #&#
\begin{pf}
Fix an arbitrary $i\in\NN_{\ep}$. Let us first prove the tightness for
$\{U^{i,n}\}_{n\geq1}$. By the nonnegativity of $U^{i,n}$'s and the
domination $U^{i,n}\leq\bU^n\to\bU\in
D([0,T],  C^+_{\mathrm{rap}})$ a.s. [recall \eqref{asconv}], by Jakubowski's
theorem [see, e.g., Theorem II.4.1 in \citet{per02}], it is enough to prove
tightness of\break  $\{U^{i,n}(\phi)\}_{n\geq1}$ in $C([0,T]\setminus\CG
_{\ep
}, \R)$, for any $\phi\in C^2_{b}(\R)$.
From~(\ref{tUVndefn}) we get
%
%e9.34 #&#
\begin{eqnarray}
\label{eq29061} U^{i,n}_t(\phi) &=& \bigl\langle
J^{x_i},\phi\bigr\rangle\1(t\geq s_i) %\\
%&&
+ M_t^{i,n,U}(\phi)
\nonumber
\\[-8pt]
\\[-8pt]
\nonumber
&& \mbox{}+ \int_0^t
U^{i,n}_s(\Delta\phi/2) \,ds -K_t^{i,n,U}(
\phi ),\qquad  t\in[0,T].
\end{eqnarray}
For any $p>2$, we use H\"older's inequality to bound the $p$th moment
of the increment of the
third term on the right-hand side of~(\ref{eq29061}),
%
%e9.35 #&#
\begin{eqnarray}
&&E \biggl[ \biggl\llvert \int_u^t
U^{i,n}_s\biggl(\frac{1}{2}\Delta\phi\biggr) \,ds
\biggr\rrvert ^p \biggr]\nonumber
\\
&&\qquad\leq \sup_{s\leq T, x\in\R}e^{\lambda p |x|}E \bigl[ \bU
^n(s,x)^{p} \bigr] \biggl[ \int_{\R}
e^{-\lambda|x|} \biggl\llvert \frac{1}{2}\Delta\phi(x) \biggr\rrvert \,dx
\biggr]^{p} (t-u)^p,
\\
\eqntext{\forall0\leq u\leq t.}
\end{eqnarray}
Now use Lemma~\ref{qmom}(b) and the Kolmogorov tightness criterion
to see that
%
%e9.36 #&#
\begin{eqnarray}\qquad
\label{29062} \biggl\{\int_0^{\cdot}
U^{i,n}_s\biggl(\frac{1}{2}\Delta\phi\biggr) \,ds
\biggr\} _{n\geq1}\qquad \mbox{ is tight in } C\bigl([0,T],\R\bigr), \forall\phi
\in C^2_{b}(\R).
\end{eqnarray}
As for the martingale $M^{i,n,U}_{\cdot}(\phi)$, we can argue exactly
as in the proof of tightness for
$\{M^{n,U}(\phi)\}_{n\geq1}$ in Lemma~\ref{lem151}(a), by using
again the domination,
%use Burkholder-Davis-Gundy inequality, and
% again the domination
$U^{i,n}(s,\cdot)\le U^n(s,\cdot)\leq\bU^n(s,\cdot), s\in[0,T]$, to
show that
% to get, that for $p>2$,
% \begin{eqnarray}
%E\left[ \left|M^{i,n,U}_{t}(\phi)-M^{i,n,U}_{u}(\phi)\right|^p\right]}\
%&\leq& C_p
% \sup_{s\leq2, x\in\R}E\left[ \bU(s,x)^{p\gamma}\right] \left\langle
%1, \left|\phi\right|^2
%&&
%   \forall0\leq u\leq t.
%This, again, by Lemma~\ref{pmom} and Kolmogorov tightness criterion
%implies that
%
%e9.37 #&#
\begin{eqnarray}
\label{eq0271} \bigl\{M^{i,n,U}_{\cdot}(\phi) \bigr
\}_{n\geq1}\qquad \mbox{is tight in } C\bigl([0,T],\R\bigr)
\end{eqnarray}
for any $\phi\in C_{b}(\R)$. As for $K^{i,n,U}$, it is dominated from
the above by $K^{n}$ and
by Lemma~\ref{lem151}(c), $\{K^n\}_{n\geq1}$
is tight in $C([0,T]\setminus\CG_{\ep}, M_F(\R))$. Therefore $\{
K^{i,n,U}\}_{n\geq1}$ is also tight in the same space.

We combine this with (\ref{29062}), (\ref{eq0271}) and (\ref
{eq29061}) to finish
the proof of tightness of $\{U^{i,n}\}_{n\geq1}$ in $C([0,T]\setminus
\CG_{\ep}, M_F(\R))$.\vspace*{1pt}
% Note that we do not have tightness
%in $C([0,2],M_F(\R)$) \,due to the jump immigration term (see the first
%term on the right-hand side of~(\ref{eq29061})), however note that
%the process is identically zero on the
%interval $[0,s_i)$.

As for $\{\tU^{i,n}\}_{n\geq1}$, we get by the same argument as above that
%
%e9.38 #&#
\begin{eqnarray}\qquad
\label{29063} \biggl\{\int_0^{\cdot}
\tU^{i,n}_s(\Delta\phi/2) \,ds \biggr\} _{n\geq1}\qquad
\mbox{is tight in } C\bigl([0,T],\R\bigr), \forall\phi \in C^2_{b}(
\R).
\end{eqnarray}
For the martingale term, fix an arbitrary $\phi\in C_{b}$. We have
again tightness of $ \{\tM^{i,n,U}_{\cdot}(\phi) \}
_{n\geq1}$
in $ C([0,T],\R)$
by the same method as for $ \{M^{i,n,U}_{\cdot}(\phi) \}
_{n\geq1} $, by using the domination,
\[
\bigl[ \bigl(\tU^n(s,\cdot)+U^n(s,\cdot)
\bigr)^{2\gamma} - U^n(s,\cdot )^{2\gamma}
\bigr]^{1/2} \sqrt{\frac{\tU^{i,n}(s,\cdot)}{\tU^n(s,\cdot)}}\leq\bU ^n(s,\cdot
)^{\gamma},\qquad  s\in[0,T].
\]

%Burkholder-Davis-Gundy inequality, and
%similar domination,
% \begin{eqnarray}
%E\left[ \left|\tM^{i,n,U}_{t}(\phi)-\tM^{i,n,U}_{u}(\phi)\right|^p
%&\leq& C_p
% \sup_{s\leq2, x\in\R}E\left[ \bU(s,x)^{p\gamma}\right] \left\langle
%1, \left|\phi\right|^2
%&&
%   \forall0\leq u\leq t,
%and again the result follows.

The tightness of $\{V^{j,n}(\phi)\}_{n\geq1}$ and $\{\tV^{j,n}(\phi
)\}
_{n\geq1}$ follows in exactly the same way.
\end{pf}

In what follows we take
any converging subsequence of the processes from
Lemmas~\ref{lem152}(a), \ref{lem151}(a) and Corollary~\ref
{cor151}. Recall that ${\cal D}$ is the countable subset of
$C^2_{b}(\R)$ which is bounded-pointwise dense in $C_b(\R)$. By
Lemma~\ref{lem152}(b) we can take a further subsequence, if needed,
so that
all the martingales from Lemma~\ref{lem152}(b) indexed by functions
from ${\cal D}$ converge in $C([0,T],\R)$.

To simplify notation we will still index
this subsequence by $n$.
Let us also switch to the Skorohod space where all the processes
mentioned in the previous paragraph
converge a.s. Since
$(\bU^{n}, \bV^{n})$ has the same law as the weakly unique in $D^{\ep
}([0,T], C^+_{\mathrm{rap}})^2$ solution to~(\ref{eq28})
[by Theorem~1.1 of~\citet{myt98w}],
we may, and shall, assume that on our probability space
$(\bU^{n}, \bV^{n})\to(\bU, \bV)$ in $D([0,T],C^+_{\mathrm{rap}})^2$, a.s.,
and, of course,
%
%e9.39 #&#
\begin{eqnarray}
\label{eq1106} U^{i,n}, \tU^{i,n}, U^{n},
\tU^n&\leq& \bU^n\qquad \forall n\geq1, i\in
\NN_{\ep},
\nonumber
\\[-8pt]
\\[-8pt]
\nonumber
V^{i,n}, \tV^{i,n}, V^{n},
\tV^n&\leq& \bV^n\qquad \forall n\geq1, i\in
\NN_{\ep}.
\end{eqnarray}
For $i\in\NN_{\ep}$, let
\[
U, V, \tU, \tV, \bU,\bV, K, U^{i}, V^{i},
\tU^{i}, \tV^{i}, K^{i,U}, K^{i,V}
\]
be the limiting points of $\{U^{n}\}_{n\geq1}, \{V^{n}\}_{n\geq1}$,
$\{\tU^{n}\}_{n\geq1}$, $ \{\tV^{n}\}_{n\geq1}$, $\{\bU^{n}\}
_{n\geq
1}$, $ \{\bV^{n}\}_{n\geq1}$,
$\{K^n\}_{n\ge1}$, $ \{U^{i,n}\}_{n\geq1}$, $\{V^{i,n}\}_{n\geq1}$,
$\{\tU^{i,n}\}_{n\geq1}$,
$\{\tV^{i,n}\}_{n\geq1}$, $\{K^{i,n,U}\}_{n\geq1}$, $\{K^{i,n,V}\}
_{n\geq1}$, respectively.
Clearly w.p.1 for all $t\in[0,T]\setminus\CG_\ep$,
%
%e9.40 #&#
%e9.41 #&#
\begin{eqnarray}
\label{Usum} U_t &=& \sum_{i \in\N_{\ep}}
U^{i}_t,\qquad\tU_t = \sum
_{i \in\N
_{\ep}} \tU ^{i}_t,
\\
V_t &=& \sum_{i \in\N_{\ep}}V^{i}_t
,\qquad\tV_t = \sum_{i \in\N
_{\ep}} \tV
^{i}_t,
\end{eqnarray}
by the corresponding equations for the approximating processes,
\[
\bU_t=U_t+\tU_t, \qquad\bV_t=V_t+
\tV_t \qquad\mbox{for all }t\in[0,T]
\]
by the same reasoning and Lemma~\ref{lem151}(a), and
\[
K=\sum_{i\in\N_\ep} K^{i,U}=\sum
_{j\in\N_\ep} K^{j,V}.
\]
%
%Now we can define
%U^{i}&=& U^{i,*}- K^i\\
%V^{i}&=& V^{i,*}- K^i\\
%Obviously the processes above belong to $D([0,2], M_F(\R))$.
By Lemma~\ref{lem151}(a) we may take versions of $U,\tU,V,\tV,\bU
,\bV
$ in $D^\ep([0,T],C^+_{\mathrm{rap}})$. We next refine the state space of
the subprocesses corresponding to the individual clusters.
%
%le9.9 #&#
\begin{lemma}
\label{lem0571}
For any $i\in\NN_{\ep}$,
\begin{eqnarray*}
&&\bigl(U^{i},\tU^i,V^i,\tV^{i},K^{i,U},K^{i,V}
\bigr)\\
&&\qquad\in \bigl(D^{\ep}\bigl([0,T],M_F(\R)\bigr)\cap
L^2(E)\bigr)^4\times D\bigl([0,T],M_F(\R)
\bigr)^2
\end{eqnarray*}
and $ (U^{i},\tU^i,V^i,\tV^{i},K^{i,U},K^{i,V} )_{i\in\N
_{\ep
}}$ satisfy (\ref{UVdefn}), (\ref{eq22}) and (\ref{tUVdefn}).
\end{lemma}
%
%pa9.subsection.subsubsection.9 #&#
\begin{pf}
% Recall that by...the solution
%is weakly unique hence the sequence of processes $\{(\bU^{n}, \bV^{n})
% weakly converging in $D([0,2],C_{\mathrm{rap}})^2$.
%Now the processes $(U^{n}, V^{n})$ are bounded from above by $(
%$(U^{n}(s,x)dxds, V^{n}(s,x)dxds)$ are
%tight in $(M_F(E)\times M_F(E))$.
Although $U^i$ (and similarly $V^i, \tU^i, \tV^i$) is defined
%in~(\ref{28061}) via the difference of limiting points of
%$U^{i,n}+K^{i,n}$ and $K^{i,n}$, we will see that in fact it can be
%also defined as a limit of
as a limit point of $\{U^{i,n}\}_{n\geq1}$ in $C([0,T]\setminus\CG
_{\ep}, M_F(\R))$, it can be also considered
as a limit of $\{U^{i,n}\}_{n\geq1}$ in the weak $L^2(E)$ topology [in
the sequel we denote the space $L^2(E)$ equipped\vadjust{\goodbreak}
with the weak topology, by $L^{2,w}(E)$].
Indeed, since
by~(\ref{eq1106}), all $U^{i,n}, \tU^{i,n}$ (resp., $V^{i,n}, \tV
^{i,n}$) are bounded from above
by $\bU^n \to\bU$ in $D([0,T], C^+_{\mathrm{rap}})$ [resp., $\bV^n\to\bV
\mbox
{ in } D([0,T],C^+_{\mathrm{rap}})$], we get that, in fact,
\[
\bigl\{U^{i,n}\bigr\}_{n\geq1}, \bigl\{\tU^{i,n}\bigr
\}_{n\geq1}, \bigl\{V^{i,n}\bigr\}_{n\geq1}, \bigl\{
\tV^{i,n}\bigr\}_{n\geq1}
\]
are all relatively compact in $L^{2,w}(E)$. This and the convergence of
$\{U^{i,n}\}_{n\geq1}$, $\{V^{i,n}\}_{n\geq1}$,
$\{\tU^{i,n}\}_{n\geq1}$,
$\{\tV^{i,n}\}_{n\geq1}$, in $C([0,T]\setminus\CG_\ep,M_F(\R))$ as
$n\rightarrow\infty$,
imply that
\[
\bigl(U^{i,n},\tU^{i,n}, V^{i,n},
\tV^{i,n}\bigr) \rightarrow\bigl(U^i,\tU^i,
V^i,\tV^i\bigr)\qquad \mbox{in } L^{2,w}(E)^4,
P\mbox{-a.s.}, \mbox{as } n\rightarrow \infty.
\]
%
%Similarly convergence in $L^{2,w}(E)$ holds for the processes $
%In what follows we will use the version of the limiting process
%defined by limit points in
%Lemma~\ref{lem152}(a),
Therefore we have
\[
U^i,\tU^i, V^i,\tV^i\in C
\bigl([0,T]\setminus\CG_{\ep}, M_F(\R)\bigr)\cap
L^2(E).
\]
From our earlier remark prior to Proposition~\ref{prop4} and
$K^{i,U},K^{u,V}\in M_F(E)$, we have
\[
\bigl(K^{i,U}, K^{i,V}\bigr)\in D\bigl([0,T],M_F(
\R)\bigr)^2.
\]

Now let us derive the semimartingale decomposition for $U^i$. Consider
the convergence of the
right-hand side of the equation for $U^{i,n}(\phi)$ in~(\ref
{tUVndefn}). By convergence of $\{U^{i,n}\}_{n\geq1}$ in $L^{2,w}(E)$
and in $C([0,T]\setminus\CG_\ep,M_F(\R))$ we
get that, for any $\phi\in C^2_{b}(\R)$ and any $t\le T$,
%
%e9.42 #&#
\begin{eqnarray}
\int_0^t \int_{\R}
U^{i,n}_s(x)\frac{\Delta}{2}\phi(x) \,dx\,ds \rightarrow\int
_0^t \int_{\R}
U^{i}_s(x)\frac{\Delta}{2}\phi(x) \,dx\,ds
\nonumber
\\[-8pt]
\\[-8pt]
\eqntext{\mbox{as } n
\rightarrow\infty.}
\end{eqnarray}
Now fix an arbitrary $\phi\in{\cal D}$. By Lemma~\ref{lem152}(b) we
may assume that $M^{i,n,U}(\phi)$
converges a.s. in $C([0,T],\R)$.
Moreover, using a bound analogous to~(\ref{eq0172}),
one can immediately get that, for any $p\geq2$,
the martingale $M_t^{i,n,U}(\phi)$ is bounded in $L^p(dP)$ uniformly in
$n$ and $t\in[0,T]$.
Hence, the limiting process is a continuous
$L^2$-martingale that we will call $M^{i,U}(\phi)$. For its quadratic
variation,
%note that as $U^{i,n}, U^{n}\leq\bU\in C_{\mathrm{rap}}$, then the
recall that the sequence $\{ (U^n)^{2\gamma-1}\}_{n\geq1}$ converges
to $U^{2\gamma-1}$ \textit{strongly} in $L^2(E)$ [by \eqref{genpc}]
and this together with convergence of $\{U^{i,n}\}_{n\geq1}$ in
$L^{2,w}(E)$ implies that,
for any $\phi\in C_{b}(\R)$ and $t\le T$, w.p.1
%
%e9.43 #&#
\begin{eqnarray}
\label{eq112}
\nonumber
\bigl\langle M^{i,n,U}(\phi)\bigr
\rangle_t &=& \int_0^t \int
_{\R} U^n(s,x)^{2\gamma-1}
U^{i,n}(s,x) \phi(x)^2 \,dx\,ds
\nonumber
\\[-8pt]
\\[-8pt]
\nonumber
&\rightarrow& \int_0^t \int
_{\R}U(s,x)^{2\gamma-1} U^{i}(s,x) \phi
(x)^2 \,dx\,ds\qquad\mbox{as } n\rightarrow\infty.
\end{eqnarray}
Hence, again by boundedness of $M_t^{i,n,U}(\phi)$, in $L^p(dP), p\geq
2$, uniformly in $t\in[0,T], n\geq1$,
we get that the limiting continuous
martingale $M^{i,U}$ has quadratic variation
\[
\bigl\langle M^{i,U}(\phi)\bigr\rangle_t=\int
_0^t \int_{\R}U(s,x)^{2\gamma-1}
U^{i}(s,x) \phi(x)^2 \,dx\,ds
\]
for all $\phi\in{\cal D}\subset C_{b}(\R)$.
%So we get, that $U^i\in D([0,2], M_F(\R))\cap L^2(E)$.
Moreover, by repeating the above argument for $V^{i,n}$ we get that
$(U^{i}, V^{i})_{i\in\NN_{\ep}}$, solves the following martingale problem:
%
%e9.44 #&#
\begin{eqnarray}\qquad\quad
\label{eq111} \cases{ %
\mbox{For
all $\phi_i, \psi_j\in{\cal D}\subset
C^2_{b}(\R)$,}
\vspace*{2pt}\cr
\displaystyle U^{i}_t(\phi_i) = \bigl\langle
J^{x_i},\phi_i\bigr\rangle\1(t\geq s_i)
%&&
+ M_t^{i,U}(\phi_i)
\vspace*{2pt}\cr
\hspace*{42pt}\displaystyle{}+ \int_0^t
U^{i}_s\biggl(\frac{1}{2}\Delta\phi_i
\biggr) \,ds - K^{i,U}_t(\phi_i)\qquad \forall t\in[0,T]
, i\in\NN_{\ep},
\vspace*{2pt}\cr
\displaystyle V^{j}_t(\psi_j) = \bigl\langle
J^{y_j},\psi_j\bigr\rangle\1(t\geq t_j)
% \\
%&&
+ M_t^{j,V}(\psi_j)
\vspace*{2pt}\cr
\hspace*{46pt}\displaystyle{}+ \int_0^t V^{j}_s
\biggl(\frac{1}{2}\Delta\psi_j\biggr) \,ds - K^{j,V}_t(
\psi_j)\qquad \forall t\in[0,T], j\in\NN_{\ep}, }
\end{eqnarray}
where $M^{i,U}(\phi_i), M^{j,V}(\psi_j)$ are martingales such that for
all $i,j\in\NN$,
%
%e9.45 #&#
\begin{eqnarray}\label{eq117}\qquad\quad
\cases{ %
\displaystyle\bigl\langle
M^{i,U}_{\cdot}(\phi_i), M_{\cdot}^{j,U}(
\phi_j)\bigr\rangle_t = \delta_{i,j} \int
_0^t \int_{\R}
U(s,x)^{2\gamma-1} U^{i}(s,x) \phi_i(x)^2
\,dx\,ds,
\vspace*{2pt}\cr
\displaystyle\bigl\langle M^{i,V}_{\cdot}(\psi_i),
M_{\cdot}^{j,V}(\psi_j)\bigr\rangle_t
= \delta_{i,j} \int_0^t \int
_{\R} V(s,x)^{2\gamma-1} V^{i}(s,x)
\psi_i(x)^2 \,dx\,ds,
\vspace*{2pt}\cr
\displaystyle\bigl\langle M_{\cdot}^{i,U}(\phi_i),
M_{\cdot}^{j,V}(\psi_j)\bigr\rangle_t
=0\qquad\forall i,j\in\NN_{\ep}.}
\end{eqnarray}
Note that the equality in (\ref{eq111})
holds for any $t$ in $[0,T]\setminus\CG_{\ep}$ since both left- and
right-hand sides are continuous processes on
$[0,T]\setminus\CG_{\ep}$; moreover the right-hand side is cadlag on
$[0,T]$. Using this and the domination $U^i_t\le\bar U_t$ and
$V^i_t\le\bar V_t$ for $t\notin\CG_\ep$, we may construct versions of
$U^i$ and $V^i$ in $D^{\ep}([0,T], M_F(\R))\cap L^2(E)$ so that
equality in (\ref{eq111})
holds for all $t$ in $[0,T]$. Clearly the martingale problem~(\ref{eq111})
can be also extended to all $\phi_i,\psi_j\in C^2_{b}(\R)$ by a
limiting procedure, again using the
$L^p(dP)$ boundedness of the martingales for any $p\geq2$.

Now let us handle the processes $(\tU^i, \tV^i), i\in\N_{\ep}$. By the
same steps that were used to treat $(U^i, V^i)_{i\in\N_{\ep}}$
we get that $(\tU^i, \tV^i)_{i\in\N_{\ep}}$ satisfies the following
martingale problem:
%
%e9.46 #&#
\begin{eqnarray}\qquad\quad
\label{eq0373} \cases{ %
\mbox{For
all $\phi_i, \psi_j\in{\cal D}\subset
C^2_{b}(\R)$,}
\vspace*{2pt}\cr
\displaystyle\tU^{i}_t(\phi_i) = \bigl\langle
J^{x_i},\phi_i\bigr\rangle\1(t\geq s_i)
%&&
+ \tM_t^{i,U}(\phi_i)
\vspace*{2pt}\cr
\hspace*{42pt}\displaystyle{} + \int_0^t
\tU^{i}_s\biggl(\frac{1}{2}\Delta\phi_i
\biggr) \,ds + K^{i,U}_t(\phi_i)\qquad \forall t\in[0,T]
, i\in\NN_{\ep},
\vspace*{2pt}\cr
\displaystyle\tV^{j}_t(\psi_j) = \bigl\langle
J^{y_j},\psi_j\bigr\rangle\1(t\geq t_j)
% \\
%&&
+\tM_t^{j,V}(\psi_j)
\vspace*{2pt}\cr
\hspace*{46pt}\displaystyle{} + \int_0^t \tV^{j}_s
\biggl(\frac{1}{2}\Delta\psi_j\biggr) \,ds + K^{j,V}_t(
\psi_j)\qquad \forall t\in[0,T], j\in\NN_{\ep},} %\label{eq111}
\end{eqnarray}
where by Lemma~\ref{lem152} $\tM^{i,U}(\phi_i), \tM^{j,V}(\psi_j)$
are continuous processes. By the same argument as
before [the uniform in $n$ and $t$, boundedness $L^p(dP), p\geq2$, of
the approximating martingales] they are
martingales and we would like to show that, for any $i,j\in\NN_{\ep}$,
%
%e9.47 #&#
\begin{eqnarray}
\label{eq0372} \cases{ %
\displaystyle\bigl\langle
\tM^{i,U}_{\cdot}(\phi_i), \tM_{\cdot}^{j,U}(
\phi _j)\bigr\rangle_t \vspace*{2pt}\cr
\displaystyle\qquad= \delta_{i,j} \int
_0^t\int_{\R}
\frac{ (\tU(s,x)+U(s,x) )^{2\gamma}-
U(s,x)^{2\gamma}}{\tU(s,x)}
\vspace*{2pt}\cr
\hspace*{75pt}\displaystyle{}\times\tU^i(s,x) \phi_i(x)^2 \,dx \,ds,
\vspace*{2pt}\cr
\displaystyle\bigl\langle\tM^{i,V}_{\cdot}(\psi_i),
\tM_{\cdot}^{j,V}(\psi _j)\bigr\rangle_t\vspace*{2pt}\cr
\displaystyle\qquad = \delta_{i,j} \int_0^t\int
_{\R}\frac{ (\tV(s,x)+V(s,x) )^{2\gamma}-
V(s,x)^{2\gamma}}{\tV(s,x)}
\vspace*{2pt}\cr
\hspace*{75pt}\displaystyle{}\times\tV^i(s,x) \psi_i(x)^2 \,dx \,ds,
\vspace*{2pt}\cr
\displaystyle\bigl\langle\tM_{\cdot}^{i,U}(\phi_i),
\tM_{\cdot}^{j,V}(\psi _j)\bigr\rangle_t
=0. }
\end{eqnarray}
As before, the orthogonality of the limiting martingales follows easily
by the uniform in $n$ and $t$, $L^p(dP),
p\geq2$, boundedness of
the approximating martingales and their orthogonality. Next we
calculate the quadratic variations. We will do it
just for $\tM^{i,U}_{\cdot}(\phi)$, for some $i\in\NN_{\ep}$.
It is enough to show that for any
$\phi\in C_{b}(\R)$ and $t\in[0,T]$,
%
%e9.48 #&#
\begin{eqnarray}
\label{eq0371}\quad
&& \int_0^t\int
_{\R}\frac{ (\tU^n(s,x)+U^n(s,x) )^{2\gamma}-
U^n(s,x)^{2\gamma}}{\tU^n(s,x)}\tU^{i,n}(s,x) \phi(x)\,dx \,ds
\nonumber
\\[-8pt]
\\[-8pt]
\nonumber
&&\qquad\rightarrow\int_0^t\int
_{\R}\frac{ (\tU(s,x)+U(s,x)
)^{2\gamma}- U(s,x)^{2\gamma}}{\tU(s,x)}\tU^i(s,x) \phi(x)\,dx \,ds,
\end{eqnarray}
in $L^{1}(dP)$, $\mbox{as }
n\rightarrow\infty$. Denote
\[
F(\tilde u, u) \equiv (\tilde u+ u )^{2\gamma}- u^{2\gamma}.
\]
Then, for any
$\phi\in C_{b}(\R)$ and $t\in[0,T]$, we get
%
%e9.49 #&#
\begin{eqnarray}
&& \biggl| \int_0^t\int
_{\R}\frac{F(\tU^n(s,x),U^n(s,x))}{\tU
^n(s,x)}\tU^{i,n}(s,x) \phi(x)\,dx \,ds
\nonumber\\
&&\hspace*{5pt}\mbox{}-\int_0^t\int
_{\R} \frac{F(\tU(s,x),U(s,x))}{\tU
(s,x)}\tU^i(s,x) \phi(x)\,dx \,ds
\biggr|
\nonumber\\
&&\qquad\leq \biggl| \int_0^t\int
_{\R} \biggl(\frac{F(\tU^n(s,x),U^n(s,x))}{\tU^n(s,x)}- \frac{F(\tU(s,x),U(s,x))}{\tU(s,x)} \biggr)
\nonumber
\\[-8pt]
\\[-8pt]
\nonumber
&&\hspace*{141pt}\qquad\quad{}\times\tU^{i,n}(s,x) \phi(x)\,dx \,ds \biggr|
\\
\nonumber
&& \qquad\quad\mbox{}+\biggl\llvert \int_0^t
\int_{\R}\frac{F(\tU(s,x),U(s,x))}{\tU
(s,x)} \bigl( \tU^i(s,x)-
\tU^{i,n}(s,x) \bigr) \phi(x)\,dx \,ds\biggr\rrvert
\\
\nonumber
&&\qquad\equiv I^{1,n}+I^{2,n}.
\end{eqnarray}
Clearly
%
%e9.50 #&#
\begin{eqnarray}
\label{eq29063} \frac{F(\tU(s,x),U(s,x))}{\tU(s,x)}\leq2\gamma\bU^{2\gamma-1}\in
L^2(E),
\end{eqnarray}
and hence by convergence
of $\tU^{i,n}$ to $\tU^i$ in $L^{2,w}(E)$, a.s., we get that
\[
I^{2,n}\rightarrow0\qquad \mbox{as } n\rightarrow\infty  \mbox{ a.s.}
\]
and by dominated convergence it is easy to get that, in fact, the
convergence is in $L^1(dP)$. As
for $I^{1,n}$, by using $\llvert \frac{ \tU^{i,n}(s,x)}{ \tU
^n(s,x)}\rrvert \leq1$ we immediately get that
\begin{eqnarray*}
I^{1,n}&\leq& \int_0^t
\int_{\R}\bigl\llvert F\bigl(\tU^n(s,x),U^n(s,x)
\bigr)- F\bigl(\tU(s,x),U(s,x)\bigr)\bigr\rrvert \phi(x)\,dx \,ds
\\
&&{}+ \int_0^t\int
_{\R}\frac{F(\tU(s,x),U(s,x))}{\tU
(s,x)}\bigl\llvert \tU (s,x)-
\tU^n(s,x)\bigr\rrvert \phi(x)\,dx \,ds.
\end{eqnarray*}
We again use (\ref{eq29063}) and convergence of $\tU^n$ and $U^n$ to
$\tU$ and $U$, respectively, in $L^p(E)$ for any
$p\geq1$,
we immediately get that, $I^{1,n}\rightarrow0$, a.s., as $n\rightarrow
\infty$. Use again the dominated
convergence theorem to get that, in fact, the convergence holds in
$L^1(dP)$, and~(\ref{eq0371}) follows.
As a result we get that $(U^i,V^i,\tU^i,\tV^i), i\in\NN_{\ep}$ solves
the martingale problem
(\ref{eq111}), (\ref{eq117}), (\ref{eq0373}), (\ref{eq0372}),
with all martingales corresponding to
different processes being orthogonal.

Now, as before, [see the proof of Lemma~\ref{lem151}(b)], the
martingales in the martingale problem can be represented as stochastic
integrals with respect to independent white
noises, and hence one immediately gets that $(U^i,V^i,\tU^i,\tV
^i)_{i\in
\NN_{\ep}}$ solves~(\ref{UVdefn}), \eqref{eq22} and (\ref{tUVdefn})
but with
$(U^i,V^i,\tU^i,\tV^i)\in(D^{\ep}([0,T],M_F(\R))\cap L^2(E))^4,
i\in
\NN_{\ep}$. Here we note that equality in \eqref{Usum} as $M_F(\R
)$-valued processes extends to all $t\in[0,T]$ by right-continuity.
%The other properties of
% $(U^i,V^i,\tU^i,\tV^i)_{i\in\NN_{\ep}}$ required in (\ref{UVdefn}-
% $U^i(s,x)V^j(s,x)=0, i,j\in\NN_{\ep}$,
%and $U^i(s,\cdot), \tU^i(s,\cdot)=0, s<s_i$, $V^i(t,\cdot), \tV^i(t,
\end{pf}

To finish the proof of Proposition~\ref{prop4} we next verify the
following lemma.

%le9.10 #&#
\begin{lemma}
\label{lem0572}
$U^i, \tU^i, V^i, \tV^i\in C([0,T]\setminus\CG_{\ep},C^+_{\mathrm{rap}})
\cap D^{\ep}([0,T], L^1(\R)),\break
\forall i\in\NN_{\ep}$.
\end{lemma}
%
%pa9.subsection.subsubsection.10 #&#
\begin{pf}
We will prove it just for $U^i$, as the proof for the other terms goes
along exactly along the same lines.
Similarly to the steps in the proof of Lemma~\ref{lem151}(b), we
first write the equation for $U^i$ in the mild form to get
\begin{eqnarray}
\label{eq24111}
U^i(t,x)&=& \int_{\R}
p_{t-s_i}(x-y)J^{x_i}_{\ep}(y) \,dy
\nonumber\\
&&{}+\int_0^t\int_{\R}
p_{t-s}(x-y) U(s,y)^{\gamma-1/2}U^i(s,y)^{1/2}
W^{i,U}(ds,dy)
\nonumber
\\[-8pt]
\\[-8pt]
\nonumber
&&{}-\int_0^t\int
_{\R} p_{t-s}(x-y) K^{i,U}(ds,dy) \\
\eqntext{\mbox{Leb-a.e. } (t,x)\in\bigl([0,T]\setminus\CG_{\ep
}\bigr)\times\R.}
\end{eqnarray}
We now argue as in the proof of part (b) of Lemma~\ref{lem151}. The
first term on the right-hand side
of~(\ref{eq24111}) clearly belongs to
$D^{\ep}([0,T], C_{\mathrm{rap}})$. Similarly by the bound
\[
U^{\gamma-1/2}\bigl(U^i\bigr)^{1/2} \leq
\bU^{\gamma} \in D\bigl([0,T], C^+_{\mathrm{rap}}\bigr),
\]
Lemma~\ref{pmom}, and Lemma~\ref{lem3061}(b), we see that the second
term on the right-hand side is in
$C([0,T], C_{\mathrm{rap}})$. As for the third term on the right-hand side,
one can use the domination
$K^{i,U}\leq K$, Lemma~\ref{lem151}(b) to get that $K^{i,U}(\{t\},
dx)=0$ for any $t\in[0,T]\setminus\CG_{\ep}$.
For $P$-a.s. $\omega$, take arbitrary $(t,x)\in([0,T]\setminus\CG
_{\ep
})\times\R$ and $\{(t_k,z_k)\}_{k\geq1}$, such that $\lim_{k\rightarrow
\infty} (t_k,z_k)=(t,x)$. Then by Lemma~\ref{lem151}(b), we get that
$\{\1(s< t_k)p_{t_k-s}(z_k-y)\}$ is uniformly integrable with respect
to $K(ds,dy)$ and hence by domination it is also uniformly integrable
with respect to $K^{i,U}(ds,dy)$. This gives continuity of the mapping
\[
(r,x)\mapsto\int_0^{r}\int
_{\R} p_{r-s}(x-y)K^{i,U}(ds,dy)
\]
on $([0,T]\setminus\CG_{\ep})\times\R$, and again by domination we may
easily show that
\[
r\mapsto\int_0^{r}\int_{\R}
p_{r-s}(\cdot-y)K^{i,U}(ds,dy)\in C\bigl([0,T]\setminus
\CG_{\ep},C^+_{\mathrm{rap}}\bigr).
\]
%
%For example, if $x_n\to x$, then $\{1(s< t)p_{t-s}(x_n-y)\}$ is
%uniformly integrable with respect to $K(ds,dy)$ and hence also with
%respect to $K^{i,U}(ds,dy)$.
All together, this gives that the right-hand side of~(\ref
{eq24111}) belongs to
$C([0,T]\setminus\CG_{\ep}, C_{\mathrm{rap}})$. Hence there is a version of
$U^i$ which is in $C([0,T]\setminus\CG_{\ep}, C^+_{\mathrm{rap}})$ as well.

Note that, in fact, the above argument also easily implies that for any
$t\in\CG_{\ep}$,
%
%e9.51 #&#
\begin{eqnarray}
\label{eq24114} U^i(r,\cdot) \rightarrow U^i(t-, \cdot)\qquad
\mbox{in } C_{\mathrm{rap}}, P\mbox{-a.s.}
\end{eqnarray}
as $r\uparrow t$, where
\begin{eqnarray}
\label{eq24113}
\nonumber
U^i(t-,x)&=& \1(t>s_i)\int
_{\R} p_{t-s_i}(x-y)J^{x_i}_{\ep}(y)
\,dy
\\
&&{} +\int_0^t\int_{\R}
p_{t-s}(x-y) U(s,y)^{\gamma-1/2}U^i(s,y)^{1/2}
W^{i,U}(ds,dy)
\\
\nonumber
&&{}-\int_0^t\int
_{\R} p_{t-s}(x-y) \bigl(K^{i,U}(ds,dy)-
\delta_{t}(ds)K^{i,U}\bigl(\{t\},dy\bigr)\bigr)
\nonumber
\end{eqnarray}
for $x\in\R$.
Indeed, for $(t,x)\in\CG_{\ep}\times\R$, take again arbitrary
$(t_k,z_k)$ such that $t_k\uparrow t$ and $z_k\rightarrow x$, as
$k\rightarrow\infty$. Again by Lemma~\ref{lem151}(b), we get that
$\{\1(s< t_k)p_{t_k-s}(z_k-y)\}$ is uniformly integrable with respect
to $(K(ds,dy)-K(\{t\},dy))$; hence by domination it is also uniformly
integrable with respect to $(K^{i,U}(ds,dy)-\delta_{t}(ds)K^{i,U}(\{t\}
,dy))$. This easily implies that
$ U^i(t_k,z_k) \rightarrow U^i(t-,x)$, where $U^i(t-,x)$
satisfies~(\ref{eq24113}), and hence (\ref{eq24114})
follows.

Clearly, (\ref{eq24114}) implies that corresponding convergence also
holds in $L^1(\R)$, and hence
to finish the proof of the lemma it is enough to show that for any
$t\in\CG_{\ep}$,
%
%e9.52 #&#
\begin{eqnarray}
\label{eq24115} U^i(r,\cdot) \rightarrow U^i(t, \cdot)\qquad
\mbox{in } L^1(\R), P\mbox{-a.s.}
\end{eqnarray}
as $r\downarrow t$. Again, as in the proof of Lemma~\ref{lem151}(b),
we will show it for $t=s_j\in\CG^{\mathrm{odd}}_{\ep}$ for some $j$.
By~(\ref{eq111}), we get that
%
%e9.53 #&#
\begin{eqnarray}
\label{24116} U^i_{s_j}(dx)=U^i_{s_j-}(dx)+
\1(s_i=s_j)J^{x_i}_{\ep}(x) \,dx -
K^{i,U}\bigl(\{s_j\},dx\bigr).
\end{eqnarray}
Recall again that $K^{i,U}(\{s_j\},dx)$ is dominated by $K(\{s_j\}
,dx)$, which, in turn, by~(\ref{eq24112})
is absolutely continuous with a density function in $C^+_{\mathrm{rap}} $.
Therefore $K^{i,U}(\{s_j\},dx)$ is also
absolutely continuous with a density function $K^{i,U}(\{s_j\},x),  x\in
\R$, bounded by a function in
$C^+_{\mathrm{rap}}$. This together with~(\ref{eq24114}), our assumptions
on $J^{x_i}_{\ep}$ and~(\ref{24116})
implies that $U^i_{s_j}(dx)$ is absolutely continuous with bounded
density function
%
%e9.54 #&#
\begin{eqnarray}
\label{eq24117} U^i_{s_j}(\cdot)\in L^1(\R).
\end{eqnarray}
For any $\eta\in(0,\ep/2)$, by combining \eqref{24116}, \eqref
{eq24114} (with $t=s_j$) and \eqref{eq24111} (with $t=s_j+\eta$),
we have
\begin{eqnarray}
\label{eq24118}
\nonumber
&&U^i(s_j+\eta,\cdot)\\
&&\qquad=
S_{\eta}U^i(s_j,\cdot)
\nonumber
\\[-8pt]
\\[-8pt]
\nonumber
&&\qquad\quad{} +\int_{s_j}^{s_j+\eta}\int_{\R}
p_{s_j+\eta
-s}(\cdot-y) U(s,y)^{\gamma-1/2}U^i(s,y)^{1/2}
W^{i,U}(ds,dy)
\\
&&\qquad\quad{}-\int_{s_j}^{s_j+\eta}\int
_{\R} p_{s_j+\eta
-s}(\cdot-y)\bigl (K^{i,U}(ds,dy)-
\delta_{s_j}(ds)K^{i,U}\bigl(\{s_j\},dy\bigr)\bigr)\nonumber
\end{eqnarray}
for $x\in\R$.
As $\eta\downarrow0$, the convergence to zero in $C_{\mathrm{rap}}$ of the
second and the third terms on the right-hand side follows easily as
in the last part of the proof of Lemma~\ref{lem151}(b). By~(\ref
{eq24117}), the first term on the right-hand side of~(\ref
{eq24118}) converges to $U^i(s_j,\cdot)$ in $L^1(\R)$ and we are done.
\end{pf}

%pa9.subsection.subsubsection.11 #&#
\begin{pf*}{Proof of Proposition~\protect\ref{prop4}}
Except for property \eqref{eq23}, Proposition~\ref{prop4} follows
from Corollary~\ref{cor151}, and
Lemmas~\ref{lem152}(a),~\ref{lem0571},~\ref{lem0572}. For
\eqref
{eq23} we note that
\[
U^i(t,x)V^j(t,x)\le U(t,x)V(t,x)=w^+(t,x)w^-(t,x)
\equiv0.
\]
\upqed\end{pf*}

%pa9.subsection.subsubsection.12 #&#
\begin{pf*}{Proof of Proposition~\protect\ref{thm11}}
As we mentioned in Remark~\ref{rem04}, since $T>1$ can be chosen
arbitrary large, it is sufficient to prove the
theorem on the time interval $[0,T]$.

Clearly, by Proposition~\ref{prop4} and the
definition of $\bU^i=U^i+\tU^i, \bV^i=V^i+\tV^i$, we immediately get
that
\[
\bigl(\bU^i,\bV^i\bigr)\in \bigl( C\bigl([0,T]\setminus
\CG_{\ep}, C^+_{\mathrm{rap}}\bigr)\cap D^{\ep}\bigl([0,T],
L^1(\R)\bigr) \bigr)^2,\qquad i\in\NN_{\ep},
\]
and satisfies~(\ref{eq25}) and (\ref{eq28}).
We saw in Section~\ref{secsetup} that \eqref{0early} and its analogue
for $(U^i,V^j)$ follow from the other properties.
%Note that by construction $\bU^i(s,\cdot)=0, s<s_i,$ and $
Then, by repeating the argument in the proof of Lemma~\ref{lem0572}
and taking into account the absence of the terms
$K^{i,U}, K^{i,V}$ at the right-hand side of the equations for $\bU^i,
\bV^i$, we immediately get that, in fact,
$(\bU^i,\bV^i)\in
D^{\ep}([0,T], C^+_{\mathrm{rap}})^2, i\in\NN_{\ep}$,
and
$\bU^i_{s_i+\cdot}\in C([0,T-s_i],C^+_{\mathrm{rap}})$, $\bV^i_{t_i+\cdot
}\in C([t_i,T-t_i],C^+_{\mathrm{rap}}), i\in\NN_{\ep} $, and part~(a)
of the
theorem follows. Part (b) follows from Lemma~\ref{lem151}(c).
\end{pf*}
\end{appendix}

% imsref loaded by akundreckaite, 2014-01-17 13:39:28

%


%suskaldyti doi

\printaddresses

\end{document}